\documentclass[12pt,twoside]{article}
\usepackage[english]{babel}
\usepackage[latin1]{inputenc}
\usepackage{amsmath}
\usepackage{amssymb,amsfonts}
\usepackage{graphicx}
\usepackage{times,amssymb,amscd}
\numberwithin{equation}{section}

\DeclareMathOperator{\arccosh}{arccosh}

\newcommand{\bA}{\mathbf{A}}
\newcommand{\bC}{\mathbf{C}}
\newcommand{\bE}{\mathbf{E}}
\newcommand{\bG}{\mathbf{G}}
\newcommand{\bH}{\mathbf{H}}

\newcommand{\bL}{\mathbf{L}}

\newcommand{\bR}{\mathbf{R}}
\newcommand{\bS}{\mathbf{S}}
\newcommand{\bV}{\mathbf{V}}

\newcommand{\bI}{\mathbf{I}}

\newcommand{\be}{\mathbf{e}}

\newcommand{\br}{\mathbf{r}}
\newcommand{\bx}{\mathbf{x}}
\newcommand{\by}{\mathbf{y}}
\newcommand{\bz}{\mathbf{z}}
\newcommand{\bg}{\mathbf{g}}
\newcommand{\bT}{\mathbf{T}}
\newcommand{\bX}{\mathbf{X}}

\newcommand{\bu}{\mathbf{u}}
\newcommand{\bv}{\mathbf{v}}

\newcommand{\BV}{\boldsymbol{V}}

\newcommand{\Bb}{\boldsymbol{b}}

\newcommand{\Be}{\boldsymbol{e}}
\newcommand{\Bu}{\boldsymbol{u}}
\newcommand{\Bv}{\boldsymbol{v}}
\newcommand{\BT}{\boldsymbol{T}}
\newcommand{\Bi}{\boldsymbol{i}}
\newcommand{\By}{\boldsymbol{y}}

\newcommand{\cA}{\mathcal{A}}
\newcommand{\cD}{\mathcal{D}}

\newcommand{\cP}{\mathcal{P}}
\newcommand{\cS}{\mathcal{S}}

\newcommand{\cC}{\mathcal{C}}
\newcommand{\cH}{\mathcal{H}}
\newcommand{\cB}{\mathcal{B}}

\newcommand{\EUC}{\mathbf E^3}
\newcommand{\SPH}{\bS^3}
\newcommand{\HYP}{\bH^3}
\newcommand{\HYN}{\bH^n}
\newcommand{\SXR}{\bS^2\!\times\!\bR}
\newcommand{\HXR}{\bH^2\!\times\!\bR}
\newcommand{\SLR}{\widetilde{\bS\bL_2\bR}}
\newcommand{\NIL}{\mathbf{Nil}}
\newcommand{\SOL}{\mathbf{Sol}}

\def\*{\star}
\def\~{\tilde}
\def\`{\bar}
\def\:{\,:\,}

\def\<{\langle}
\def\>{\rangle}
\def\({\left(}
\def\){\right)}
\def\[{\left[}
\def\]{\right]}

\def\0{{0^\prime}}
\def\1{{1^\prime}}
\def\2{{2^\prime}}
\def\3{{3^\prime}}

\def\Isom{\mathop{\rm Isom}\nolimits}

\def\pol{(_\*)}
\def\spr{\<\;,\;\>}

\def\om{\omega}

\mathchardef\Ga="0100
\mathchardef\Pi="0105

\mathchardef\wu="0875
\mathchardef\wv="0876
\mathchardef\wV="0856
\mathchardef\wze"0930
\mathchardef\we="0865
\mathchardef\wff="0866
\mathchardef\wx="0878
\mathchardef\wy="0879
\mathchardef\wa="0861
\mathchardef\wb="0862
\mathchardef\wr="0872
\mathchardef\wal="080B
\mathchardef\wbt="080C
\mathchardef\wsi="081B

\def\bX{{\bf X}}
\def\bG{{\bf G}}
\def\bE{{\bf E}}
\def\bH{{\bf H}}
\def\bS{{\bf S}}
\def\bV{{\bf V}}
\def\bT{{\bf T}}
\def\bC{{\bf C}}

\def\bL{{\bf L}}
\def\bx{{\bf x}}
\def\by{{\bf y}}
\def\be{{\bf e}}
\def\bu{{\bf u}}
\def\bv{{\bf v}}

\def\br{{\bf r}}

\def\b#1{{\bf#1}}

\def\PS{{\mathcal P}{\mathcal S}^3}
\def\P{{\mathcal P}^3}
\def\A{{\mathcal A}^3}

\begin{document}
\pagestyle{myheadings}
\markboth{\centerline{Jen\H o Szirmai}}
{Classical notions and problems in Thurston geometries}
\title
{Classical Notions and Problems in Thurston Geometries
\footnote{
\newline
}}

\author{Jen\H o Szirmai \\
\normalsize Budapest University of Technology and \\
\normalsize Economics, Institute of Mathematics, \\
\normalsize Department of Geometry \\
\normalsize Mûegyetem rkp. 3., H-1111 Budapest, Hungary \\
\normalsize szirmai@math.bme.hu
\date{\normalsize{\today}}}
\maketitle
\begin{abstract}
Of the Thurston geometries, those with constant curvature geometries (Euclidean $\EUC$, hyperbolic $\HYP$, spherical $\SPH$)  have been extensively studied,  
but the other five geometries, $\HXR$, $\SXR$, $\NIL$, $\SLR$, $\SOL$ 
have been thoroughly studied only from a differential geometry and topological 
point of view. However, classical concepts highlighting the beauty and underlying 
structure of 
these -- such as geodesic curves and spheres, the lattices,  
the geodesic triangles and their surfaces, their interior 
sum of angles and similar statements to those 
known in constant curvature geometries can be formulated. These 
have not been the focus of attention. In this survey, we summarize our results
on this topic and pose additional open questions.
\end{abstract}
\newtheorem{Theorem}{Theorem}[section]
\newtheorem{conj}[Theorem]{Conjecture}
\newtheorem{corollary}[Theorem]{Corollary}
\newtheorem{lemma}[Theorem]{Lemma}
\newtheorem{exmple}[Theorem]{Example}
\newtheorem{definition}[Theorem]{Definition}
\newtheorem{rmrk}[Theorem]{Remark}
\newtheorem{proposition}[Theorem]{Proposition}
\newenvironment{remark}{\begin{rmrk}\normalfont}{\end{rmrk}}
\newenvironment{example}{\begin{exmple}\normalfont}{\end{exmple}}
\newenvironment{acknowledgement}{Acknowledgement}

\section{Introduction}
\label{section1}
W. Thurston's results are essential for understanding the geometric structure of our world, 
where the eight so-called Thurston geometries play the leading role.
The importance of these geometries is emphasized by the Thurston's famous theorem:

Let $(X;\bG)$ be a $3$-dimensional homogeneous geometry, where $X$ is a simply connected
Riemannian space with a maximal group $\bG$ of isometries, acting transitively on $X$ with
compact point stabilizers. $\bG$ is maximal means that no proper extension of $\bG$ can act on
the Riemannian space $X$ in the same way. We recall the
\begin{Theorem}[Thurston, \cite{M97, S, T}] Any $3$-dimensional 
homogeneous geometry $(X;\bG)$ above that ad
mits a compact quotient is equivalent (equivariant) to one of the geometries $(X;\bG) =
Isom X)$ where the space $X$ is one of $\EUC$, $\HYP$, $\SPH$, $\HXR$, $\SXR$, $\NIL$, $\SLR$ or $\SOL$. 
\end{Theorem}
Therefore, there are eight so-called Thurston geometries, that were described in 
\cite{S,T,KPa20}. 
Among them $\EUC$, $\SPH$ and $\HYP$ are the classical spaces of constant zero, positive and negative curvature, respectively.  
Further geometries $\SXR$, $\HXR$ denote the direct product geometries where 
$\mathbf{S}^2$ is the spherical and $\mathbf{H}^2$ is the hyperbolic base plane and 
the real line $\mathbf{R}$ is with usual metric. 
Then $\SLR$ and $\NIL$ are obtained as twisted products 
of $\mathbf{R}$ with $\mathbf{H}^2$ and $\mathbf{E}^2$, respectively; 
and finally $\SOL$ geometry is a 
twisted product of the Minkowski plane $\mathbf{M}^2$ as fibre, 
with $\mathbf{R}$ as base. 
In each of them there exists an infinitesimal (positive definite) Riemannian metric 
that is invariant under certain translations, guaranteeing homogeneity at every point.
These translations in general commute only in $\EUC$, but a discrete (discontinuous) translation group, taken as a 
lattice, can be defined with compact fundamental domain in analogy to the Euclidean case, but with some different properties. 
The additional symmetries can define crystallographic groups, giving nice tilings, packings, material structures, etc.

I mention here only the packing and covering problems. In addition to 
pure mathematical curiosity, the study of sphere packings and coverings 
{\it generalized Kepler problem} 
is important because
it is possible that under different conditions (e.g. strong magnetic field) 
materials cannot be realized in the usual Euclidean space
but in one of the other Thurston geometries.
The structures of substances formed under these conditions 
may differ from the Euclidean case and can follow, for example, the geometry of non-constant curvature spaces, and in these new geometries their atoms can be modelled by  
$\HXR$, $\SXR$, $\NIL$, $\SLR$ or $\SOL$ spheres. For example, 
in $\NIL$ geometry we can define lattices and corresponding lattice-like ball packings 
where we found
geodesic ball packings with kissing number 14 that is denser than the 
densest Euclidean case (see Section 7, \cite{MSz06}, \cite{Sz07}). 
(The density is $\approx 0.78085.$)

A unified approach to Thurston geometries enabling the investigations in this direction were 
made possible by the nice paper of E.~Molnár \cite{M97} where he showed that the Thurston geometries can be uniformly modeled in the projective $3$-space $\cP^3$ 
(or in the projective $3$-sphere $\cP\cS^3$). 
The projective spherical model is based on linear algebra over the real vector space $\bV^4$ (for points) and its dual $\BV_4$ (for planes), up to a positive real factor, 
so that the proper dimension is indeed three. A plane $\rightarrow$ 
point polarity or $\BV_4 \times \BV_4 \rightarrow \bR$ scalar product 
(by specified signature) induces the invariant 
metric in a unified way. In our work we
will use these projective models of Thurston geometries.

The constant curvature geometries $\EUC$, $\HYP$, $\SPH$ have been extensively studied from the point of view of elementary geometry, differential geometry 
and topology. In this article we focus on results obtained in the other five non-constant curvature Thurston geometries $\HXR$, $\SXR$, $\NIL$, $\SLR$, $\SOL$.
These spaces have been investigated from the perspective 
of differential geometry and topology but 
few results are stated in connection with their internal 
structure in the classical sense. A lot of elementary notions are missing, problems are not formulated, theorems are not proved.
Hence, we focus on non-constant curvature Thurston geometries in this survey and 
we emphasize some surprising facts in our paper. 

Now, we review the notions of distance, angle, sphere, geodesic triangle, 
their surfaces and congruences in each afore mentioned geometry. 
Then Apollonius surfaces, Dirichlet-Voronoi cells, 
and in a separate chapter we review the the concepts of sphere (ball) packings and their densities 
and the corresponding results so far. 

{\it Furthermore, we emphasize the results related to the projective models of the considered geometries. 
In our opinion, these models are suitable for the elementary examination and visualization of the above geometries.}
\begin{rmrk} 
There is another way to define distance using the concept 
of so-called {\it translational distance}. 
We introduced this concept in paper \cite{MoSzi10}, 
but in this survey we summarize the results related 
to the concept of geodesic distance. Note that translation distance and geodesic 
distance 
are the same in the Euclidean geometry $\EUC$, Bolyai-Lobachevsky hyperbolic geometry $\HYP$, spherical $\SXR$ and $\HXR$ spaces, 
but give different values in the $\NIL$, $\SOL$ and $\SLR$ geometries 
(see \cite{CsSz16, MSzV, MSzV14, Sz19, Sz12-1, Sz13-2, VSz19}. 
\end{rmrk}
As the reader will see the above results and their visualizations will open a new window towards other (geometric) worlds. 
\section{On the projective models of Thurston geometries}
We first summarize the key information about projective models of Thurston geometries (see \cite{M97, MSz, MSz06}). 

All the Thurston 3-geometries will be uniformly modelled in the 
{\it projective spherical space} $\cP\cS^3$ that can be embedded into 
the affine hence into {\it Euclidean 4-space}.
Our main tool will be a $4$-dimensional vector space $\bV^4$ over the real numbers $\bR$ with basis
$\{\be_0,\be_1,\be_2,\be_3\}$, (that is not assumed to be orthonormal).

$\bV^4$ is the embedding
real vector space with its affine image $\cA^4(O,\bV^3,\BV_3)$. Let 
$\{O;~\be_0,\be_1,\be_2, \be_3\}$, be a
{\it coordinate system} in the affine $4$-space $\cA^4=\mathbf{E}^4$ with origin $O$ and a
(not necessarily orthonormal) vector basis $\{\be_0,\be_1,\be_2,\be_3\}$ for $\bV^4$, 
where our affine model plane
$\cA^3=\bE^3\subset \cP^3= \cA^3 \cup (\Bi)$ is placed to the point 
$E_0(\be_0)$ with equation $x^0=1$. 
Here any non-zero vector $\bx=x^0\be_0+x^1\be_1+x^2\be_2+x^3\be_3=:x^i\be_i$ (the index sum convention of Einstein-Schouten will be used)
represents a point $X(\bx)$ of $\cA^3$, but also a point of the {\it projective sphere} $\cP\cS^3$ after having introduced the following {\it positive equivalence}.  
For non-zero vectors
\begin{equation}
\begin{gathered}
\bx~ \sim ~c~\bx ~ \text{with} ~ 0 < c  \in\bR ~ \text{represent the same point} ~ X = (\bx~ \sim ~c~\bx) \ \
 \text{of} \ \ \cP\cS^3;  \\
\bz~ \sim ~ 0~\be_0+z^1\be_1+z^2\be_2+z^3\be_3 ~ \text{will be an ideal point}~ ~ (\bz)~  \text{of} ~ \cP\cS^3.
\end{gathered} \tag{2.1}
\end{equation}
We write: $(\bz) \in (\Bi)$, where $(\Bi)$ is the ideal plane (sphere) to $\cA^3$, extending the affine space $\cA^3$ into
the projective sphere $\cP\cS^3$. Here $(\bz)$ and $(-\bz)$, and in general $(\bx)$  and $(-\bx)$, are opposite points of $\cP\cS^3$. 
Then {\it identification} of the opposite point pairs of $\cP\cS^3$ leads to the projective space $\cP^3$.
Thus the embedding $\cA^3=\bE^3\subset \cP^3 \subset \cP\cS^3$ can be formulated in the vector space $\bV^3$ in a unified way. 

The dual (form) space $\BV_4$ to $\bV^4$ is defined as the set of {\it real valued linear functionals} or forms on $\bV^4$.
That means that we pose the following requirements for any form $\Bu\in \BV_3$
\begin{equation}
\begin{gathered}
\Bu~:~ \bV^4 ~\ni~ \bx \mapsto \bx\Bu \in \bR ~ \text{with linearity} \\
(a\bx+b\by)\Bu=a(\bx \Bu)+b(\by \Bu) ~ \text{for any} ~ \bx,~\by \in \bV^4 \ \
 \text{and for any} \ \ a,b\in\bR.
\end{gathered} \tag{2.2}
\end{equation}
{\bf We emphasize our convention}. {\it The vector coefficients are written from the left, 
then linear forms act on vectors on the right (as
an easy associativity law, similar conventions will be used also later on)}.

This ``built in" linear structure allows us to define the addition $\Bu+\Bv$ of two linear forms $\Bu$, $\Bv$, and the multiplication $\Bu c$ of a linear form
$\Bu$ by a real factor $c$, both resulting in linear forms of $\BV_4$. Moreover, we can define for any basis $\{\be_i\}$ in $\bV^4$ the
{\it dual basis} $\{\Be^j\}$ in $\BV_4$ by the Kronecker symbol $\delta_i^j$:
\begin{equation}
\be_i \Be^j=\delta_i^j = \left\{ \begin{array}{ll}
1 & \text{if}~ i = j \\
0 & \text{if}~ i \ne j
\end{array} \right. ~ ~ ~ i,j=0,1,2,3. \tag{2.3}
\end{equation}
Furthermore, we see that the general linear form $\Bu:=\Be^0u_0+\Be^1u_1+\Be^2u_2:=\Be^ju_j$ takes on the vector
$\bx:=x^0 \be_0+x^1 \be_1+x^2 \be_2 +x^3 \be_3:=x^i \be_i$ the real value
\begin{equation}
(x^i \be_i)(\Be^j u_j) = x^i ( \be_i \Be^j )u_j = x^i \delta_i^j u_j = x^i u_j : = x^0 u_0+x^1 u_1 + x^2 u_2+x^3 u_3. \tag{2.4}
\end{equation}
Thus, a linear form $\Bu \in \BV_4$ describes a {\it 3-dimensional subspace} 
$u$, i.e. a vector hyperplane of $\bV^4$ through the origin.
Moreover, forms
\begin{equation}
\Bu \sim \Bu  k ~ \text{with} ~0<k\in \bR ~ \text{represent the same oriented hyperplane of} ~ \bV^4. \tag{2.5}
\end{equation}
The positive equivalence class of forms $(\Bu)$ determines an open {\it half-space} $(\Bu)^+$ of $\bV^4$, i.e. the vector classes $(\bx)$ for which
\begin{equation}
(\Bu)^+~:=~\{(\bx)~:~\bx~\Bu >0 \}. \tag{2.6}
\end{equation}
This also gives a corresponding half-sphere of $\cP\cS^3$, and a corresponding half-hyperplane of $\cA^3$.
Note that {\it this is not so for the projective plane} $\cP^2$ which is not orientable,
because the equivalence mapping $\bx \mapsto -\bx$ has negative determinant in $\bV^4$!

In order to advance on our main goal, we introduce a {\it bijective linear mapping} $\bT$ of $\bV^4$ onto itself, i.e.
\begin{equation}
\begin{gathered}
\bT~:~\bV^3 \ni \bx \mapsto \bx\bT=:\by\in \bV^3~ ~ \text{with requirements} \\
x^i \be_i \mapsto (x^i \be_i) \bT=x^i(\be_i \bT)= x^i t_i^j \be_j =: y_j\be_j, ~ \det(t_i^j) \ne 0.
\end{gathered} \tag{2.7}
\end{equation}
Assume that $\bT$ has the above matrix $(t_i^j)$ with respect to basis $\{\be_i\}$ of $\bV^3$ $(i,j=0,1,2)$.
Then $\bT$ defines a {\it projective point transformation} $\tau(\bT)$ of $\cP\cS^2$
onto itself, which {\it preserves all the incidences of subspaces} of $\bV^3$ and hence {\it incidences of points and lines of} $\cP\cS^2$,
respectively. The matrix $(t_i^j)$ and its positive multiples $(c t_i^j)=(t_i^j c)$ with $0 < c \in \bR$ (and only these mappings) define the
same point transformation $\tau(\bT \sim \bT c)$ of $\cP\cS^2$ by the above requirements.
As usual, we define the composition, or product, of transforms $\bT$ and $\mathbf{W}$ of vector space $\bV^3$ in this order (right action on $\bV^3)$
by
\begin{equation}
\bT \mathbf{W}~:~\bV^4 \ni \bx \rightarrow (\bx\bT)\mathbf{W}=\by\mathbf{W}=\bz =: \bx(\bT \mathbf{W})
\tag{2.8}
\end{equation}
with matrices $(t_i^j)$ and $(w_j^k)$ to basis $\{\be_i\}$ $(i,j,k=0,1,2,3)$ as follows by our index conventions:
\begin{equation}
\be_i(\bT \mathbf{W})=(\be_i\bT)\mathbf{W}=(t_i^j\be_j)\mathbf{W}=(t_i^j)(w_j^k\be_k)=(t_i^jw_j^k)\be_k,
\tag{2.9}
\end{equation}
with summation (from 0 to 3) for the occurring equal upper and lower indices. 

The {\it inverse matrix class} of above $(t_i^j)$ is denoted by $(T_j^k) \sim \frac{1}{c} T_j^k$, 
with $t_i^j T_j^k=\delta_i^k$ {\it induces the corresponding 
linear transform $\BT$  of the dual $\BV_4$
(i.e. for lines) onto itself, and its inverse} $\BT^{-1}$, i.e. 
\begin{equation}
\begin{gathered}
\BT~:~\bV^4 \ni \Bv \mapsto \BT\Bv=:\Bu \in \BV_4~ ~ \text{such that} \\
\by\Bv=(\bx\bT)\Bv=\bx(\BT\Bv)=\bx\Bu,~ ~ \text{specifically} \\
0=\bx\Bu=(\bx \bT)(\BT^{-1}\Bu)= \By \Bv,~ \text{so} ~ 0=\by \Bv=(\by \bT^{-1})(\BT \Bv)=\bx\Bu \\
X~{\rm I}~u \leftrightarrow Y := X \tau ~ {\rm I} ~ v := \tau u ~ \text{holds}
\end{gathered} \tag{2.10}
\end{equation}
for the $\tau$-images of points and lines, respectively. We can see that the induced action on the dual $\BV_3$ is a
{\it left action} and so is the induced action on the lines of $\cP\cS^2$.

Let us introduce the positive equivalence in $\bV^4$ for non-zero vectors
\begin{equation}
\begin{gathered}
\bx~ \sim ~c~\bx ~ \text{with} ~ 0 < c  \in\bR ~ \text{defines the same point} ~ X(\bx) ~ 
\text{of} ~ \cP\cS^3
\end{gathered} \tag{2.11}
\end{equation}
whose coordinates in $\bx=x^i\be_i$, with respect to basis $\{\be_i \}$ $(i=0,1,2,3)$, can be written in matrix form
\begin{equation}
\bx=(x^0,x^1,x^2,x^3)
\begin{pmatrix} \mathbf{e}_0 \\ \be^1 \\ \be^2 \\ \be_3 \end{pmatrix}.~ ~\text{A form}~ ~
\Bu=\begin{pmatrix}
\Be^0&\Be^1&\Be^2&\Be^3
\end{pmatrix}
\begin{pmatrix} u_0 \\ u_1 \\ u_2 \\ u_3 \end{pmatrix} \tag{2.12}
\end{equation}
in the dual space $\BV_4$, again up to positive equivalence, describes an oriented plane ($2$-sphere) of $\cP\cS^3$ with the dual basis
$\{\Be^j\}, ~ \be_i\Be^j=\delta_i^j$ (the Kronecker symbol) $(i,j=0,1,2,3)$. Equalities
\begin{equation}
(x^i \be_i)(\Be^j u_j) = x^i ( \be_i \Be^j )u_j = x^i \delta_i^j u_j = x^i u_j  \tag{2.13}
\end{equation}
express the incidence $X~{\rm I}~u$. Formula (2.13) describes the set of 
varying points $X(\bx)$ on the fixed plane $u(\Bu)$, and at the same time
the set of planes $u(\Bu)$ incident to the fixed point $X(\bx)$.
The projective transform $\tau(\bT,\BT^{-1})$ with inverse matrix pair $(t_i^j)$ to $\bT$ of $\bV^4$ and $(t_i^j)^{-1}~ \sim ~ (T_j^k)$
to $\BT^{-1}$ of $\BV_4$ -- with respect to the dual basis pair $\{\be_i\}$, $\{\Be^j\}$, as in formulas (2.10) --
can be described in matrix form. First for points it is:
\begin{equation}
(x^0,x^1,x^2,x^3)
\begin{pmatrix}
t_0^0&t_0^1&t_0^2&t_0^3 \\
t_1^0&t_1^1&t_1^2&t_1^3 \\
t_2^0&t_2^1&t_2^2&t_2^3 \\
t_3^0&t_3^1&t_3^2&t_3^3
\end{pmatrix}
\begin{pmatrix} \be_0 \\ \be_1 \\ \be_2 \\ \be_3\end{pmatrix} ~\sim~
(y^0,y^1,y^2,y^3)
\begin{pmatrix} \be_0 \\ \be_1 \\ \be_2 \\ \be_3\end{pmatrix} \tag{2.14}
\end{equation}
and for planes $\Be^iT_i^k u_k \sim \Be^iv_i$.
Here $\tau(\bT,\BT^{-1})$ preserves the incidence by $0=(\bx\Bu)=(\by\Bv)$.
These are, again up to positive equivalence, related to a coordinate simplex $E_0E_1E_2E_3$ with the unit point
$E(\be=\be_0+\be_1+\be_2+\be_3$ and to $e^0 e^1 e^2 e^3$  with the unit plane $e(\Be=\Be^0+\Be^1+\Be^2+\Be^3$, where
$e^i=(E_jE_kE_l)$ with $\{0,1,2,3\}=\{i,j,k,l\}$.
	
Spherical space geometry $\SPH$ is defined with the additional polarity 
$\Pi(_*)$ or scalar product $\langle ~,~ \rangle$ in
$\BV_4$, with positive diagonal matrix $(\pi^{ij})$ and we have $\Be_*^i=\be^i=\pi^{ij}\be_j$,
\begin{equation}
(\Be^0,\Be^1,\Be^2,\Be^3) \mathop{\longrightarrow}_{\displaystyle *}
\begin{pmatrix} \mathbf{e}^0 \\ \be^1 \\ \be^2 \\ \be^3 \end{pmatrix}=
\begin{pmatrix}
1&0&0&0 \\
0&1&0&0 \\
0&0&1&0 \\
0&0&0&1
\end{pmatrix}
\begin{pmatrix} \be_0 \\ \be_1 \\ \be_2 \\ \be_3 \end{pmatrix}. \tag{2.15}
\end{equation}
For Euclidean $\EUC$ geometry and hyperbolic $\HYP$ (Bolyai-Lobachevsky) 
geometry the corresponding $\pi^{ij}$ matrices
$(i,j\in\{0,1,2,3\})$ are respectively:
\begin{equation}
\EUC ~:~
\begin{pmatrix}
0&0&0&0 \\
0&1&0&0 \\
0&0&1&0 \\
0&0&0&1
\end{pmatrix},
\quad \HYP ~:~ \begin{pmatrix}
-1&0&0&0 \\
0&1&0&0 \\
0&0&1&0 \\
0&0&0&1
\end{pmatrix}. \tag{2.16}
\end{equation}
For the other Thurston geometries we include Table 1 from \cite{MSz}, where additional 
information for the transform groups are also indicated.

\newpage
\centerline{\bf Table 1}
\vspace{3mm}
\centerline{The eight Thurston geometries modelled in $\cP\cS^3$ by a polarity}
\centerline{or scalar product and its isometry group.}
%
%
{\footnotesize{
$$
\vbox{\halign{\vrule width0pt height10.5pt depth3.5pt
\vrule\vrule
\vrule\hfil~#~\hfil&
  \vrule\hfil~#~\hfil&
  \vrule\hfil~#~\hfil&
  \vrule\hfil~#~\hfil
\vrule\vrule\vrule\cr
\noalign{\hrule}
\noalign{\hrule}
\noalign{\hrule}
&
  Signature of&
  &
  The group $G=\Isom\bX$ as\cr
Space&
  polarity $\Pi \pol$ &
  Domain of proper points&
  a special collineation\cr
$\bX$&
  or scalar prod-&
  of $\bX$ in $\PS$ $(\bV^4(\bR),\,\BV_4)$&
  group of $\PS$\cr
&
  uct $\spr$ in $\BV_4$&
  &
  \cr
\noalign{\hrule}
\noalign{\hrule}
\noalign{\hrule}
$\bS^3$&
  $(+\,+\,+\,+)$&
  $\PS$&
  Coll $\PS$ preserving $\Pi\pol$\cr
\noalign{\hrule}
$\bH^3$&
  $(-\,+\,+\,+)$&
  $\{(\bx)\in\P\:\<\bx,\bx\><0\}$&
  Coll $\P$ preserving $\Pi\pol$\cr
\noalign{\hrule}
&
  $(-\,-\,+\,+)$&
  Universal covering of $\cH:=$&
  Coll $\PS$ preserving $\Pi\pol$\cr
$\SLR$&
  with skew&
  $:=\{[\bx]\in\PS\:\<\bx,\bx\><0\}$&
  and fibres with 4 parameters. \cr
&
  line fibering&
  by fibering transformations&
  \cr
\noalign{\hrule}
$\bE^3$&
  $(0\,+\,+\,+)$&
  $\A=\P\setminus\{\om^\infty\}$ where&
  Coll $\P$ preserving $\Pi\pol$,\cr
&
  &
  $\om^\infty:=(\Bb^0)$, \ $\Bb^0_\*=\b0$&
  generated by plane reflections\cr
\noalign{\hrule}
&
  $(0\,+\,+\,+)$&
  &
  $G$ is generated by plane reflec-\cr
$\SXR$&
  with $O$-line&
  $\A\setminus\{O\}$&
  tions and sphere inversions,\cr
&
  bundle&
  $O$ is a fixed origin&
  leaving invariant the $O$-\cr
&
  fibering&
  &
  concentric 2-spheres of $\Pi\pol$\cr
\noalign{\hrule}
&
  $(0\,-\,+\,+)$&
  &
  $G$ is generated by plane reflec-\cr
&
  with $O$-line&
  $\cC^+=\{X\in\A\:$&
  tions and hyperboloid inver-\cr
$\HXR$&
  bundle&
  $\<\overrightarrow{OX},\overrightarrow{OX}\><0$, half cone$\}$&
  sions, leaving invariant the\cr
&
  fibering&
  by fibering&
  $O$-concentric half-hyperboloids\cr
&
  &
  &
  in the half-cone $\cC^+$ by $\Pi\pol$\cr
\noalign{\hrule}
&
  $(0\,-\,+\,+)$&
  $\A=\cP^3\setminus \phi$&
  Coll. of $\cA^3$ preserving \cr
$\SOL$&
  and parallel &
  &
  $\Pi(_*)$ and the \cr
&
  plane fibering&
  &
  fibering with 3 parameters  \cr
&
  with an ideal plane $\phi$&
  &
    \cr
\noalign{\hrule}
&
 Null-polarity $\Pi(_\star)$&
  $\A=\cP^3\setminus \phi$&
  Coll. of $\cA^3$ preserving \cr
$\NIL$&
  with parallel&
  &
  $\Pi(_\star)$ with \cr
&
  line bundle fibering &
  &
  $4$ parameters \cr
&
  $F$ with its polar &
  &
  \cr
&
ideal plane $\phi$ &
&
\cr
\noalign{\hrule}
\noalign{\hrule}
\noalign{\hrule}
}}
$$}}
%
\section{$\SXR$ and $\HXR$ spaces}
Of the rich literature not directly related to the projective model of two considered geometries examined, only a few would now be highlighted. 
In the nice papers \cite{MT14, R13-1, R13} the authors investigated the topic of special surfaces of the above geometries among others the surfaces of constant main curvature 
and related to these their minimal surfaces (see the further references given in them). In \cite{N17} L.~N{\'e}meth studied Pascal pyramids in the space $\HXR$ space 
whose examination has led to nice combinatorial and number theory correlations.
\subsection{Geodesic curves in $\SXR$ geometry}
In this section we recall important notions and results from the papers \cite{M97, PSSz10, PSSz11-2, PSz12, Sz13-1, Sz11-1,  
Sz11-2, Sz20, Sz21, W06}. 

The well-known infinitesimal arc-length square at any point of $\SXR$ is the following
\begin{equation}
   \begin{gathered}
     (ds)^2=\frac{(dx)^2+(dy)^2+(dz)^2}{x^2+y^2+z^2}.
       \end{gathered} \tag{3.1}
     \end{equation}
We shall apply the usual geographical coordiantes $(\phi, \theta), ~ (-\pi < \phi \le \pi, ~ -\frac{\pi}{2}\le \theta \le \frac{\pi}{2})$ 
of the sphere with the fibre coordinate $t \in \bR$. We describe points in the above coordinate system in our model by the following equations: 
\begin{equation}
x^0=1, \ \ x^1=e^t \cos{\phi} \cos{\theta},  \ \ x^2=e^t \sin{\phi} \cos{\theta},  \ \ x^3=e^t \sin{\theta} \tag{3.2}.
\end{equation}
Then we have $x=\frac{x^1}{x^0}=x^1$, $y=\frac{x^2}{x^0}=x^2$, $z=\frac{x^3}{x^0}=x^3$, i.e. the usual Cartesian coordinates.
By \cite{M97} we obtain that in this parametrization the infinitesimal arc-length square 
at any point of $\SXR$ is the following
\begin{equation}
   \begin{gathered}
      (ds)^2=(dt)^2+(d\phi)^2 \cos^2 \theta +(d\theta)^2.
       \end{gathered} \tag{3.3}
     \end{equation}
The geodesic curves of $\SXR$ are generally defined as having locally minimal arc length between any two (near enough) points. 
The system of equations of the parametrized geodesic curves $\gamma(t(\tau),\phi(\tau),\theta(\tau))$ in our model can be determined by the 
general theory of Riemann geometry (see \cite{KN}, \cite{Sz11-2}).

Then by (3.1-2) we obtain the system of equations of a geodesic curve in our Euclidean model (see \cite{Sz11-1} and Fig.~1.):
\begin{equation}
  \begin{gathered}
   x(\tau)=e^{\tau \sin{v}} \cos{(\tau \cos{v})}, \\ 
   y(\tau)=e^{\tau \sin{v}} \sin{(\tau \cos{v})} \cos{u}, \\
   z(\tau)=e^{\tau \sin{v}} \sin{(\tau \cos{v})} \sin{u},\\
   -\pi < u \le \pi,\ \ -\frac{\pi}{2}\le v \le \frac{\pi}{2}. \tag{3.4}
  \end{gathered}
\end{equation}
\subsection{Geodesic curves in $\HXR$ geometry}
In this section we recall the important notions and results from the papers \cite{M97,  PSSz11-2, Sz12-5}.

The points of $\HXR$, form an open cone in projective space $\mathcal{P}^3$, as follows:
\begin{equation}
\HXR:=\big\{ X(\bx=x^i \be_i)\in \mathcal{P}^3: -(x^1)^2+(x^2)^2+(x^3)^2<0<x^0,~x^1 \big\}. \notag
\end{equation}
E. Moln\'ar \cite{M97} found the infinitesimal arc length square, at any point of $\HXR$ as follows
\begin{equation}
   \begin{gathered}
     (ds)^2=\frac{1}{(-x^2+y^2+z^2)^2}\cdot [(x)^2+(y)^2+(z)^2](dx)^2+ \\ + 2dxdy(-2xy)+2dxdz (-2xz)+ [(x)^2+(y)^2-(z)^2] (dy)^2+ \\ 
     +2dydz(2yz)+ [(x)^2-(y)^2+(z)^2](dz)^2.
       \end{gathered} \tag{3.5}
     \end{equation}
This simplifies in the cylindrical coordinates $(t, r, \alpha)$, $(r \ge 0, ~ -\pi < \alpha \le \pi)$ 
with fibre coordinate $t \in \bR$. Points in our model are then 
\begin{equation}
x^0=1, \ \ x^1=e^t \cosh{r},  \ \ x^2=e^t \sinh{r} \cos{\alpha},  \ \ x^3=e^t \sinh{r} \sin{\alpha}  \tag{3.6}.
\end{equation}
Then we have $x=\frac{x^1}{x^0}=x^1$, $y=\frac{x^2}{x^0}=x^2$, $z=\frac{x^3}{x^0}=x^3$, i.e. the usual Cartesian coordinates.
We obtain by \cite{M97} that in this parametrization the infinitesimal arc length square in (3.5)
at an arbitrary point of $\HXR$ is 
\begin{equation}
   \begin{gathered}
      (ds)^2=(dt)^2+(dr)^2 +\sinh^2{r}(d\alpha)^2.
       \end{gathered} \tag{3.7}
     \end{equation}
The geodesic curves of $\HXR$ are generally defined as having locally minimal arc 
length between any two (near enough) points. 
The systems of equations of the parametrized geodesic curves $\gamma(t(\tau),r(\tau),\alpha(\tau))$ in our model can be determined by the 
general theory of Riemann geometry (see \cite{Sz12-5}).

Then by (3.6-7) we obtain the system of equations for a geodesic curve in our model 
(see \cite{Sz12-5} and Fig.~2.):
\begin{equation}
  \begin{gathered}
   x(\tau)=e^{\tau \sin{v}} \cosh{(\tau \cos{v})}, \\ 
   y(\tau)=e^{\tau \sin{v}} \sinh{(\tau \cos{v})} \cos{u}, \\
   z(\tau)=e^{\tau \sin{v}} \sinh{(\tau \cos{v})} \sin{u},\\
   -\pi < u \le \pi,\ \ -\frac{\pi}{2}\le v \le \frac{\pi}{2}. \tag{3.8}
  \end{gathered}
\end{equation}
\subsection{Distances and spheres}

Let $X$ be one of the two geometries, $X \in \{ \SXR, \HXR \}$. Using the geodesic curves we introduced in 
\cite{M97, PSSz10, Sz11-2, PSz12, Sz13-1, Sz11-1, PSSz11-2, 
Sz12-5} the notions of geodesic distances in $X$. 
\begin{definition}
The distance $d^{X}(P_1,P_2)$ between the points $P_1$ and $P_2$ is 
defined by the arc length of the geodesic curve 
from $P_1$ to $P_2$.
\end{definition}
\begin{definition}
 The geodesic sphere of radius $\rho$ (denoted by $S^{X}_{P_1}(\rho)$) with centre at 
 point $P_1$ is defined as the set of all points 
 $P_2$ in the space with the condition $d^{X}(P_1,P_2)=\rho$. Moreover, we require that the geodesic sphere is a simply connected 
 surface without self-intersection in $X$ space (Fig.1,~2).
 \end{definition}
 \begin{definition}
 The body of the geodesic sphere with centre $P_1$ and radius $\rho$ in space $X$ is called 
 geodesic ball, and denoted by $B_{P_1}(\rho)$,
 i.e. $Q \in B_{P_1}(\rho)$ iff $0 \leq d(P_1,Q) \leq \rho$.
 \end{definition}
 \begin{proposition}
 A geodesic sphere and ball of radius $\rho$ exists in the $\SXR$ space if and only if $\rho \in [0,\pi].$
 \end{proposition}
 \begin{proposition}
 $S(\rho)$ is a simply connected surface in $\HXR$ for $\rho >0$.
 \end{proposition}
 \begin{figure}[ht]
 \centering
 \includegraphics[width=5cm]{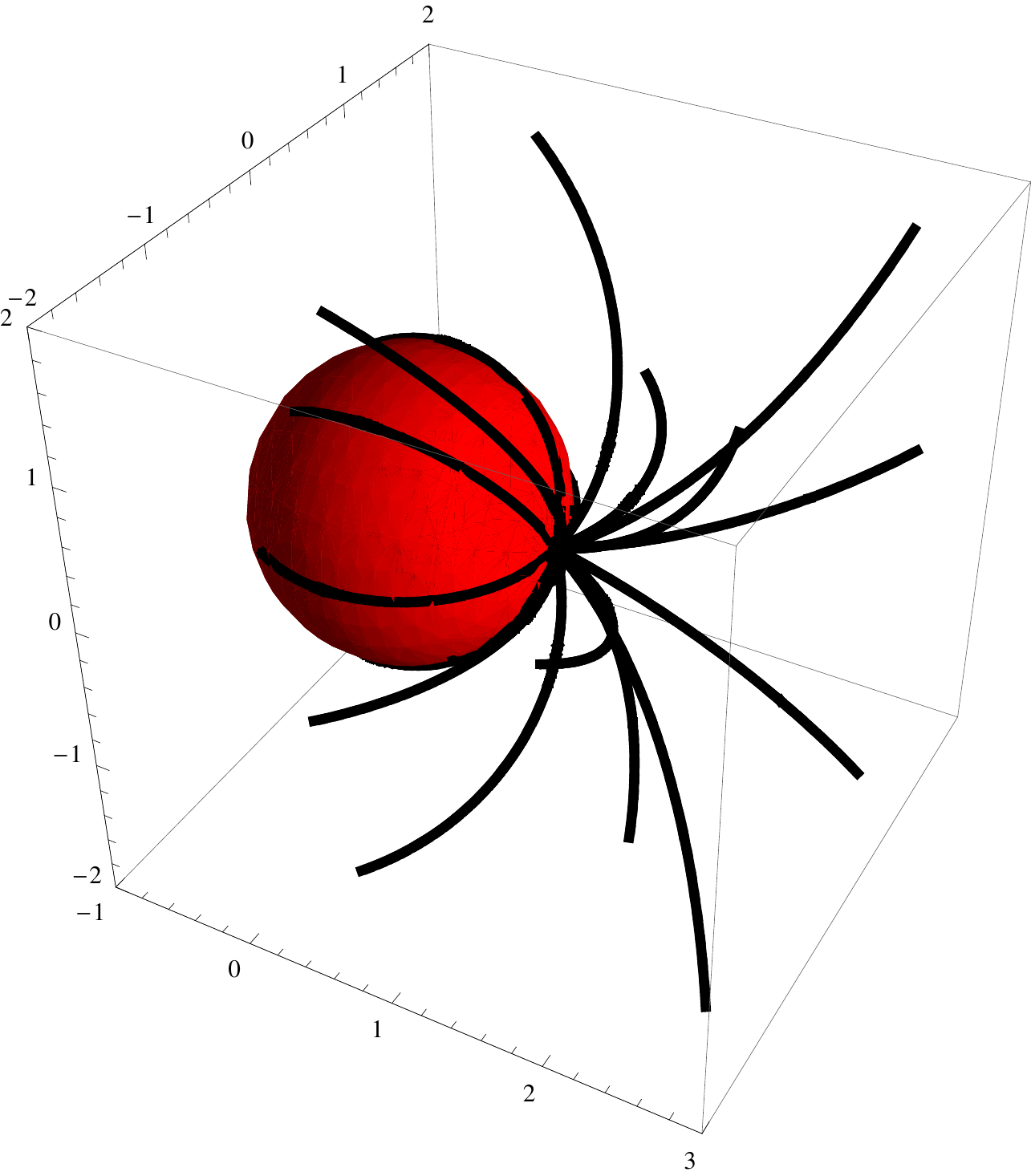} \includegraphics[width=5cm]{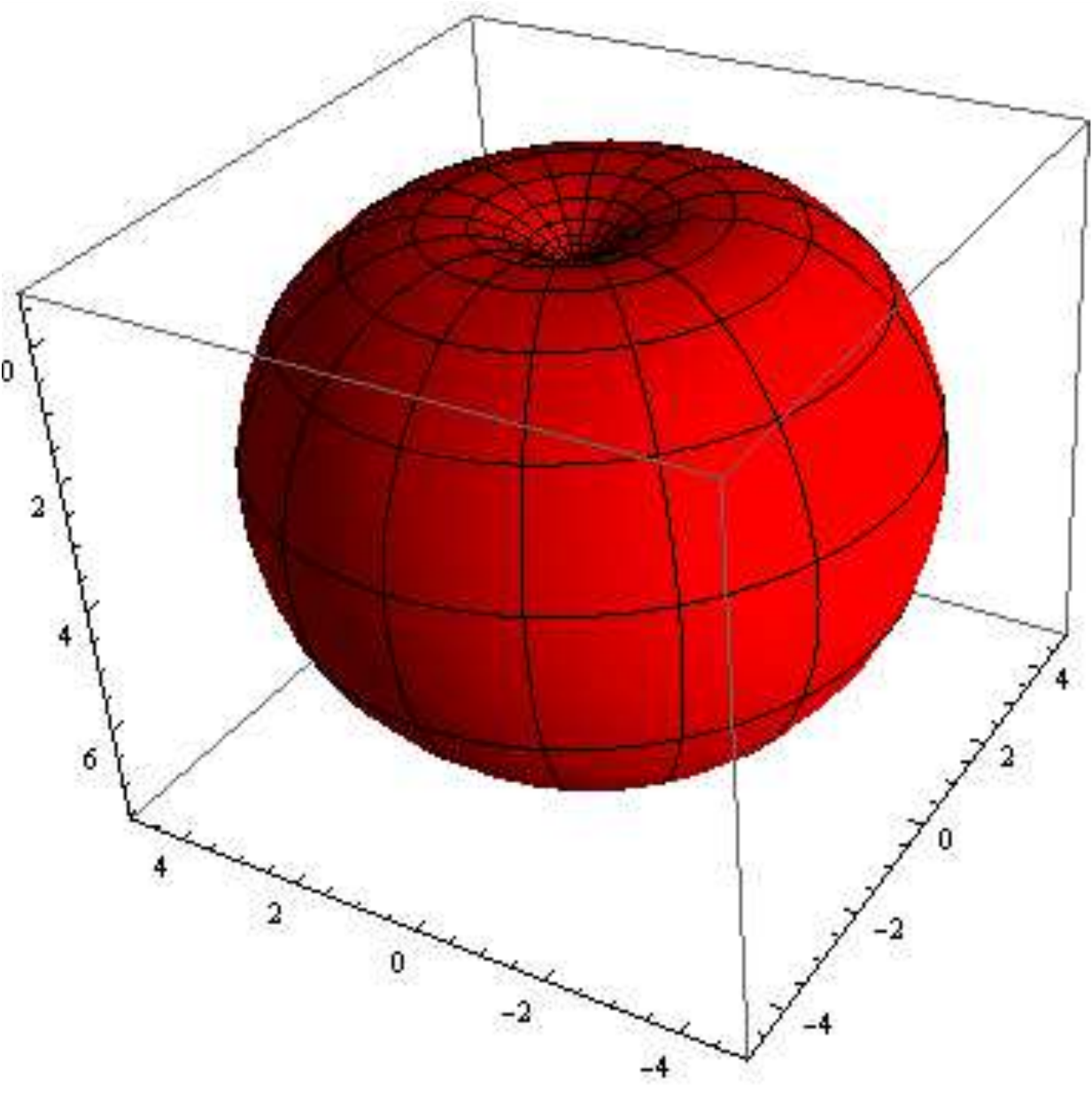}
 \caption{a. Geodesics with varying parameters. b. The geodesic sphere with radius $2$ centered at $(1,1,0,0)$.}
 \label{pic:spider}
 \end{figure}
 \begin{figure}[ht]
 \centering
 \includegraphics[width=5cm]{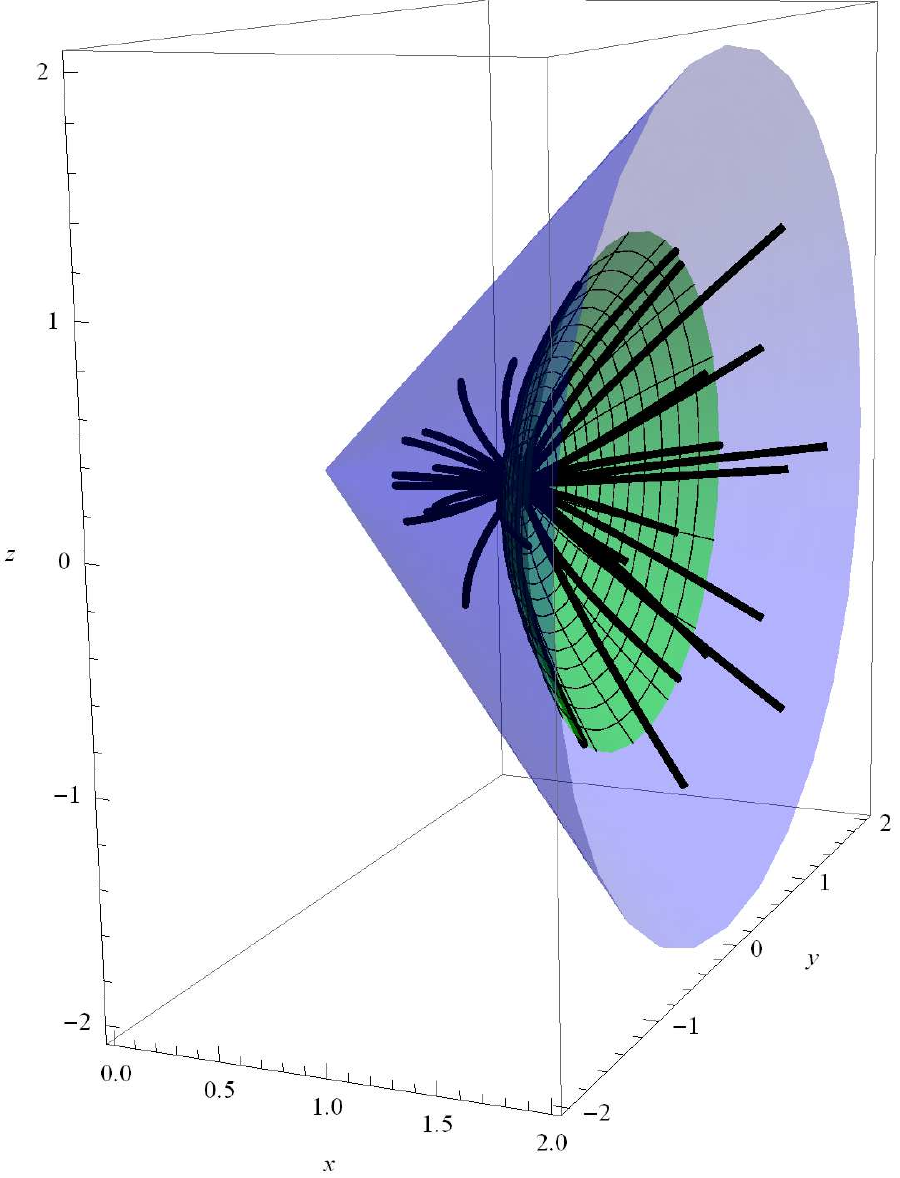} \includegraphics[width=5cm]{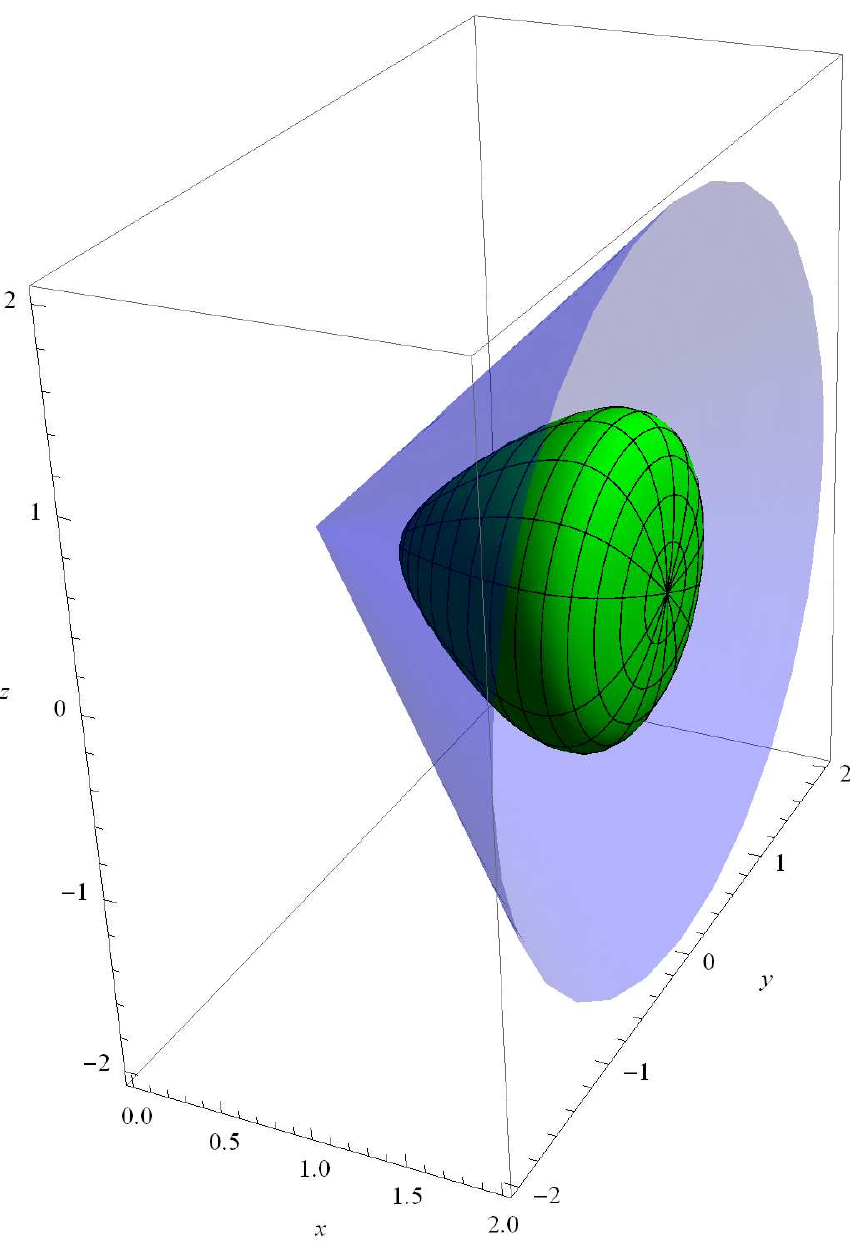}
 \caption{a. Geodesics with varying parameters  b. 
 The "base-hyperboloid" in the cone and a geodesic sphere with radius $\frac{2}{3}$ centered at $(1,1,0,0)$ in $\HXR$.}
 \label{pic:spider}
 \end{figure}
\subsection{Geodesic triangles and their interior angle sums}
We recall the important notions related to the interior angle sums of geodesic triangles in the above geometries elaborated in \cite{Sz20}. 

A geodesic triangle in Riemannian geometry and more generally in metric geometry a
is figure consisting of three different points together with their pairwise-connecting geodesic curves.
The points are known as the vertices, while the geodesic curve segments are known as the sides of the triangle.

In the geometries of constant curvature $\EUC$, $\HYP$, $\SPH$ the well-known sums of the interior angles of geodesic 
triangles characterize the space. It is related to the Gauss-Bonnet theorem which states that the integral of the Gauss curvature
on a compact $2$-dimensional Riemannian manifold $M$ is equal to $2\pi\chi(M)$ where $\chi(M)$ denotes the Euler characteristic of $M$.
This theorem has a generalization to any compact even-dimensional Riemannian manifold (see e.g. \cite{Ch}, \cite{KN}).

Therefore, it is interesting to investigate the interior angle sums of geodesic triangles in Thurston geometries.

We consider $3$ points $A_1$, $A_2$, $A_3$ in the projective model of $X$ space (see Section 2) $(X\in\{\SXR, \HXR \}$.
The {\it geodesic segments} $a_k$ connecting the points $A_i$ and $A_j$
$(i<j,~i,j,k \in \{1,2,3\}, k \ne i,j$) are called sides of the {\it geodesic triangle} with vertices $A_1$, $A_2$, $A_3$ (see Fig.~4).

In Riemannian geometries the infinitesimal arc length square (see (3.1) and (3.5)) is used to define the angle $\theta$ between two geodesic curves.
If their tangent vectors at their common point are $\bu$ and $\bv$ and $g_{ij}$ are the components of the metric tensor then
\begin{equation}
\cos(\theta)=\frac{u^i g_{ij} v^j}{\sqrt{u^i g_{ij} u^j~ v^i g_{ij} v^j}} \tag{3.9}
\end{equation}
Considering a geodesic triangle $A_1A_2A_3$ we can assume by the homogeneity of 
the geometries considered that one vertex 
coincides with the point $A_1=(1,1,0,0)$ and the other two vertices are $A_2=(1,x_2,y_2,z_2)$ and $A_3=(1,x_3,y_3,z_3)$. 

We denote the {\it interior angles} of geodesic triangles that are denoted at vertex 
$A_i$ by $\omega_i$ $(i\in\{1,2,3\})$.

We analysed in \cite{Sz20} the interior angle sums of geodesic triangles in both geometries. 
The answer this question is a consequence of the comparison theorems
in Riemannian geometry (Toponogov and Alexandrov's theorems, see \cite{CE}), since the sectional curvature of $\SXR$
is non-negative and the sectional curvature of $\HXR$ is non-positive. 
\subsubsection{Interior angle sums in $\SXR$ geometry}
We gave a new direct approach to this guestion 
based on the projective model and the isometry groups 
of $\SXR$ and $\HXR$ in \cite{Sz20}. 

We directly obtained from equations (3.4) of geodesic curves the following
\begin{lemma}[\cite{Sz20}]
Let $P$ be an arbitrary point and $g(A_1,P)$ ($A_1=(1,1,0,0)$) be a geodesic curve in 
the projective model of $\SXR$ geometry. 
The points of the geodesic curve $g(A_1,P)$ and the centre of the model $E_0$ lie in a 
plane in the Euclidean sense (see Fig.~3). 
\end{lemma}
\begin{figure}[ht]
\centering
\includegraphics[width=12cm]{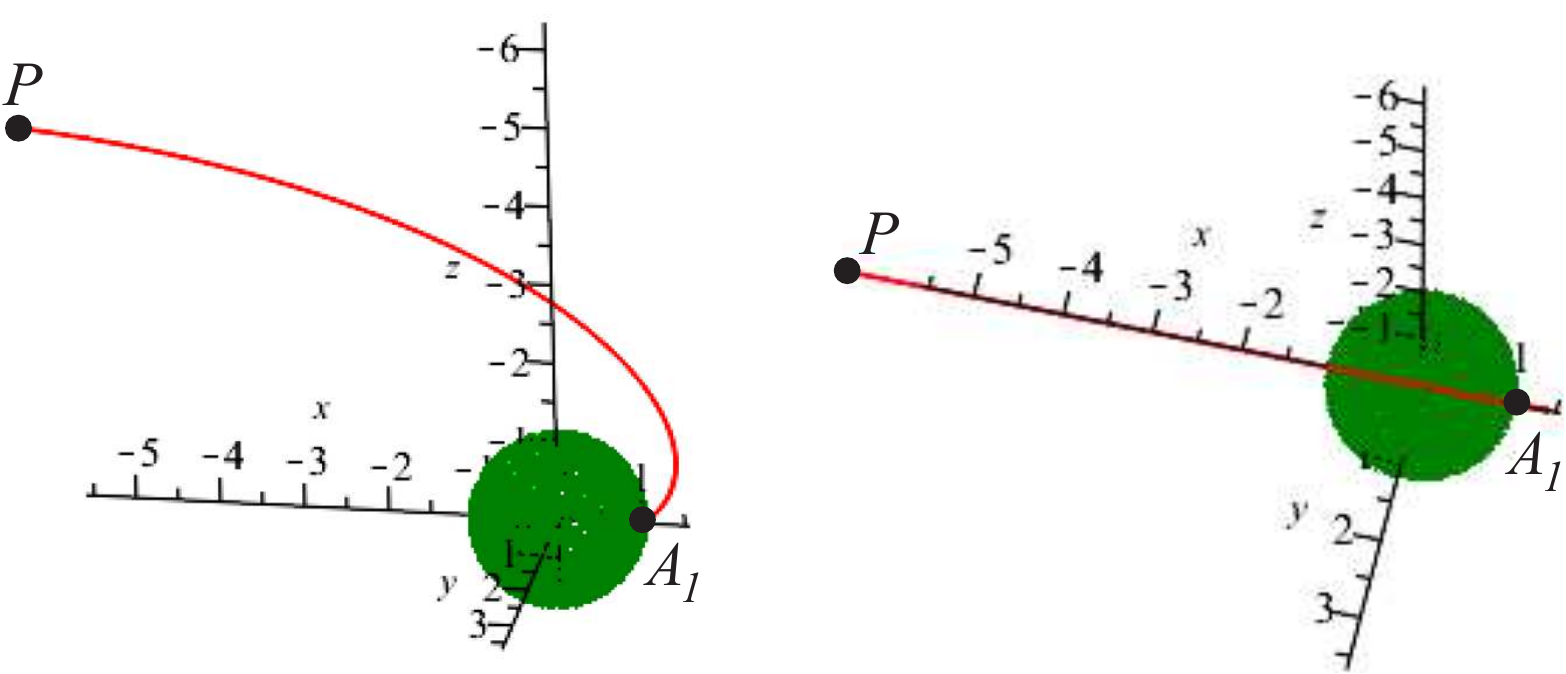}
\caption{Geodesic curve $g(A_1,P)$ ($A_1=(1,1,0,0)$ and $P \in \SXR$) with ``base plane", 
the plane of the geodesic curve contains the origin $E_0=(1,0,0,0)$ of the model \cite{Sz20}.}
\label{}
\end{figure}
\begin{Theorem}[\cite{Sz20}]
If the Euclidean plane of the vertices of an $\SXR$ geodesic triangle $A_1A_2A_3$ contains the centre of model $E_0$ then its 
interior angle sum is equal to $\pi$.
\end{Theorem}
\begin{figure}[ht]
\centering
\includegraphics[width=13cm]{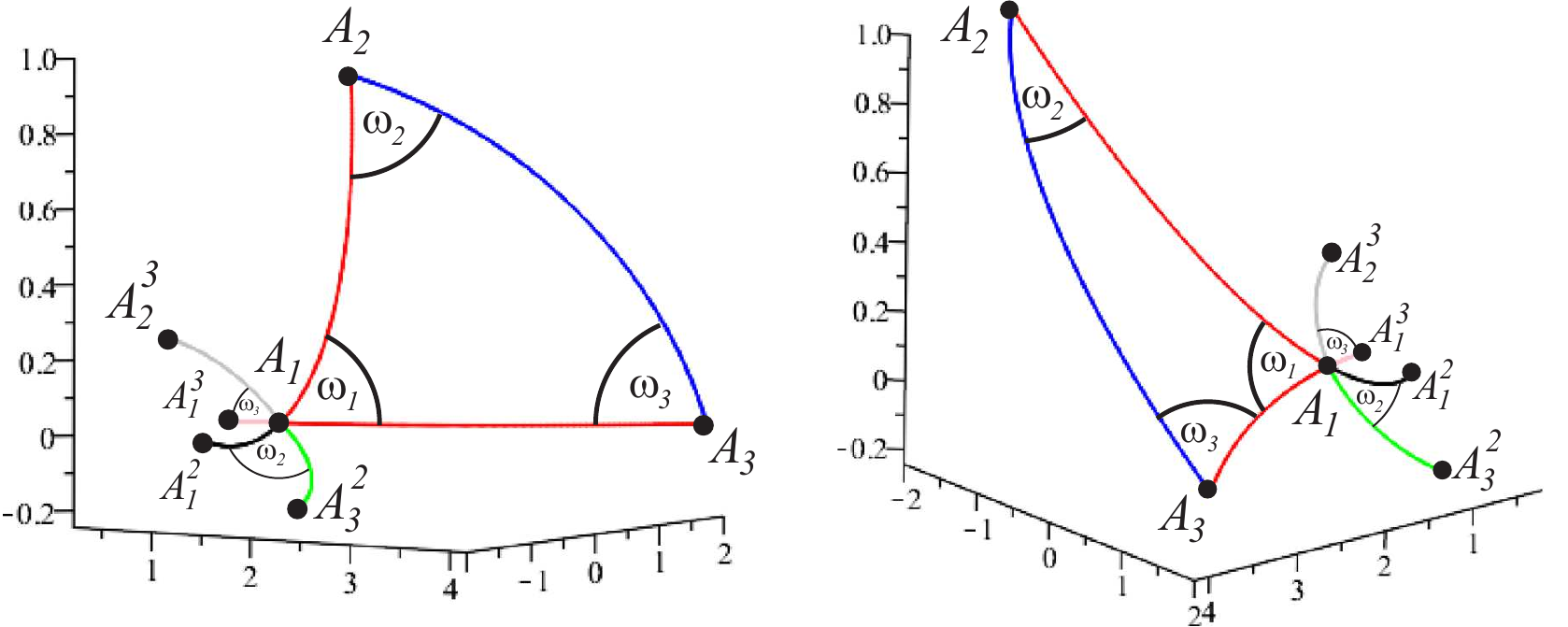}
\caption{Geodesic triangle with vertices $A_1=(1,1,0,0)$, $A_2=(1,3,-2,1)$, $A_3=(1,2,1,0)$ in $\SXR$ geometry, and transformed images of its geodesic side segments \cite{Sz20}.}
\label{}
\end{figure}
We could determine the interior angle sum of arbitrary geodesic triangle (for some numerical 
results see \cite{Sz20}). 

\medbreak

We obtained the following 
\begin{Theorem}[\cite{Sz20}]
If the Euclidean plane of the vertices of a $\SXR$ geodesic triangle $A_1A_2A_3$ does not contain the centre of model $E_0$ then its 
interior angle sum is greater than $\pi$.
\end{Theorem}
\begin{rmrk}
It is well known that if the vertices $A_1,A_2,A_3$ lie on a sphere of radius $R\in \bR^+$ centred at $E_0$ then the interior angle sum of spherical triangle 
$A_1A_2A_3$ is greater than $\pi$.
\end{rmrk}
Let $\Delta^{\SXR}(t)$ $(t\in\bR^+)$ denote the above geodesic triangle with {\it interior angles} $\omega_i(t)$ at 
the vertex $A_i$ $(i\in\{1,2,3\})$.  

The interior angle sum function $S(\Delta(t))=\sum_{i=1}^3(\omega_i(t))$ can be determined relative to the parameters
$x_2,y_2,z_2,x_3,y_3 \in \bR$ by the formulas (3.4), (3.6), (3.7) and by Lemma 1.
Analyzing the above complicated continuous functions of single real variable $t$ we get that its maximum is achieved at a point $t_0 \in (0, \infty)$ depending on given parameters. 
Moreover, $S(\Delta^{\SXR}(t))$ is stricly increasing on the interval $(0,t_0)$, stricly decreasing on the interval $(t_0,\infty)$ and
$$
\lim_{t \rightarrow 0}S(\Delta^{\SXR}(t))=\pi, ~ ~ ~ ~ \lim_{t \rightarrow \infty} S(\Delta^{\SXR}(t))=\pi. \hspace{1cm} 
$$ 
In summary we have the following
\begin{Theorem}[\cite{Sz20}]
The sum of the interior angles of a geodesic triangle of $\SXR$ space is greater than or equal to $\pi$. 
\end{Theorem}
\subsubsection{Interior angle sums in $\HXR$ geometry}
Similarly to the $\SXR$ space we investigated the interior angles of a geodesic triangle $A_1A_2A_3$ 
and its interior angle sum $\sum_{i=1}^3(\omega_i)$ in the $\HXR$ space (see \cite{Sz20} and Fig.~6).
\begin{lemma}[\cite{Sz20}]
Let $P$ be an arbitrary point and $g(A_1,P)$ ($A_1=(1,1,0,0)$) is a geodesic curve in the projective
model of $\HXR$ geometry. 
The points of the geodesic curve $g(A_1,P)$ and the centre of the model $E_0$ lie in a plane in Euclidean sense (see Fig.~5).
\end{lemma}
\begin{figure}[ht]
\centering
\includegraphics[width=12cm]{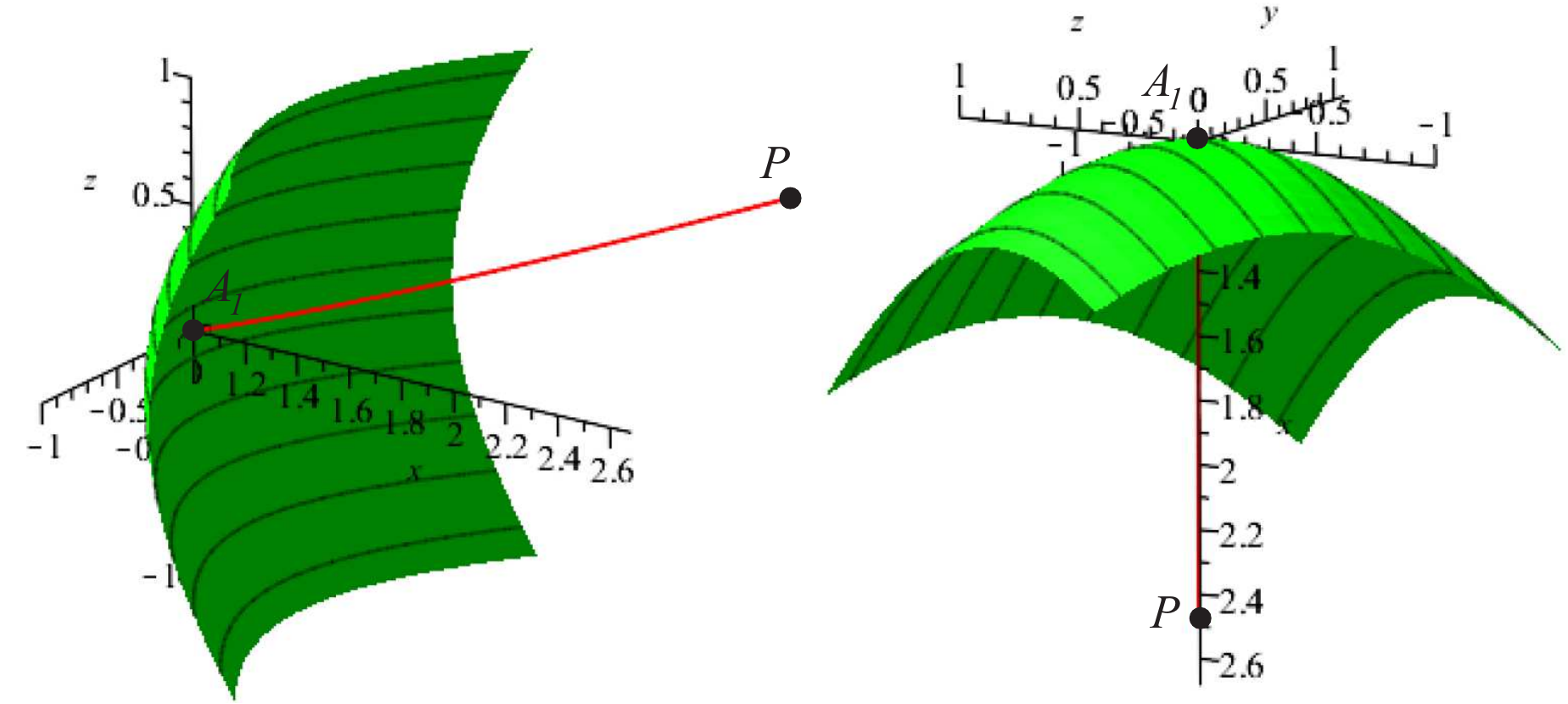}
\caption{Geodesic curve $g(A_1,P)$ ($A_1=(1,1,0,0)$ and $P \in \HXR$) with ``base plane" (the "upper" sheet of the two-sheeted hyperboloid), the plane of a geodesic curve contains 
the origin $E_0=(1,0,0,0)$ of the model \cite{Sz20}.}
\label{}
\end{figure}
\begin{rmrk}
More information about the isometry group of $\HXR$ and about its discrete subgroups can be found in \cite{Sz12-5}.  
\end{rmrk}
\begin{Theorem}[\cite{Sz20}]
If the Euclidean plane of the vertices of a $\HXR$ geodesic triangle $A_1A_2A_3$ contains the centre of model $E_0$ then its 
interior angle sum is equal to $\pi$. 
\end{Theorem}
\begin{figure}[ht]
\centering
\includegraphics[width=13cm]{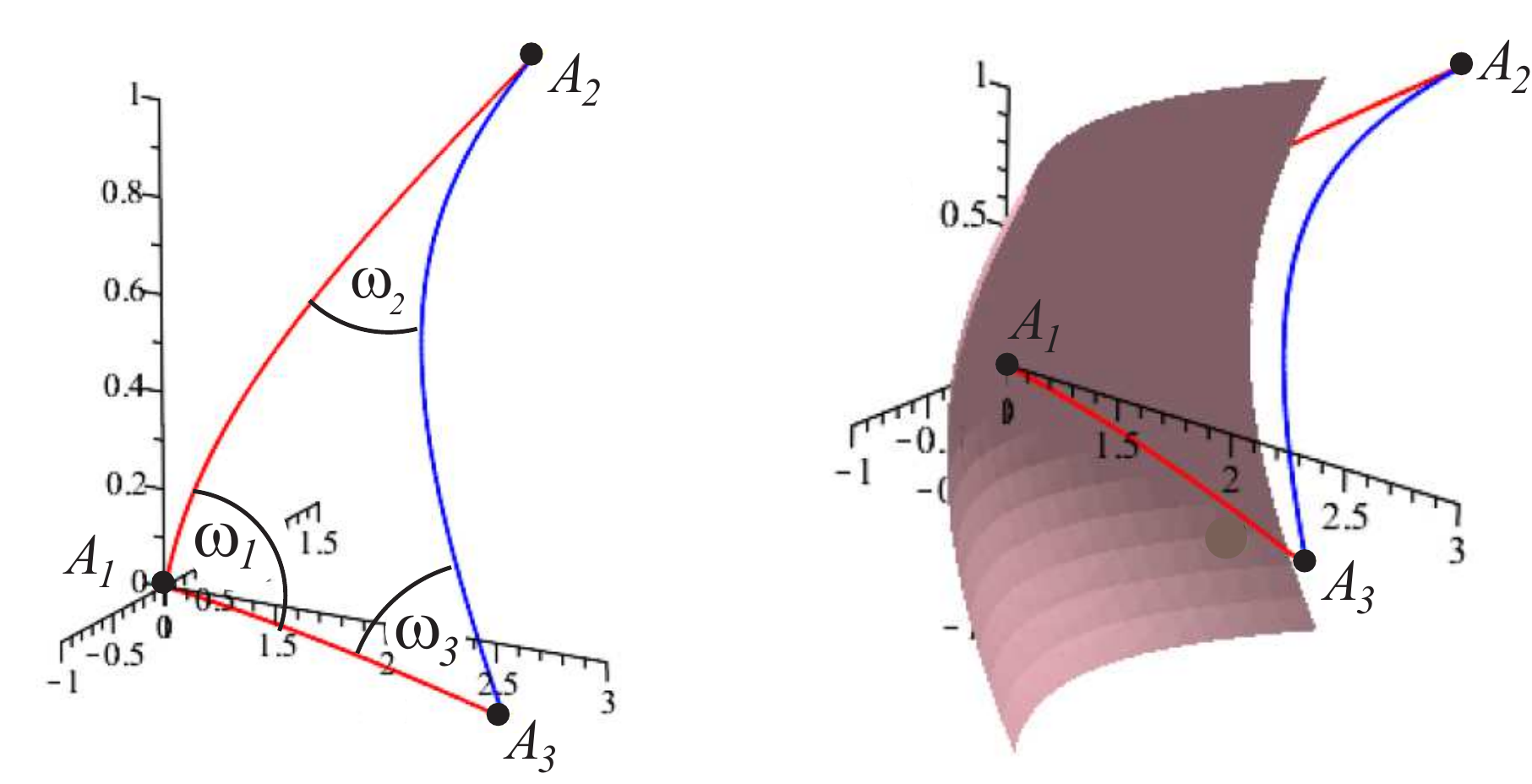}
\caption{Geodesic triangle with vertices $A_1=(1,1,0,0)$, $A_2=(1,2,3/2,1)$, $A_3=(1,3,-1,0)$ in $\HXR$ geometry \cite{Sz20}.}
\label{}
\end{figure}
We could determine the interior angle sum of arbitrary $\HXR$ geodesic triangle.

We obtain the following 
\begin{Theorem}[\cite{Sz20}]
If the Euclidean plane of the vertices of a $\HXR$ geodesic triangle $A_1A_2A_3$ does not contain the centre of model $E_0$ then its 
interior angle sum is less than $\pi$.
\end{Theorem}
\begin{rmrk}
It is well known that if the vertices $A_1,A_2,A_3$ lie on an "upper" sheet of the two-sheeted hyperboloid (in the hyperboloid model of the hyperbolic plane geometry where 
the straight lines of hyperbolic 2-space are modeled by geodesics on the hyperboloid) centred at $E_0$ then the interior angle sum of hyperbolic triangle 
$A_1A_2A_3$ is less than $\pi$.
\end{rmrk}
Let $\Delta(t)$ $(t\in\bR^+)$ denote the above geodesic triangle with {\it interior angles} $\omega_i(t)$ at 
the vertex $A_i$ $(i\in\{1,2,3\})$.  

The interior angle sum function $S(\Delta^{\HXR}(t))=\sum_{i=1}^3(\omega_i(t))$ can be determined relative to the parameters
$x_2,y_2,z_2,x_3,y_3 \in \bR$. 
$$
\lim_{t \rightarrow 0}S(\Delta^{\HXR}(t))=\pi, ~ ~ ~ ~ \lim_{t \rightarrow \infty} S(\Delta^{\HXR}(t))=\pi. \hspace{1cm} 
$$ 

In summary we obtained the following
\begin{Theorem}[\cite{Sz20}]
The sum of the interior angles of a geodesic triangle of $\HXR$ space is less than or equal to $\pi$. 
\end{Theorem}
\subsection{Surfaces of geodesic triangles }
We consider $3$ points $A_0$, $A_1$, $A_2$ in the projective model of the space $X$ (see Section 2, subsections 3.1-2) $(X\in\{\SXR, \HXR \})$.
The {\it geodesic segments} $a_k$ connecting the points $A_i$ and $A_j$
$(i<j,~i,j,k \in \{0,1,2\}, k \ne i,j$) are called sides of the {\it geodesic triangle} with vertices $A_0$, $A_1$, $A_2$.

The definition of the surface of a geodetic triangle in the space $X$ is not straightforward. The usual geodesic triangle surface definition in these geometries is 
not possible because the geodesic curves starting from different vertices and ending at points of the corresponding opposite edges define different 
surfaces, i.e. 
{\it geodesics starting from different vertices and ending at points on the corresponding 
opposite side usually do not intersect.}
Therefore, we introduced the definition (see [\cite{Sz21}) of the surface $\mathcal{S}_{A_0A_1A_2}$ of the geodesic triangle using the notion of the generalized Apollonius surfaces:
\begin{definition}
The Apollonius surface $\mathcal{A}\cS^X_{P_1P_2}(\lambda)$ in the Thurston 
geometry $X$ is  
the set of all points of $X$ whose geodesic distances from 
two fixed points are in a constant ratio $\lambda\in\mathbf{R}^+_0$ where 
$X\in \EUC,\SPH,\HYP,\SXR,\HXR,\NIL,\SLR,\SOL.
$ i.e. $\mathcal{A}\cS^X_{P_1P_2}(\lambda)$ of two arbitrary points 
$P_1,P_2 \in X$ consists of all points $P'\in X$,
for which $d^X(P_1,P')=\lambda \cdot d^X(P',P_2)$ 
($\lambda\in [0,\infty$)) where $d^X$ is the corresponding distance function of $X$. If $\lambda=0$, 
then $\mathcal{A}\cS^X_{P_1P_2}(0):=P_1$ and it is clear, that in case $\lambda \to \infty$ then $d(P',P_2) \to 0$ therefore we say $\mathcal{A}\cS^X_{P_1P_2}(\infty):=P_2$.
\end{definition}
\begin{definition}
\begin{enumerate}
\item We consider the geodesic triangle $A_0A_1A_2$ in the projective model of $X$ space $(X\in\{\SXR, \HXR \})$ and consider the Apollonius surfaces 
$\mathcal{A}\cS^{X}_{A_0A_1}$ $(\lambda_1)$ and $\mathcal{A}\cS^{X}_{A_2A_0}(\lambda_2)$ ($\lambda_1,\lambda_2 \in [0,\infty)$, $\lambda_1^2+\lambda_2^2>0$). It is clear, that if 
$Y \in \mathcal{C}(\lambda_1,\lambda_2):=\mathcal{A}\cS^{X}_{A_0A_1}(\lambda_1)\cap \mathcal{A}\cS^{X}_{A_2A_0}(\lambda_2)$ then 
$\frac{d^X(A_0,Y)}{d^X(Y,A_1)}=\lambda_1$ and  $\frac{d^X(A_2,Y)}{d^X(Y,A_0)}=\lambda_2$ $\Rightarrow$ $\frac{d^X(A_2,Y)}{d^X(Y,A_1)}=\lambda_1 \cdot \lambda_2$
for parameters $\lambda_1,\lambda_2 \in (0,\infty)$ and if $\lambda_1=0$ then $\mathcal{C}(\lambda_1,\lambda_2)=A_0$, if $\lambda_2=0$ then $\mathcal{C}(\lambda_1,\lambda_2)=A_2$ 
\item 
\begin{equation}
\begin{gathered}
P^X(\lambda_1,\lambda_2):=\{ P \in X~|~P \in \mathcal{C}(\lambda_1,\lambda_2)~{\text{and}}~ d^X(P,A_0)=\min_{Q \in \mathcal{C}(\lambda_1,\lambda_2)}({d^X(Q,A_0)}) \\
~\text{with given real parameters}~\lambda_1,\lambda_2 \in [0,\infty),~ \lambda_1^2 + \lambda_2^2 > 0 \}
\end{gathered} \tag{3.10}
\end{equation}
\item The surface $\mathcal{S}^X_{A_0A_1A_2}$ of the geodesic triangle $A_0A_1A_2$ is 
\begin{equation}
\mathcal{S}_{A_0A_1A_2}^X:=\{P^X(\lambda_1,\lambda_2) \in X,~ \text{where}~\lambda_1,\lambda_2 \in [0,\infty),~ \lambda_1^2 + \lambda_2^2 > 0 \}. \tag{3.11}
\end{equation}
\end{enumerate}
\end{definition}
This method clearly leads to the following implicit equation of the Apollonius surfaces $\cA\cS^X_{P_1P_2}(\lambda)$ of two proper points 
$P_1(1,a,b,c)$ and $P_2(1,d,e,f)$ with given ratio $\lambda \in \bR^+_0$, 
in $X$ geometry:
\begin{Theorem}[\cite{Sz21}]
The implicit equation of the Apollonius surfaces $\cA\cS^X_{P_1P_2}(\lambda)$ 
of two proper points $P_1(1,a,b,c)$ and $P_2(1,d,e,f)$ with given ratio 
$\lambda \in \bR^+_0$, in the $X$ geometry is:
\begin{equation}\label{heqeq2}
\begin{gathered}
\cA\cS^X_{P_1P_2}(\lambda)(x,y,z) \Rightarrow \\
4 \omega_X^2{\left(\frac{ax \pm by \pm cz}{\sqrt{a^2 \pm b^2 \pm c^2}\sqrt{x^2 \pm y^2 \pm z^2}} \right)}+{\log^2{\Big(\frac{a^2 \pm b^2 \pm c^2}{x^2 \pm y^2 \pm z^2}\Big)}}=\\
=\lambda^2 \Big[4 \omega_X^2{\left(\frac{dx \pm ey \pm fz}{\sqrt{d^2 \pm e^2 \pm f^2}\sqrt{x^2 \pm y^2 \pm z^2}} \right)}+{\log^2{\Big(\frac{d^2 \pm e^2 \pm f^2}
{x^2 \pm y^2 \pm z^2}\Big)}}\Big],
\end{gathered} \notag
\end{equation}
where if $X=\SXR$ then all $\pm$ signs are $+$, $\omega_X(x)=\arccos(x)$ 
and if $X=\HXR$ then the all $\pm$ signs are $-$, $\omega_X(x)=\arccosh(x)$ (see Fig.~7,~8).  
\end{Theorem}
\begin{figure}[ht]
\centering
\includegraphics[width=13cm]{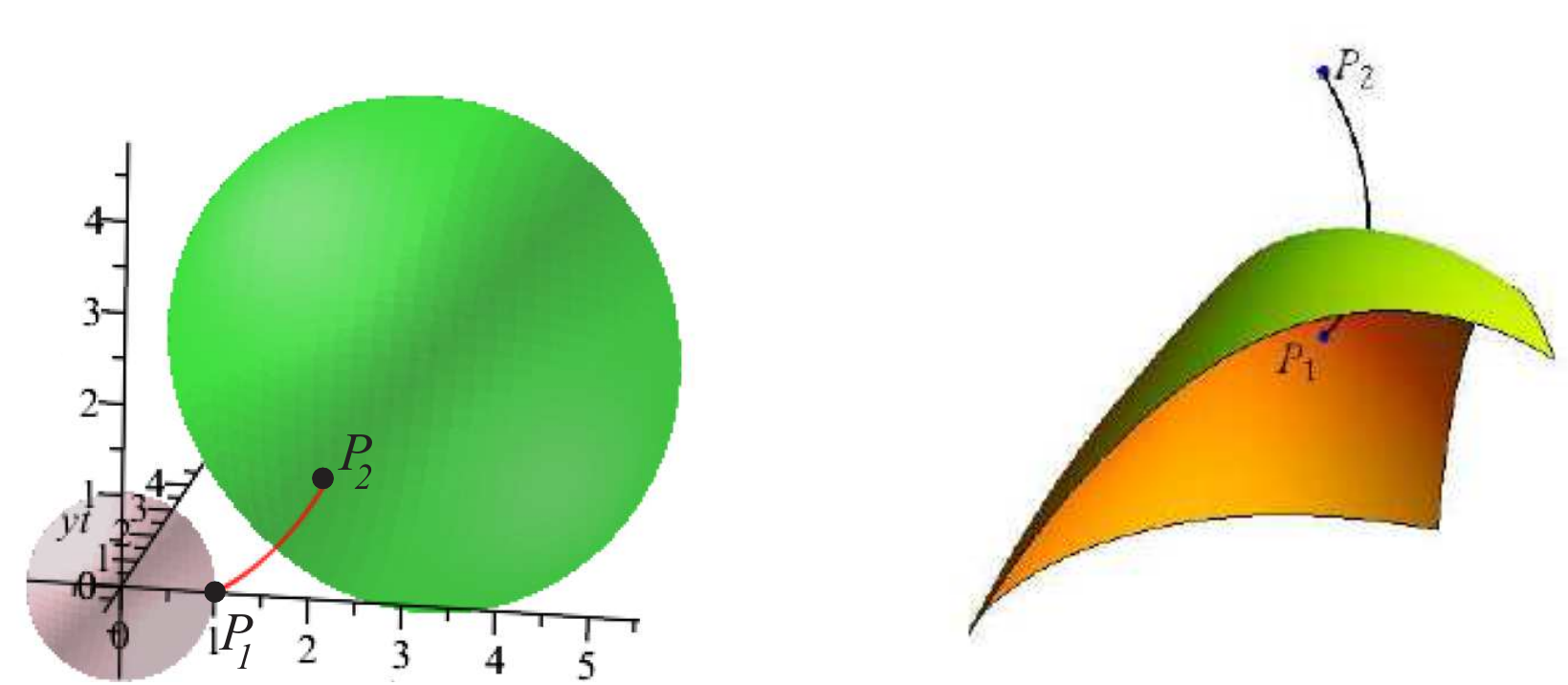}
\caption{The Apollonius surface $\mathcal{A}\cS^{\SXR}_{P_1P_2}(\lambda)$ where $P_1=(1,1,0,0)$, $P_2=(1,2,1,1)$, $\lambda=2$ (left) and $\lambda=1$ (right, equidistance surface (see 
\cite{Sz21, PSSz10})}
\label{}
\end{figure}
\begin{figure}[ht]
\centering
\includegraphics[width=13cm]{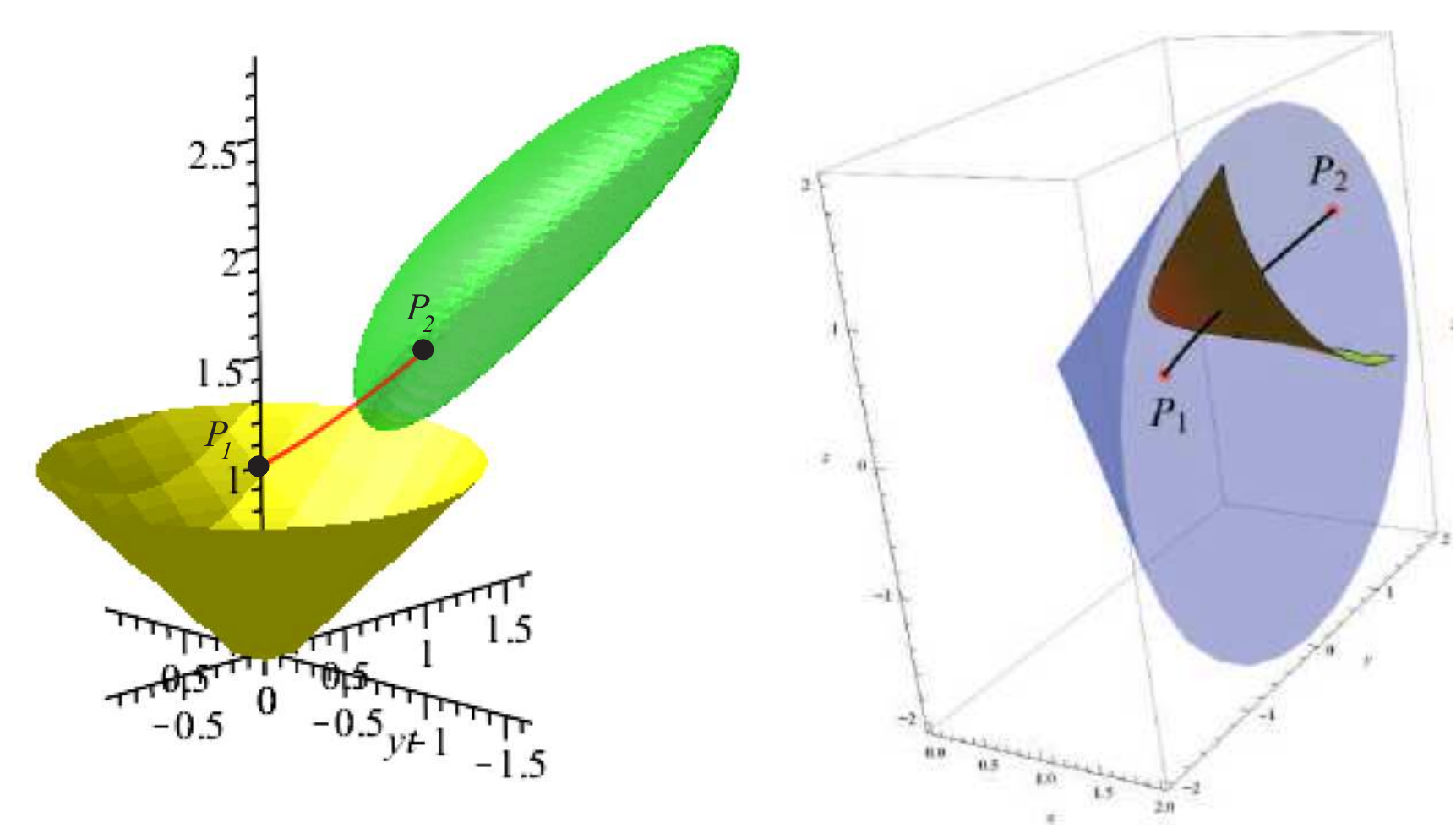}
\caption{The Apollonius surface $\mathcal{A}\cS^{\HXR}_{P_1P_2}(\lambda)$ where $P_1=(1,1,0,0)$, $P_2=(1,3/2,1,-1/2)$, $\lambda=2$ (left) and $P_1=(1,1,0,0)$, 
$P_2=(1,2,1,1)$, $\lambda=1$ (right, equidistance surface) (see \cite{Sz21, PSSz11-2}}
\label{}
\end{figure}
We used the statement of the following lemma
\begin{lemma}[\cite{Sz20}]
Let $P$ be an arbitrary point and $g^{X}(P_1,P)$ $(X \in \{ \SXR, \HXR \}$, $P_1=(1,1,0,0)$) is a geodesic curve in the considered model of $X$ geometry. 
The points of the geodesic curve $g^X(P_1,P)$ and the centre of the model $E_0=(1,0,0,0)$ lie in a Euclidean plane. 
\end{lemma}
\subsection{Geodesic tetrahedra and their circumscribed spheres}
We consider $4$ points $A_0$, $A_1$, $A_2$, $A_3$ in 
the projective model of the space $X$ (see Section 2, subsections 3.1-2 $X\in\{\SXR, \HXR \})$.
These points are the vertices of a {\it geodesic tetrahedron} in the space $X$ if any two {\it geodesic segments} connecting the points $A_i$ and $A_j$
$(i<j,~i,j \in \{0,1,2,3\})$ do not have common inner points and any three vertices do not lie on the 
same geodesic curve.
Now, the geodesic segments $A_iA_j$ are called edges of the geodesic tetrahedron $A_0A_1A_2A_3$.
\begin{figure}[ht]
\centering
\includegraphics[width=12cm]{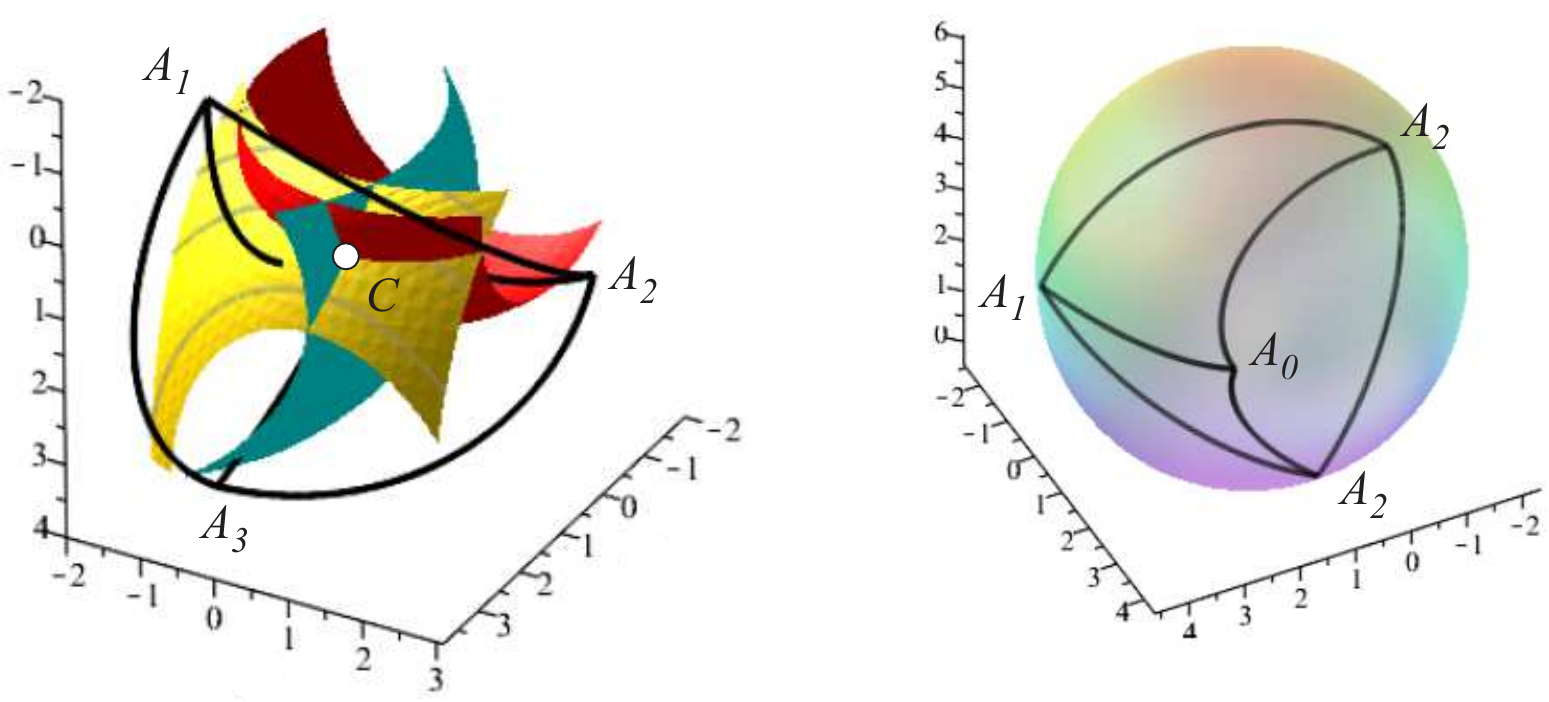}
\caption{Geodesic $\SXR$ tetrahedron with vertices $A_0=(1,1,0,0)$, $A_1=(1,-2, -1/2,3)$, $A_2=(1,1,3,0)$, $A_3=(1,4,-1,2)$
and its circumscibed sphere of radius $r \approx 1.30678$ with circumcenter $C=(1,\approx 0.64697, \approx 0.51402, \approx 0.15171)$ \cite{Sz21}.}
\label{}
\end{figure}
The circumscribed sphere of a geodesic tetrahedron is a geodesic sphere that touches each of the tetrahedron's vertices.
As in the Euclidean case the radius
of a geodesic sphere circumscribed around a tetrahedron 
$T$ is called the circumradius of $T$, and the center point of this sphere 
is called the circumcenter of $T$ (Fig.~9,~10).
For any $\SXR$ geodesic tetrahedron there exists a unique geodesic surface on 
which all four vertices lie. If its radius less or equal to $\pi$ then 
the above surface is a geodesic sphere (called circumscibed sphere, see Definition 3.2). 
\begin{figure}[ht]
\centering
\includegraphics[width=12cm]{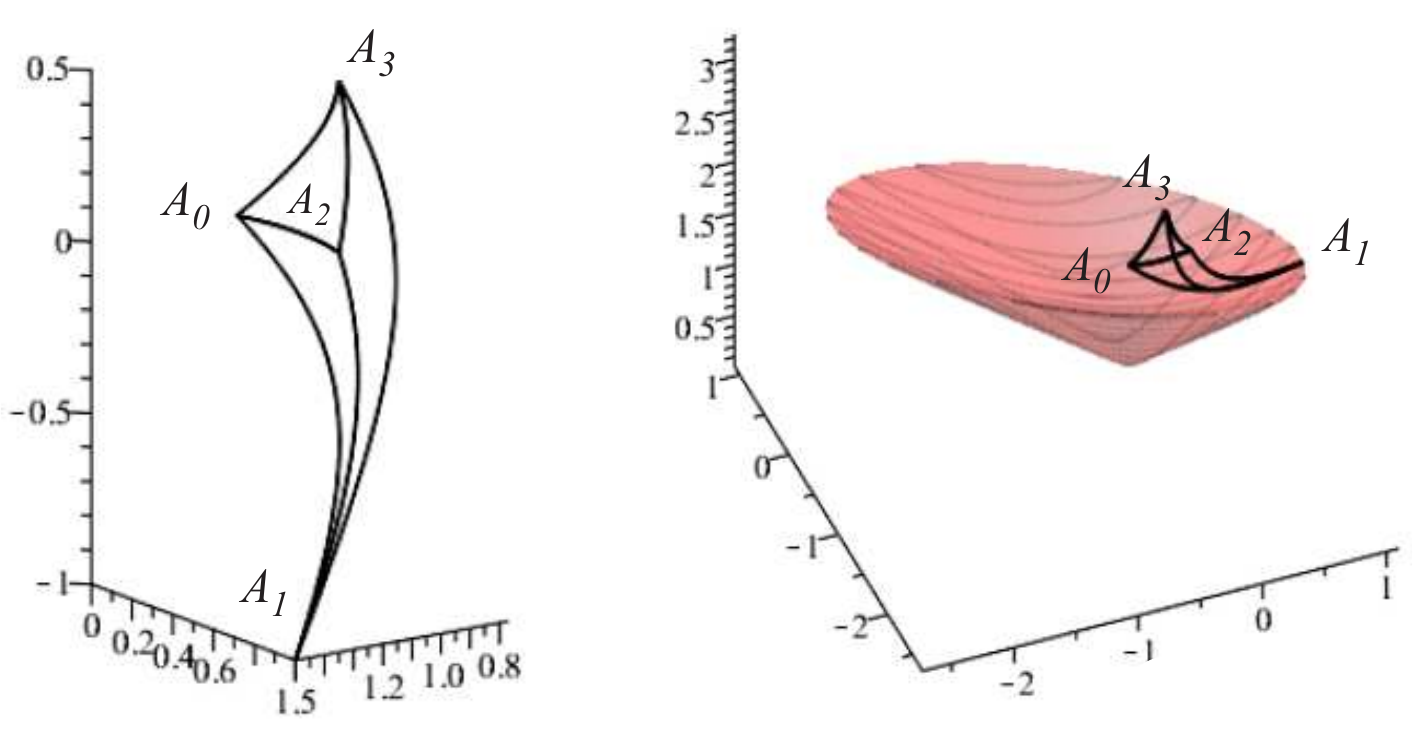}
\caption{Geodesic $\HXR$ tetrahedron with vertices $A_0=(1,1,0,0)$, $A_1=(1,3/2,1,-1)$, $A_2=(1,1,1/2,0)$, $A_3=(1,1,1/2,1/2)$
and its circumscibed sphere of radius $r \approx 2.89269$ with circumcenter $C=(1,\approx 0.07017, \approx -0.02714, \approx -0.02640)$ \cite{Sz21} }
\label{}
\end{figure}
\begin{rmrk}
If the common point of bisectors lies at infinity then the 
vertices of tetrahedron lie on a horosphere-like surface and if the 
common point is an outer point then the vertices of 
the tetrahedron are on a hypershere-like surface. 
These surfaces will be examined in detail in a forthcoming paper.
\end{rmrk}
\begin{figure}[ht]
\centering
\includegraphics[width=13cm]{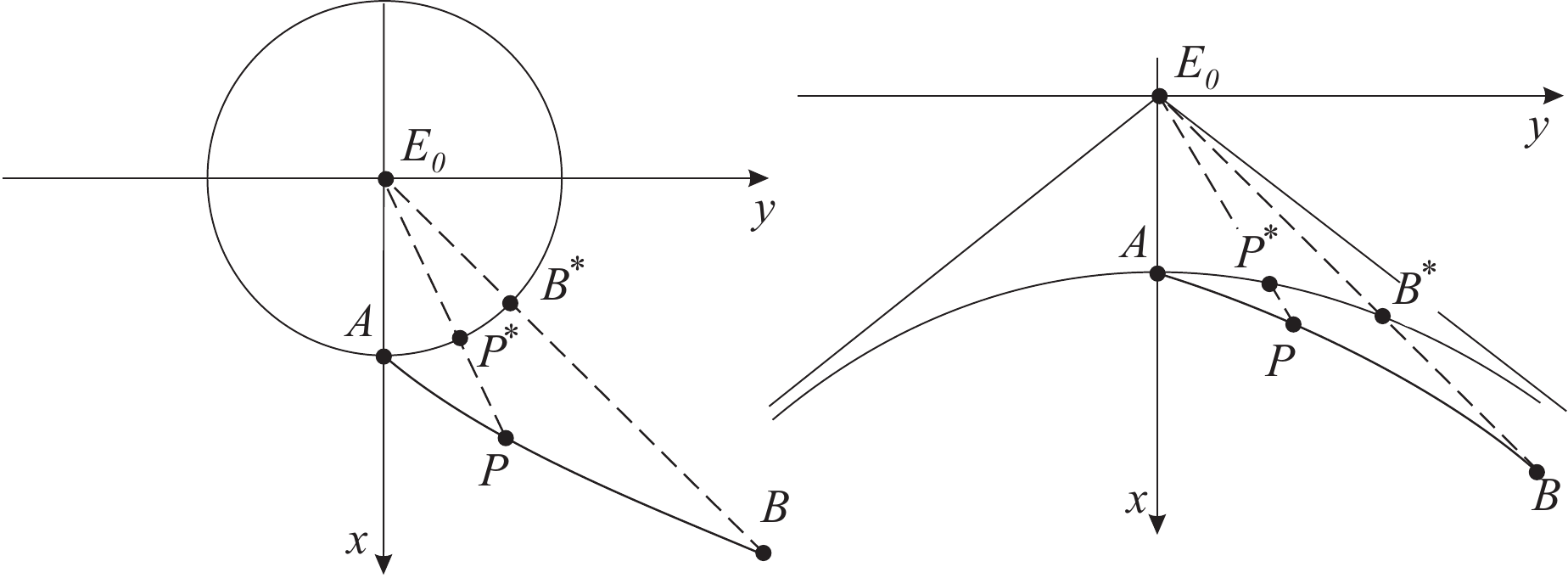}
\caption{Simple ratios in $\SXR$ and $\HXR$ spaces}
\label{}
\end{figure}
For any $\HXR$ geodesic tetrahedron there exists a unique geodesic surface on which all four vertices lie. If its centre is a proper point of $\HXR$ space then 
the above surface is a geodesic sphere (called a circumscibed sphere). 
\subsection{Menelaus' and Ceva's theorems in $\SXR$ and $\HXR$ spaces} 
First we recall the definition of simple ratios in the sphere $\bS^2$ and the 
plane $\bH^2$ (see \cite{PS14}).
The models of the above plane geometries of constant curvature are 
embedded in the models of the previously described geometries $\SXR$ and $\HXR$ 
as ``base planes" and are used hereinafter for our discussions. 

A spherical triangle is the space enclosed by arcs of great circles on the surface of a sphere, 
subject to the constraint that these arcs and the 
further circular arcs in the spherical 
plane are always less or equal than a semicircle.
\begin{definition}
If $A$, $B$ and $P$ are distinct points on a line in 
the $Y\in\{\bH^2,\bS^2\}$ space, then 
their simple ratio is
$s^Y(A,P,B) =  w^Y(d^Y(A,P))/w^Y(d^Y(P,B))$, if $P$ is between $A$ and $B$, and
$s^Y(A,P,B) = -w^Y(d^Y(A,P))/w^Y(d^Y(P,B))$, otherwise where
$w^Y(x):=sin(x)$ if $Y=\bS^2$ and $w^Y(x):=sinh(x)$ if $Y=\bH^2$. 
\end{definition}

\begin{rmrk} Basic properties of simple ratio:
\begin{enumerate}
\item $s^Y(A,P,B) = -s^Y(B,P,A)$,
\item if $P$ is between $A$ and $B$, then $s^Y(A,P,B) \in (0,1)$,
\item if $P$ is on $AB$, beyond $B$, then $s^Y(A,P,B) \in (-\infty,-1)$,
\item if $P$ is on $AB$, beyond $A$, then $s^Y(A,P,B) \in (-1,0)$.
\end{enumerate}
\end{rmrk}
Note that the value of $s^Y(A,P,B)$ determines the position of $Y$ relative to $A$ and $B$.

With this definition, the corresponding sine rule of the geometry $Y$ leads to 
Menelaus's and Ceva's theorems \cite{K,PS14}:

\begin{Theorem} [Menelaus's Theorem for triangles in the plane $Y$]
If is a line $l$ not passing through any vertex of an triangle $ABC$ such that
$l$ meets $BC$ in $Q$, $AC$ in $R$, and $AB$ in $P$,
then $$s^Y(A,P,B)s^Y(B,Q,C)s^Y(C,R,A) = -1.$$ 
\end{Theorem}

\begin{Theorem}[Ceva's Theorem for triangles in the plane $Y$]
If $T$ is a point not on any side of a triangle $ABC$ such that
$AT$ and $BC$ meet in $Q$, $BX$ and $AC$ in $R$, and $CX$ and $AB$ in $P$,
then $$s^Y(A,P,B)s^Y(B,Q,C)s^Y(C,R,A) = 1.$$ 
\end{Theorem}
\subsection{Generalizations of Menelaus' and Ceva's theorems}
\subsubsection{Geodesic triangle in general position}
First we consider a {\it general location geodesic triangle} 
$A_0A_1A_2$ in the projective model of the space $X$ (see subsections 3.1-2) $(X\in\{\SXR, \HXR \})$.
Without loss of generality, we can assume that $A_0=(1,1,0,0)$ and $A_2$ lies in the coordinate plane $[x,y]$. 
The geodesic lines that contain the sides $A_0A_1$ and $A_0A_2$ of the given triangle can be characterized directly by the corresponding parameters $v$ and $u$ (see (3.4) and (3.8)).

The geodesic curve including the side segment $A_1A_2$ is also determined by 
one of its endpoints and its parameters, however in order to determine the corresponding parameters of this 
geodesic line we use {\it orientation preserving isometric transformations} $\bT^{X}(A_2)$, as elements of the isometry group of the geometry $X$, that
maps $A_2=(1,x_2,y_2,0)$ onto $A_0=(1,1,0,0)$ (up to a positive determinant factor). 

We extend the definition of the simple ratio to the $X\in\{\SXR, \HXR \}$ spaces.
If $X=\SXR$ then is clear that the space contains its ``base sphere" (unit sphere centred in $E_0$) which is a geodesic surface. 
Therefore, similarly to the spherical spaces we assume that the geodesic arcs  
are always less than or equal to a semicircle.
\begin{definition}
If $A$, $B$ and $P$ are distinct points on a non-fiber-like geodesic curve in the $X \in \{\SXR, \HXR \}$ space , then 
their simple ratio is
$$s_g^X(A,P,B) =  w^X\Big({d^X(A,P)}{\cos(v)}\Big)/w^X\Big({d^X(P,B)}{\cos(v)}\Big),$$ if $P$ is between $A$ and $B$, and
$$s_g^X(A,P,B) = -w^X\Big({d^X(A,P)}{\cos(v)}\Big)/w^X\Big({d^X(P,B)}{\cos(v)}\Big),$$ 
otherwise where
$w^X(x):=sin(x)$ if $X=\SXR$, $w^Y(x):=sinh(x)$ if $X=\HXR$ and $v$ is 
the parameter of the geodesic curve containing points $A,B$ and $P$ (see Fig.~11).  
\end{definition}
\begin{Theorem}[Ceva's Theorem for triangles in general location, \cite{Sz21}]
If $T$ is a point not contained in any side of a geodesic triangle $A_0A_1A_2$ in $X\in\{\SXR, \HXR\}$ such that
the curves $A_0T$ and $g_{A_1A_2}^X$ meet in $Q$, $A_1T$ and $g_{A_0A_2}^X$ in $R$, and $A_2T$ and $g_{A_0A_1}^X$ in $P$, $(A_0T, A_1T, A_2T \subset \mathcal{S}^X_{A_0A_1A_2})$
then $$s_g^X(A_0,P,A_1)s_g^X(A_1,Q,A_2)s_g^X(A_2,R,A_0) = 1.$$ 
\end{Theorem}

\begin{Theorem}[Menelaus's theorem for triangles in general location, \cite{Sz21}]
If $l$ is a line not through any vertex of a geodesic triangle 
$A_0A_1A_2$ lying in a surface $\mathcal{S}^X_{A_0A_1A_2}$ in the $X\in\{\SXR, \HXR\}$ 
geometry such that
$l$ meets the geodesic curves $g_{A_1A_2}^X$ in $Q$, $g_{A_0A_2}^X$ in $R$, 
and $g_{A_0A_1}^X$ in $P$,
then $$s^X_g(A_0,P,A_1)s^X_g(A_1,Q,A_22)s^X_g(A_2,R,A_0) = -1.$$
\end{Theorem} 
\subsubsection{Fibre type triangle}
We consider a {\it fibre type geodesic triangle} $A_0A_1A_2$ in the projective model of the space $X$.
Without limiting generality, we can assume that $A_0=(1,1,0,0)$, $A_1=(1,x_1,y_1,0)$ and $A_2=(1,x_2,y_2,0)$ lie in the coordinate plane $[x,y]$. 
The geodesic lines that contain the sides $A_0A_1$ and $A_0A_2$ 
of the given triangle can be characterized directly by the corresponding parameters 
$v$ and $u=0$ similar to the above case.

We extend the definition of the simple ratio to the $X\in\{\SXR, \HXR \}$ spaces.

If $X=\SXR$ it is clear that the $\SXR$ space contains its ``base sphere" (unit sphere centred in $E_0$) which is a geodesic surface. 
Therefore, similarly to the spherical spaces we assume that the geodesic arcs  
are always less or equal than a semicircle.
\begin{definition}
If $A$, $B$ and $P$ are distinct points on fibrum-like geodesic 
curve in the $X \in \{\SXR, \HXR \}$ space , then 
their simple ratio is
$$s^X_f(A,P,B) =  d^X(A,P)/d^X(P,B)$$ if $P$ is between $A$ and $B$, and
$$s^X_f(A,P,B) = -d^X(A,P)/d^X(P,B)$$ (see Fig.~10).  
\end{definition}
\begin{Theorem}[Ceva's Theorem in the $X$ geometry for triangles in fibre types, \cite{Sz21}]
If $T$ is a point not on any side of a geodesic triangle $A_0A_1A_2$ in 
$X\in\{\SXR, \HXR\}$ such that
the geodesic curves $g_{A_0T}^X$ and $g_{A_1A_2}^X$ meet in $Q$, $g_{A_1T}^X$ and $g_{A_0A_2}^X$ in $R$, and $g_{A_2T}^X$ and $g_{A_0A_1}^X$ in $P$, 
$(g_{A_0T}^X, g_{A_1T}^X, g^X_{A_2T} \subset \mathcal{S}^X_{A_0A_1A_2}$)
then $$s^X_f(A_0,P,A_1)s^X_f(A_1,Q,A_2)s^X_f(A_2,R,A_0) = 1.$$
\end{Theorem}
\begin{Theorem}[Menelaus's theorem in the $X$ geometry for triangles in fibre types, \cite{Sz21}]
If $l$ is a line not through any vertex of an geodesic triangle $A_0A_1A_2$ lying in its surface $\mathcal{S}^X_{A_0A_1A_2}$ in $X\in\{\SXR, \HXR\}$ geometry such that
$l$ meets geodesic curves $g_{A_1A_2}^X$ in $Q$, $g_{A_0A_2}^X$ in $R$, and $g_{A_0A_1}^X$ in $P$,
then $$s^X_f(A_0,P,A_1)s^X_f(A_1,Q,A_2)s^X_f(A_2,R,A_0) = -1$$.
\end{Theorem} 
Thus, we can formulate similar theorems for the fibre-like geodesic triangle as 
for the corresponding Euclidean triangles therefore 
the Ceva's and Menelaus' theorems in the geometry $X$ follow from the well-known 
corresponding Euclidean cases. 
\section{$\NIL$ space}
in her dissertatation \cite{B} K.~Brodaczewska studied some aspects of the elementary geometry 
and in \cite{NSV20} the authors discussed the visualization of the $\NIL$ geometry.
In papers \cite{MSz06, PSSz11-1, SchSz11, SchM15} we investigated the equidistant surfaces, 
parallelohedra, crystallography, and an other possible model (so called linear model)
of $\NIL$ geometry. In \cite{I} J.~Inoguchi classified the minimal translation surfaces of the considered space.

In this Section we summarize the relevant notions and notations (see \cite{M97}, \cite{Sz07}).
\subsection{Basic notions of $\NIL$ geometry}
$\NIL$ geometry is a homogeneous 3-space derived from the famous real matrix group $\mathbf{L(R)}$, used by W.~Heisenberg in his electro-magnetic studies.
The Lie theory with the method of projective geometry make possible to describe this topic.

The left (row-column) multiplication of Heisenberg matrices
     \begin{equation}
     \begin{gathered}
     \begin{pmatrix}
         1&x&z \\
         0&1&y \\
         0&0&1 \\
       \end{pmatrix}
       \begin{pmatrix}
         1&a&c \\
         0&1&b \\
         0&0&1 \\
       \end{pmatrix}
       =\begin{pmatrix}
         1&a+x&c+xb+z \\
         0&1&b+y \\
         0&0&1 \\
       \end{pmatrix}
      \end{gathered} \tag{4.1}
     \end{equation}
defines the ``translations" $\mathbf{L}(\mathbf{R})= \{(x,y,z): x,~y,~z\in \mathbf{R} \}$ on the points of
$\NIL= \{(a,b,c):a,~b,~c \in \mathbf{R}\}$.
These translations are not commutative, in general. The matrices $\mathbf{K}(z) \vartriangleleft \mathbf{L}$ of the form
     \begin{equation}
     \begin{gathered}
       \mathbf{K}(z) \ni
       \begin{pmatrix}
         1&0&z \\
         0&1&0 \\
         0&0&1 \\
       \end{pmatrix}
       \mapsto (0,0,z)
      \end{gathered}\tag{4.2}
     \end{equation}
constitute the one parametric centre, i.e. each of its elements commutes with all elements of $\mathbf{L}$.
The elements of $\mathbf{K}$ are called {\it fibre translations}. $\NIL$ geometry of the Heisenberg group can be projectively
(affinely) interpreted by the "right translations"
on points as the matrix formula
     \begin{equation}
     \begin{gathered}
       (1;a,b,c) \to (1;a,b,c)
       \begin{pmatrix}
         1&x&y&z \\
         0&1&0&0 \\
         0&0&1&x \\
         0&0&0&1 \\
       \end{pmatrix}
       =(1;x+a,y+b,z+bx+c)
      \end{gathered} \tag{4.3}
     \end{equation}
shows, according to (4.1). Here we consider $\mathbf{L}$ the projective collineation
group with right actions in homogeneous coordinates.

E. Moln\'ar \cite{M97} derived the well-known infinitesimal arc length square, invariant under translations $\bL$ at any point of $\NIL$ as follows
\begin{equation}
   \begin{gathered}
      (dx)^2+(dy)^2+(-xdy+dz)^2=\\
      (dx)^2+(1+x^2)(dy)^2-2x(dy)(dz)+(dz)^2=:(ds)^2
       \end{gathered} \tag{4.4}
     \end{equation}
Hence we obtain the symmetric metric tensor field $g$ on $\NIL$ with components $g_{ij}$, and its inverse:
\begin{equation}
   \begin{gathered}
       g_{ij}:=
       \begin{pmatrix}
         1&0&0 \\
         0&1+x^2&-x \\
         0&-x&1 \\
         \end{pmatrix},  \quad  g^{ij}:=
       \begin{pmatrix}
         1&0&0 \\
         0&1&x \\
         0&x&1+x^2 \\
         \end{pmatrix} \\
         \text{where} \ \det(g_{ij})=1.
        \end{gathered} \tag{4.5}
     \end{equation}
The translation group $\mathbf{L}$ defined by formula (4.3) can be extended to a larger group 
$\mathbf{G}$ of collineation,
preserving the fibering, that will be equivalent to the (orientation preserving) isometry group of $\NIL$.
In \cite{M06} E.~Moln\'ar has shown that
a rotation by angle $\omega$
about the $z$-axis at the origin, as isometry of $\NIL$, leaves invariant the Riemann
metric everywhere, and is a quadratic mapping in $x,y$ to $z$-image $\overline{z}$ as follows:
     \begin{equation}
     \begin{gathered}
       \br(O,\omega):(1;x,y,z) \to (1;\overline{x},\overline{y},\overline{z}); \\
       \overline{x}=x\cos{\omega}-y\sin{\omega}, \ \ \overline{y}=x\sin{\omega}+y\cos{\omega}, \\
       \overline{z}=z-\frac{1}{2}xy+\frac{1}{4}(x^2-y^2)\sin{2\omega}+\frac{1}{2}xy\cos{2\omega}.
      \end{gathered} \tag{4.6}
     \end{equation}
This rotation formula, however, is conjugate by the quadratic mapping 
     \begin{equation}
     \begin{gathered}
       \mathcal{M}:~x \to x'=x, \ \ y \to y'=y, \ \ z \to z'=z-\frac{1}{2}xy  \ \ \text{to} \\
       (1;x',y',z') \to (1;x',y',z')
       \begin{pmatrix}
         1&0&0&0 \\
         0&\cos{\omega}&\sin{\omega}&0 \\
         0&-\sin{\omega}&\cos{\omega}&0 \\
         0&0&0&1 \\
       \end{pmatrix}
       =(1;x",y",z"), \\
       \text{with} \ \ x" \to \overline{x}=x", \ \ y" \to \overline{y}=y", \ \ z" \to \overline{z}=z"+\frac{1}{2}x"y",
      \end{gathered} \tag{4.7}
     \end{equation}
i.e. to the linear rotation formula. This quadratic conjugacy modifies the $\NIL$ translations 
in (4.3), as well.
This is characterized by the following important classification theorem.
\begin{Theorem}[\cite{M06}]
\begin{enumerate}
\item Any group of $\NIL$ isometries, containing a 3-dimensional translation lattice,
is conjugate by the quadratic mapping in (4.7) to an affine group of the affine (or Euclidean) space $\bA^3=\EUC$
whose projection onto the (x,y) plane is an isometry group of $\bE^2$. Such an affine group preserves a plane
$\to$ point null-polarity.
\item The involutive line reflection about the $y$ axis
     \begin{equation}
     \begin{gathered}
       (1;x,y,z) \to (1;-x,y,-z),
      \end{gathered} \notag
     \end{equation}
preserves the Riemann metric, and its conjugates by the above isometries in {$1$} (those of the identity component)
are also {$\NIL$}-isometries. Orientation reversing $\NIL$-isometries do not exist.
\end{enumerate}
\end{Theorem}
\begin{rmrk}
We obtain a new projective model for $\NIL$ geometry from the projective model, derived from quadratic mapping $\mathcal{M}$.
This is the {\it linearized model of $\NIL$ space} (see \cite{B}, \cite{M06}) 
that seems to be more advantageous for future study. 
But we will continue to use the classical Heisenberg model in this survey.
\end{rmrk}
\subsection{Geodesic curves, spheres and their properties} \label{subsection2}
The geodesic curves of the $\NIL$ geometry are generally defined as having locally minimal arc length between their any two (near enough) points.
The system of equations of the parametrized geodesic curves $g(x(t),y(t),z(t))$ in our model 
can be determined by the
Levy-Civita theory of Riemann geometry.
We can assume, that the starting point of a geodesic curve is the origin because we can 
transform a curve to have
arbitrary starting point by translation;
\begin{equation}
\begin{gathered}
        x(0)=y(0)=z(0)=0; \ \ \dot{x}(0)=c \cos{\alpha}, \ \dot{y}(0)=c \sin{\alpha}, \\ \dot{z}(0)=w; \ - \pi \leq \alpha \leq \pi. \notag
\end{gathered}
\end{equation}
The arc length parameter $s$ is introduced by
\begin{equation}
 s=\sqrt{c^2+w^2} \cdot t, \ \text{where} \ w=\sin{\theta}, \ c=\cos{\theta}, \ -\frac{\pi}{2}\le \theta \le \frac{\pi}{2}, \notag
\end{equation}
i.e. unit velocity can be assumed.

The system of equatios for helix-like geodesic curves (see Fig.~12) $g(x(t),y(t),z(t))$ if $0<|w| <1 $ is:
\begin{equation}
\begin{gathered}
x(t)=\frac{2c}{w} \sin{\frac{wt}{2}}\cos\Big( \frac{wt}{2}+\alpha \Big),\ \
y(t)=\frac{2c}{w} \sin{\frac{wt}{2}}\sin\Big( \frac{wt}{2}+\alpha \Big), \notag \\
z(t)=wt\cdot \Big\{1+\frac{c^2}{2w^2} \Big[ \Big(1-\frac{\sin(2wt+2\alpha)-\sin{2\alpha}}{2wt}\Big)+ \\
+\Big(1-\frac{\sin(2wt)}{wt}\Big)-\Big(1-\frac{\sin(wt+2\alpha)-\sin{2\alpha}}{2wt}\Big)\Big]\Big\} = \\
=wt\cdot \Big\{1+\frac{c^2}{2w^2} \Big[ \Big(1-\frac{\sin(wt)}{wt}\Big)
+\Big(\frac{1-\cos(2wt)}{wt}\Big) \sin(wt+2\alpha)\Big]\Big\}. \tag{4.8}
\end{gathered}
\end{equation}
In the cases with when $w=0$ the geodesic curve is the following:
\begin{equation}
x(t)=c\cdot t \cos{\alpha}, \ \ y(t)=c\cdot t \sin{\alpha}, \ \ z(t)=\frac{1}{2} ~ c^2 \cdot t^2 \cos{\alpha} \sin{\alpha}. \tag{4.9}
\end{equation}
\begin{figure}[ht]
\centering
\includegraphics[width=7cm]{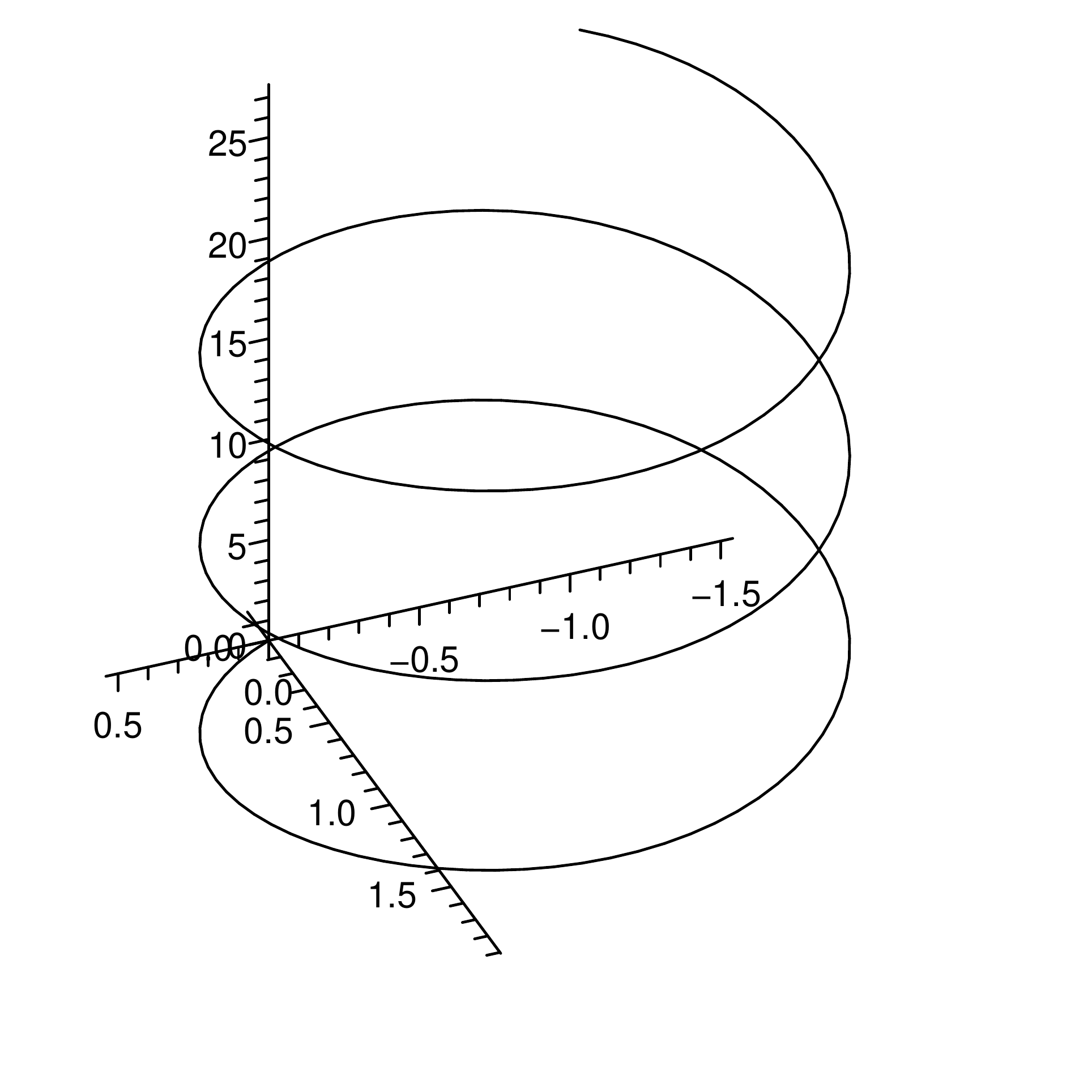}
\caption{Geodesic curve with parameters $\alpha=\frac{\pi}{6}$ and $\beta=\frac{\pi}{4}$.}
\label{}
\end{figure}
The cases $|w|=1$ are trivial: $(x,y)=(0,0), \ z=w \cdot t$.
\begin{definition}
The distance $d(P_1,P_2)$ between the points $P_1$ and $P_2$ is defined by the arc length of geodesic curve
from $P_1$ to $P_2$.
\end{definition}
\begin{definition}
 The geodesic sphere of radius $R$ with centre at the point $P_1$ is defined as the set of all points 
 $P_2$ in the space with the condition $d(P_1,P_2)=R$. Moreover, we require that the geodesic sphere is a simply connected 
 surface without self-intersection 
 in $\NIL$ space (Fig.~13).
 \end{definition}
 \begin{definition}
 The body of the geodesic sphere with centre $P_1$ and radius $R$ in $\NIL$ space is called a 
 geodesic ball, denoted by $B_{P_1}(R)$,
 i.e., $Q \in B_{P_1}(R)$ iff $0 \leq d(P_1,Q) \leq R$.
 \end{definition}
 We proved in \cite{Sz13-2,Sz07} the following important theorems:
 \begin{figure}[ht]
 \centering
 \includegraphics[width=15cm]{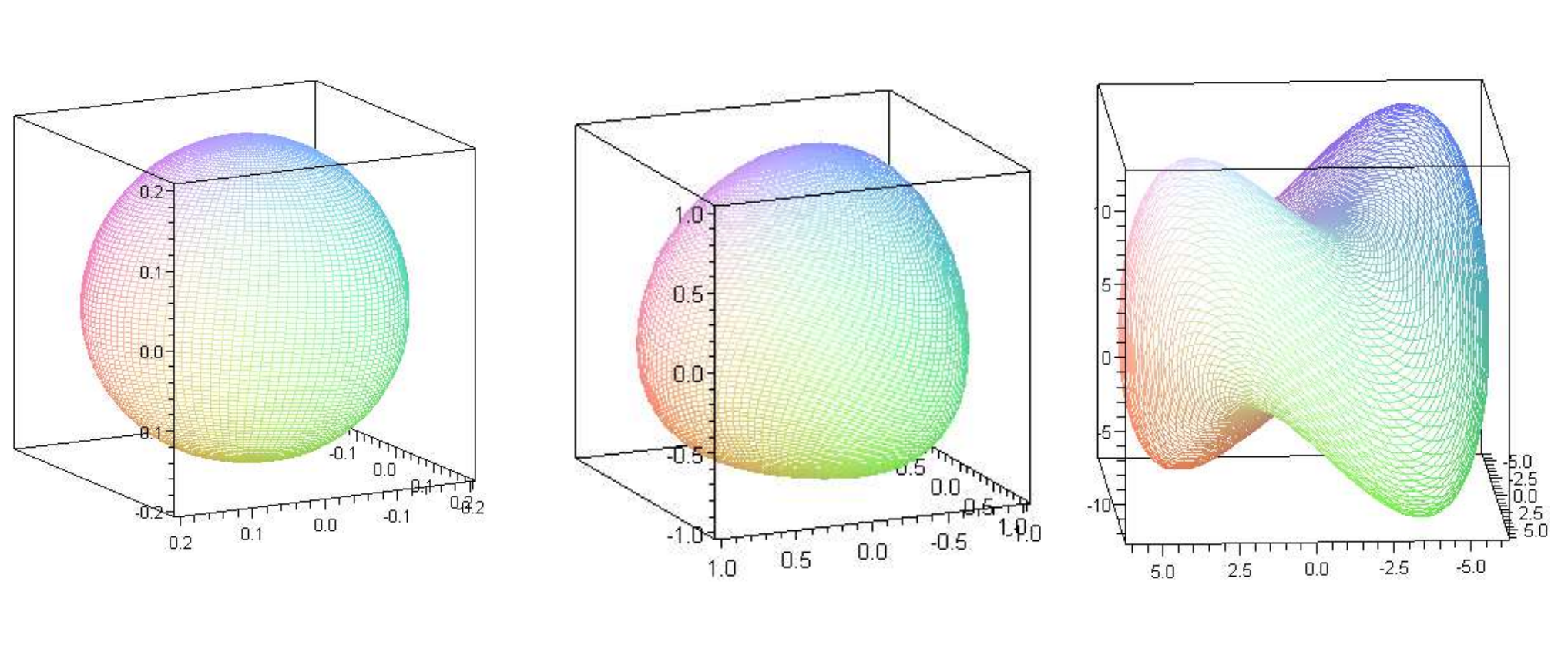}
 \caption{Nil geodesic spheres of radii: $R=0.2$, $R=1$, $R=6$.}
 \label{}
 \end{figure}
 \begin{Theorem}[\cite{Sz07}]
 The geodesic sphere and ball of radius $R$ exists in $\NIL$ space if and only if $R \in [0,2\pi].$
 \end{Theorem}
 \begin{Theorem}[\cite{Sz13-2}]
 The geodesic $\NIL$ ball $B(S(R))$ is convex in the affine-Euclidean sense in our model if and only if $R \in [0,\frac{\pi}{2}]$. 
\end{Theorem}
Next, we recall some important properties of geodesic curves and spheres proved in \cite{Sz22}. 

\begin{enumerate}
\item Consider points $P(x(t),y(t),z(t))$ lying on a sphere $S$ of radius $R$  
centred at the origin. The coordinates of $P$ are given by parameters $(\alpha\in[-\pi, \pi), ~ \theta\in [-\frac{\pi}{2},\frac{\pi}{2}], ~ R>0)$. 

From the equations (4.8) and (4.9) we directly obtain the following
\begin{lemma}[\cite{Sz22}]
\begin{enumerate}
\item
$x(t)^2+y(t)^2=\frac{4c^2}{w^2}\sin^2{\frac{wt}{2}},$
that is if $\theta \ne \pm \frac{\pi}{2}$ and $t=R$ is given and $\alpha\in[-\pi,\pi)$ then the endpoints $P$ of the geodesic curves lie on a cylinder of radius 
$r=\left| \frac{4c}{w}\sin{\frac{wR}{2}}\right|$ with axis $z$. Therefore, we obtain the following connection between parameters $\theta$ and $R$:
\begin{equation}
R=2\cdot \arcsin \left[  \frac{\sqrt{x^2(R)+y^2(R)}}{2\cdot \cot{\theta}}\right] \frac{1}{\sin{\theta}} \tag{4.10}
\end{equation}
\item
If $\theta = \pm \frac{\pi}{2}$ then the endpoints $P(x(R),y(R),z(R))$ of the geodesics $g(x(t),y(t),z(t))$ lie on the 
$z$-axis thus their orthogonal projections onto the $[x,y]$-plane is the origin and $x(R)=y(R)=0$, $z(R)=d(O,R)=R$.
\item
Moreover, the cross section 
of the spheres $S$ with the plane $[x,z]$ is given by the following system of equation:
     \begin{equation}
     \begin{gathered}
    X(R,\theta)=\frac{2c}{w} \sin{\frac{wR}{2}}=\frac{2\cos{\theta}}{\sin{\theta}} \sin{\frac{R \sin{\theta}}{2}}, \\ 
    Z(R,\theta)=wR+\frac{c^2R}{2w} - \frac{c^2}{2w^2}\sin{wR}= \\ 
    R\sin{\theta}+\frac{R\cos^2{\theta}}{2\sin{\theta}} - \frac{\cos^2{\theta}}{2\sin^2{\theta}}\sin(R\sin{\theta}), \ \ (\theta\in [-\frac{\pi}{2},\frac{\pi}{2}]\setminus\{0\}); \\
    \text{if} \ \theta=0 \ \text{then} \ X(R,0)=R, \ Z(R,0)=0. \tag{4.11}
   \end{gathered}
   \end{equation}
\end{enumerate}
\end{lemma}
\item 
In \cite{Sz22} we introduced  the usual notion of the fibre projection $\mathcal{P}$, a 
projection parallel to fibre lines (parallel to $z$-axis) onto the $[x,y]$ plane.
The image of a point $P$ is the intersection with the $[x,y]$ base plane of the line parallel to fibre line passing through $P$, $\mathcal{P}(P)=P^*$.

Analysing the parametric equations of the geodesic curves $g(x(t),y(t),z(t))$ 
with starting points at the origin we gound the following
\begin{lemma}[\cite{Sz22}]
If $0 <|w| <1$ for geodesic curve $g(x(t),y(t),z(t))$ $(t\in [0,R])$ then 
the fibre projection $\mathcal{P}$ of the geodesic curves onto the $[x, y]$ plane is an 
Euclidean circular arc that is contained in a circle with equation 
\begin{equation}
\Big(x(t)+\frac{c}{w}\sin{\alpha}\Big)^2+\Big(y(t)-\frac{c}{w}\cos{\alpha}\Big)^2=\Big(\frac{c}{w}\Big)^2=\cot^2{\theta}. \tag{4.12}
\end{equation}
If $w=0$ then fibre projection $\mathcal{P}$ of the geodesic curves $g(x(t),y(t),z(t))$ $(t\in [0,R])$ onto the $[x, y]$ plane is a segment with starting point at the origin
where it is contained by the straight line with equation
\begin{equation}
y=\tan{\alpha}\cdot x.\tag{4.13}
\end{equation}
If $w=1$ then the fibre projection $\mathcal{P}$ of the geodetic curves $g(x(t),y(t),z(t))$ $(t\in [0,R])$ onto 
the $[x, y]$ plane is the origin.
\end{lemma}
From the equation (4.12) we directly have the following
\begin{corollary}[\cite{Sz22}]
\begin{enumerate}
\item If we know the equation of the circle that contains the orthogonal projected image $OP^*$ of a geodesic curve segment 
$g_{OP}=g(x(t),$ $y(t),z(t))$ $(t\in [0,R])$ onto the $[x,y]$ plane where $0 <|w| <1$ is a known real number and the 
coordinates of $P^*=(x(R), y(R),0)$ then the parametric equation of the geodesic curve segment 
$g_{OP}$ is uniquely determined. 
That means that there is a {\it one-to-one correspondence between the circular arcs $OP^*$ 
and the geodesic curve segments $OP$} in the above sense. 
\item If $w=0$ then the fibre projection $\mathcal{P}$ of the geodetic curve is a segment 
with starting point at the origin and
it is contained in the straight line $y=\tan{\alpha}\cdot x$, therefore in this situation there is
a {\it one-to-one correspondence between the projected image $OP^*$ and the geodesic 
curve segments $OP$}. 
\item If $w=1$ then the fibre projection $\mathcal{P}$ of the geodetic curve 
is the origin so here it is also a {\it one-to-one correspondence between the projected image and the above geodesic curves.}
\end{enumerate}
\end{corollary}
\end{enumerate}
\subsection{Geodesic triangles and their interior angle sums} \label{section3}
Similarly to the $\SXR$ and $\HXR$ geometries in subsection 3.4 or more generally in Riemannian geometries the angle $\theta$ of two intersecting curves can be determined by
the metric tensor of the considered geometry (see (4.5)) using the formula (3.9).
If their tangent vectors at their common point are $\bu$ and $\bv$ and $g_{ij}$ are the components of the metric tensor. 

It is clear by the above definition of the angles and by the metric tensor (4.5), that
the angles are the same as the Euclidean angles at the origin by a pull back translation. 

We note here that the angle of two intersecting geodesic curves depend on the orientation of the tangent vectors. We will consider
the {\it interior angles} of the triangles that are denoted at the vertex $A_i$ by $\omega_i$ $(i\in\{1,2,3\})$.

A geodesic triangle is called fibre-like if one of its edges lies on a fibre line. In this section we study the right-angled fibre-like
triangles. We can assume without loss of generality
that the vertices $A_1$, $A_2$, $A_3$ of a fibre-like right-angled triangle (see Fig.~14.a-b) have the following coordinates:
\begin{equation}
A_1=(1,0,0,0),~A_2=(1,0,0,z^2),~A_3=(1,x^3,0,z^3=z^2) \tag{4.14}
\end{equation}
\begin{figure}[ht]
\centering
\includegraphics[width=13cm]{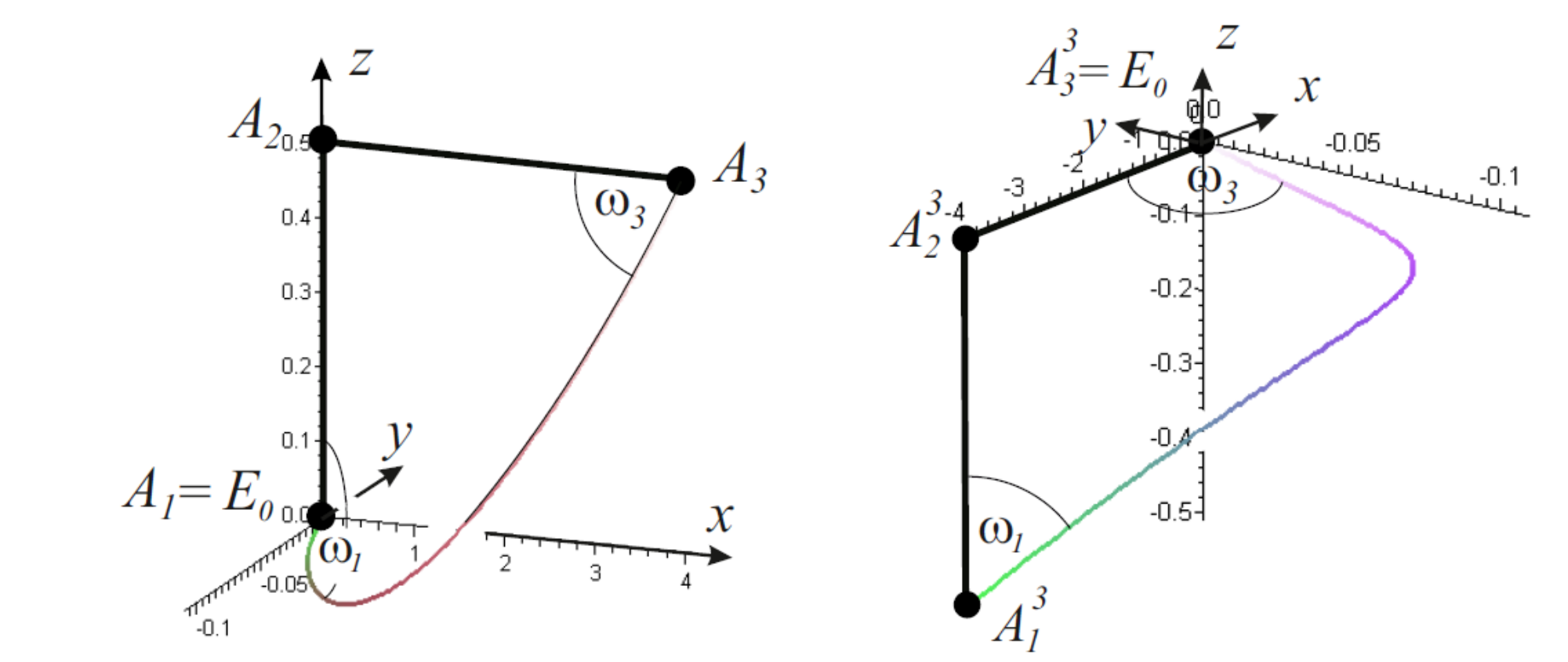}

a. \hspace{6cm} b.
\caption{a.~Fibre-like geodesic triangle $A_1A_2A_3$, where $A_1=(1,0,0,0)$,~$A_2=(1,0,0,\frac{1}{2})$,~$A_3=(1,4,0,\frac{1}{2})$. ~b. Its translated image
$A_1^3A_2^3A_3^3$ where $A_1^3=(1,-4,0,-\frac{1}{2})$,~$A_2^3=(1,-4,0,0)$,~$A_3^3=(1,0,0,0)$ (see \cite{Sz16}).}
\label{}
\end{figure}

{\it In order to determine the interior angles,
we defined \emph{translations} $\bT_{A_i}$, $(i\in \{2,3\})$ as elements of the isometry group of $\NIL$, that
maps the origin $E_0$ onto $A_i$ (see Fig.~14,~15).}
Our aim is to determine angle sum $\sum_{i=1}^3(\omega_i)$ of the interior angles of the above right-angled
fibre-like geodesic triangle $A_1A_2A_3$.
We have seen that $\omega_2=\frac{\pi}{2}$ and the angle of geodesic curves with a common point at the origin $E_0$ is the same as the
Euclidean one. Therefore it can be determined in the usual Euclidean sense. 
Hence, $\omega_1$ is equal to the angle $\angle (g(E_0, A_3),g(E_0, A_2))$ 
where $g(E_0, A_3)$ and $g(E_0, A_2)$ are oriented geodesic curves.
Moreover, the translation $\bT_{A_3}$ is an isometry
in $\NIL$ geometry thus
$\omega_3$ is equal to the angle $\angle (g(A_3^3, A_1^3),g(A_3^3, A_2^3)) $ 
where $g(A_3^3, A_1^3)$ and $g(A_3^3, A_2^3)$ are also oriented geodesic curves $(E_0=A_3^3)$.

We denote the oriented unit tangent vectors of the geodesic curves $g(E_0, A_i^j)$ with $\mathbf{t}_i^j$ where
$(i,j)\in \{(1,3),(2,3),(3,0),(2,0)\}$ and $A_3^0=A_3$, $A_2^0=A_2$.
The Euclidean coordinates of $\mathbf{t}_i^j$ are :
\begin{equation}
\mathbf{t}_i^j=(\cos(\theta_i^j) \cos(\alpha_i^j), \cos(\theta_i^j) \sin(\alpha_i^j), \sin(\theta_i^j)). \tag{4.15}
\end{equation}
\begin{lemma}[\cite{Sz16}]
The sum of the interior angles of a fibre-like right-angled geodesic triangle is greater or equal to $\pi$.
\end{lemma}
\begin{conj}[\cite{Sz16}]
The sum of the interior angles of any fibre-like geodesic triangle is greater or equal to $\pi$.
\end{conj}
We fix the coordinates the $z^2=z^3\in \mathbf{R}$ of $A_2$ and $A_3$ and study the interior angle sum
$\sum_{i=1}^3(\omega_i(x^3))$ of the right-angled
geodesic triangle $A_1A_2A_3$ if the $x^3$ coordinate of $A_3$ tends to zero or infinity.
E.g. $\lim_{x^3 \to 0}(\omega_1(x^3))=0$ 
because the geodesic line $g(E_0,A_3)$ tends to the geodesic line $g(E_0,A_2)$ therefore their angle $\omega_1$ tends to the zero, 
and $\omega_3$ tends to $\frac{\pi}{2}$ (see Fig.~14). Similarly to this from the system of equations 
(4.8) gives the following results
\begin{lemma}[\cite{Sz16}]
If coordinates $z^2=z^3 \in \mathbf{R}$ are fixed then
$$\lim_{x^3 \to 0}(\omega_1(x^3))=0, ~
\lim_{x^3 \to 0}(\omega_3(x^3))=\frac{\pi}{2} ~ \Rightarrow ~ \lim_{x^3 \to 0} \Bigg(\sum_{i=1}^3(\omega_i(x^3))\Bigg)=\pi,$$
$$\lim_{x^3 \to \infty}(\omega_1(x^3))=\frac{\pi}{2}, ~
\lim_{x^3 \to \infty}(\omega_3(x^3))=0 ~ \Rightarrow ~ \lim_{x^3 \to \infty} \Bigg(\sum_{i=1}^3(\omega_i(x^3))\Bigg)=\pi.$$
\end{lemma}
\subsubsection{Hyperbolic-like right angled geodesic triangles}
A geodesic triangle is hyperbolic-like if its vertices lie in the base plane of the model.
In this section we recall the results of \cite{Sz16} about the interior angle sum of right-angled hyperbolic-like
triangles. We can assume without loss of generality
that the veritices $A_1$, $A_2$, $A_3$ of a hyperbolic-like right-angled triangle (see Fig.~15) $T_g$ have the following coordinates:
\begin{equation}
A_1=(1,0,0,0),~A_2=(1,0,y^2,0),~A_3=(1,x^3,y^2=y^3,0) \tag{4.15}
\end{equation}
\begin{figure}[ht]
\centering
\includegraphics[width=12cm]{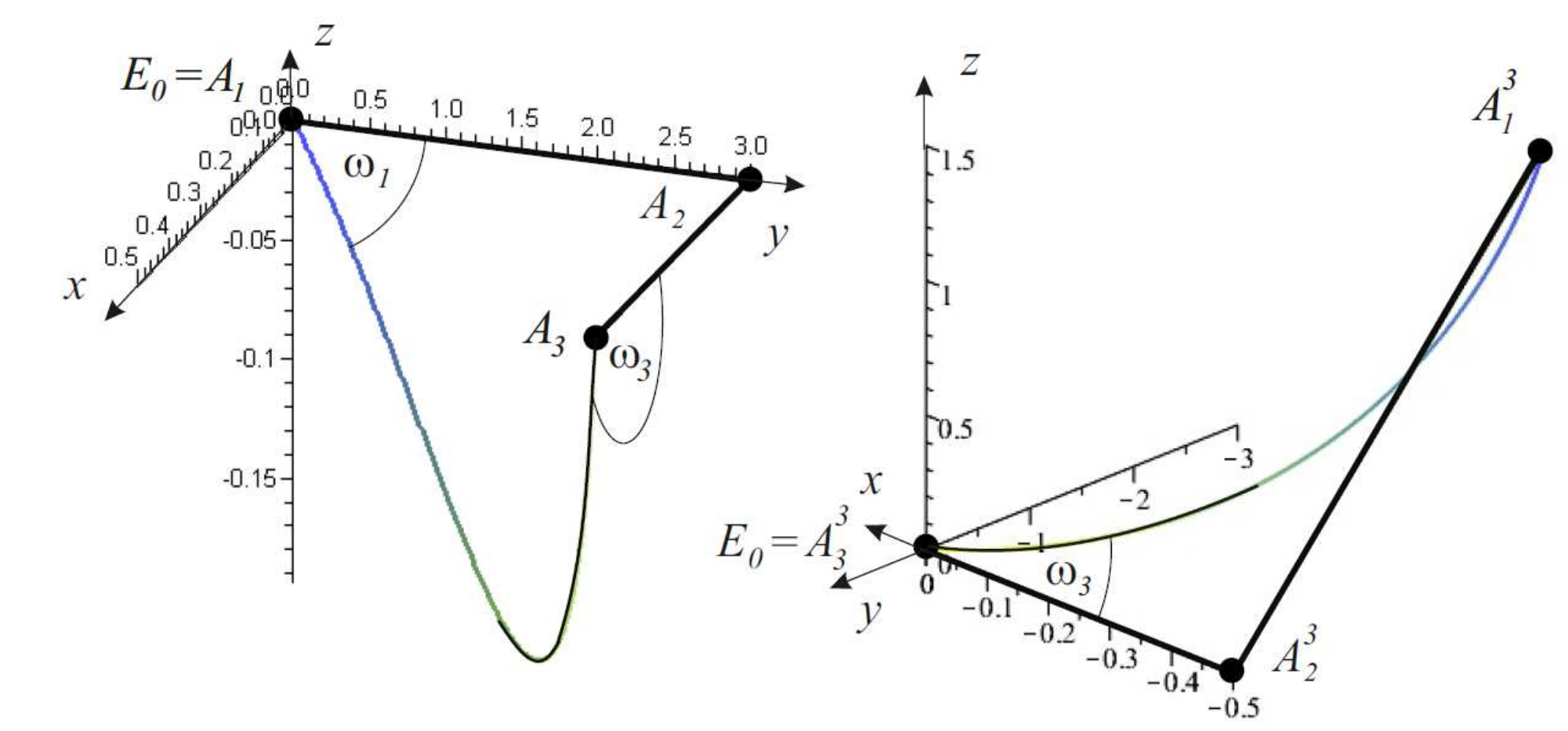}

a.~\hspace{5cm} ~b.
\caption{Hyperbolic-like geodesic triangle $A_1A_2A_3$, where $A_1=(1,0,0,0)$,~$A_2=(1,0,3,0)$,~$A_3=(1,\frac{1}{2},3,0)$.  ~b. Its translated image
$A_1^3A_2^3A_3^3$ where $A_1^3=(1,-\frac{1}{2},-3,\frac{3}{2})$,~$A_2^3=(1,-\frac{1}{2},0,0)$,~$A_3^3=(1,0,0,0)$ (see \cite{Sz16})}
\label{}
\end{figure}
First we fix the $x^3\in \mathbf{R}$ coordinate of the vertex $A_3$ and study the the interior angle sum
$\sum_{i=1}^3(\omega_i(y^2=y^3))$ of the right-angled
geodesic triangle $A_1A_2A_3$ if $y^2=y^3$ the coordinates of vertices $A_2$ and $A_3$ tend to zero or infinity.
From the system of equations (4.8) we obtain the following results
\begin{lemma}[\cite{Sz16}]
If the coordinate $x^3 \in \mathbf{R}$ is fixed then
$$\lim_{y^2=y^3 \to 0}(\omega_1(y^2))=\frac{\pi}{2}, ~
\lim_{y^2=y^3 \to 0}(\omega_3(y^2))=0 ~ \Rightarrow ~ \lim_{y^2=y^3 \to 0} \Bigg(\sum_{i=1}^3(\omega_i(y^2))\Bigg)=\pi,$$
$$\lim_{y^2=y^3 \to \infty}(\omega_1(y^2))=0, ~
\lim_{y^2=y^3 \to \infty}(\omega_3(y^2))=\frac{\pi}{2}~ \Rightarrow ~ \lim_{y^2=y^3 \to \infty} \Bigg(\sum_{i=1}^3(\omega_i(y^2))\Bigg)=\pi.$$
\end{lemma}
Secondly we fix the $y^2=y^3\in \mathbf{R}$ coordinates of the vertices $A_2$ and $A_3$ and study the internal angle sum
$\sum_{i=1}^3(\omega_i(x^3))$ of the right-angled
geodesic triangle $A_1A_2A_3$ if $x^3$ coordinate of vertex $A_3$ tends to zero or infinity.
From the system of equations (4.8) we obtain the following
\begin{lemma}[\cite{Sz16}]
If the coordinates $y^2=y^3 \in \mathbf{R}$ are fixed then
$$\lim_{x^3 \to 0}(\omega_1(x^3))=0, ~
\lim_{x^3 \to 0}(\omega_3(x^3))=\frac{\pi}{2} ~\Rightarrow ~ \lim_{x^3 \to 0} \Bigg(\sum_{i=1}^3(\omega_i(x^3))\Bigg)=\pi,$$
$$\lim_{x^3 \to \infty}(\omega_1(x^3))=\frac{\pi}{2}, ~
\lim_{x^3 \to \infty}(\omega_3(x^3))=0 ~ \Rightarrow ~ \lim_{x^3 \to \infty} \Bigg(\sum_{i=1}^3(\omega_i(x^3))\Bigg)=\pi.$$
\end{lemma}
We can determine the interior angle sum of arbitrary hyperbolic-like geodesic triangle similarly 
as in the fibre-like case.
\medbreak
Finally, we have the following 
\begin{lemma}[\cite{Sz16}]
The interior angle sums of hyperbolic-like right-angled geodesic triangles are less than or equal to $\pi$.
\end{lemma}
\begin{conj}[\cite{Sz16}]
The sum of the interior angles of any hyperbolic-like geodesic triangle is less than or equal to $\pi$.
\end{conj}
\subsubsection{Geodesic triangles with interior angle sum $\pi$}
In the above sections we discussed the fibre- and hyperbolic-like geodesic triangles and proved that there are right-angled
geodesic triangles whose angle sum $\sum_{i=1}^{3}(\omega_i)$ is greater, less than or equal to $\pi$, but $\sum_{i=1}^{3}(\omega_i)=\pi$ is
realized if one of the vertices of a geodesic triangle $A_1A_2A_3$ tends to the infinity. 
In \cite{Sz16} we proved in \cite{Sz16} the following
\begin{lemma}[\cite{Sz16}]
There exists a geodesic triangle $A_1A_2A_3$ with interior angle sum $\pi$ where its vertices are {\it proper}
(i.e. $A_i \in \NIL$, vertices are not at infinity, $i \in \{1,2,3\}$).
\end{lemma}
We summarize the lemmas of this Section as follows
\begin{Theorem}[\cite{Sz16}]
The sum of the interior angles of a geodesic triangle of $\NIL$ space can be greater, less or equal to $\pi$.
\end{Theorem}

\subsection{On Menelaus' and Ceva's theorems in $\NIL$ space} 
As in previous $\SXR$ and $\HXR$ spaces, the question arises as to what the surface of a geodetic triangle will be and which elementary theorems may be true in this geometry. 
In the paper \cite{Sz22}, as in the previously discussed spaces, we introduced the concept of the surface of $\NIL$ geodetic triangles using Apollonius surfaces 
and examined the theorems of Menelaus and Ceva.
To discuss on Menelaus' and Ceva's theorems, 
we had to define what we consider to be a {\it line} on the surface of a geodesic triangle and the definition of a simple ratio.

Let $\mathcal{S}_{A_0A_1A_2}$ be the surface of the geodesic triangle 
$A_0A_1A_2$ and $P_1$, $P_2 \in \mathcal{S}_{A_0A_1A_2}$ be any two points. 
Natural requirements for a {\it line} passing through points $P_1$ and $P_2$ 
on $\mathcal{S}_{A_0A_1A_2}$ are that:
\begin{enumerate}
\item Two surface points uniquely determine a {\it line} (connecting curve) $\mathcal{G}^{\mathcal{S}_{A_0A_1A_2}}_{P_1P_2}$.
\item Any two points on a surface {\it line} $\mathcal{G}^{\mathcal{S}_{A_0A_1A_2}}_{P_1P_2}$ 
define the same line.
\item The surface {\it line} determined by two points of a geodesic curve lying on 
the surface $\mathcal{S}_{A_0A_1A_2}$ coincides with the geodesic curve.   
\end{enumerate}
\begin{rmrk}
An obvious option for definition of a line (connecting curve) $\mathcal{G}^{\mathcal{S}_{A_0A_1A_2}}_{P_1P_2}$ would be the fibre projection of the geodesic curve $g_{P_1P_2}$ 
into the surface $\mathcal{S}_{A_0A_1A_2}$ but it is clear, that this definition does not satisfy requirement 2. 
\end{rmrk}
We consider a {\it geodesic triangle $A_0A_1A_2$} in the projective model of $\NIL$ space 
(see Subsection 4.1).
Without loss of generality, we can assume that $A_0=(1,0,0,0)$. 
The geodesic lines that contain the sides $A_0A_1$ and $A_0A_2$ of the given 
triangle can be characterized directly 
by the corresponding parameters $\theta_i$ and $\alpha_i$ $(i=1,2)$
(see (4.8) and (4.9)).
The geodesic curve including the side segment $A_1A_2$ is also determined by one of 
its endpoints and its parameters. 
In order to determine the corresponding parameters of this 
geodesic line we use {\it for example a $\NIL$ translation} $\bT(A_1)$, 
as elements of the isometry group of $\NIL$ geometry, that
maps $A_1=(1,x_1,y_1,z_1)$ onto $A_0=(1,0,0,0)$ (up to a positive determinant factor). 
\begin{rmrk} 
By the results of Theorem 4.4, we may assume that the surface $\mathcal{S}_{A_0A_1A_2}$ of the geodesic triangle 
$A_0A_1A_2$ is contained in a geodesic 
$\NIL$ sphere of radius $\pi$. 
\end{rmrk}
We generalized the notion of simple ratio to the point triples lying on geodesic lines of $\NIL$ space:
\begin{definition}
Let $A$, $B$, and $P$ be distinct points on a geodesic curve in $\NIL$ space then 
their simple ratio is
$$s^N(A,P,B) =  {d(A,P)}/{d(P,B)}$$ if $P$ is between $A$ and $B$, and
$$s^N(A,P,B) = -{d(A,P)}/{d(P,B)}$$ 
otherwise, where
$d$ is the distance function of $\NIL$ geometry.   
\end{definition}
Let $A$, $B$ and $P$ be distinct points on a non-fibre-like geodetic curve in the $\NIL$ and let 
$A^*$, $B^*$ and $P^*$ be
their projected images by $\mathcal{P}$. 
\begin{lemma}[\cite{Sz22}]
The Euclidean length $\mathcal{C}(A^*,P^*)$ of circle arc or line 
segment $\overset{\LARGE\frown}{A^*P^*}$ satisfies the following equations 
\begin{equation}
\begin{gathered}
\mathcal{C}(A^*,P^*)=d(A,P)\cdot \cos{\theta},\ \ \ \mathcal{C}(P^*,B^*)=d(P,B)\cdot \cos{\theta}.
\end{gathered} \tag{4.16}
\end{equation}
{Therefore, the projection $\mathcal{P}$  preserves the ratio of lengths in the above sense.}
\end{lemma}

Lemma 4.7, Corollary 4.8 and the above projection $\mathcal{P}$ were used 
to define the surface line but the definition is technical due to the complex structure of the 
geometry and therefore not detailed here (see \cite{Sz22}). 
The main results are summarized as follows 
\begin{corollary}[\cite{Sz22}]
Menelaus' theorem does not hold in $\NIL$ geometry. 
However, as can be seen above, the Menelaus' condition plays an important role 
in defining {\it lines} on the surface of a given triangle.
\end{corollary}
Using the above Menelaus' condition, similar to the Euclidean proof, we obtain the $\NIL$ Cava's theorem:
\begin{Theorem}[\cite{Sz22}]
If $T$ is a point which does not lie on any side of a geodesic triangle 
$A_0A_1A_2$ in $\NIL$ space such that
the curves $\mathcal{G}^{\mathcal{S}_{A_0A_1A_2}}_{A_0T}$ and $g_{A_1A_2}$ intersect at $P_{12}$, 
$\mathcal{G}^{\mathcal{S}_{A_0A_1A_2}}_{A_1T}$ and $g_{A_0A_2}$ at $P_{02}$, 
and $\mathcal{G}^{\mathcal{S}_{A_0A_1A_2}}_{A_2T}$ and $g_{A_0A_1}$ at $P_{01}$, 
then $$s^N(A_0,P_{01},A_1)s^N(A_1,P_{12},A_2)s^N(A_2,P_{02},A_0) = 1.$$ 
\end{Theorem}
Using the Lemma 4.8 follows that the corresponding Ceva's theorem is also true for the 
{\it projected configuration} i.e. for the triangle $A_0^*A_1^*A_2^*$ and the points $T^*$, $P_{01}^*$, $P_{12}^*$, $P_{02}^*$. 
\begin{Theorem}[\cite{Sz22}]
If $T^*$ is a point not on any side of {\it circle arc} triangle (the projected image of a geodesic triangle in general type) $A_0^*A_1^*A_2^*$ in the base plane of the $\NIL$ space such that
the arcs (or line segments)  $\overset{\LARGE\frown}{A_0^*T^*}$ and $\overset{\LARGE\frown}{A_1^*A_2^*}$ meet in $P_{12}^*$, $\overset{\LARGE\frown}{A_1^*T^*}$ and 
$\overset{\LARGE\frown}{A_0^*A_2^*}$ in $P_{02}^*$, and $\overset{\LARGE\frown}{A_2^*T^*}$ and $\overset{\LARGE\frown}{A_0^*A_1^*}$ in $P_{01}^*$, 
then $$s^c(A_0^*,P_{01}^*,A_1^*)s^c(A_1^*,P_{12}^*,A_2^*)s^c(A_2^*,P_{02}^*,A_0^*) = 1.$$ 
\end{Theorem}
\begin{rmrk}
Using the previous notions and theorems, as in the Euclidean case, we can define, 
for example, the circumscribed 
circle of a geodesic triangle and its centre, the centroid of a geodesic triangle 
as the point where the three medians of the triangle meet. 
A median of a geodesic triangle $A_0A_1A_2$ in the $\NIL$ space is a surface {\it line} 
$\mathcal{S}_{A_0A_1A_2})$ from one vertex to the mid point on the 
opposite side of the triangle contained in. We will examine these in a forthcoming paper. 
\end{rmrk}
\section{$\SLR$ geometry}
The basic concepts of the model of $\SLR$ geometry can be found in \cite{M97}. 

In \cite{DESS09} the authors considered the geodesics and geodesic spheres that 
gave exact solutions of ODE system that describes geodesics. 
Moreover, geodesic spheres are determined and a visualization of $\SLR$ geometry is also given. 
In \cite{EH14} Z. Erjavec and D. Horvat investigated and characterized the 
non-geodesic biharmonic curves and proved the statement that 
only proper biharmonic curves are helices. Also, the explicit parametric equations 
of proper biharmonic helices were found.
In \cite{E15} the author derived the equation of minimal surface and gave fundamental 
examples of minimal surfaces.
In \cite{E19} Z. Erjavec discussed the so-called Killing magnetic curves. 
In \cite{NSV20} the authors discussed visualization methods of the considered geometry.

In \cite{CsSz16}, we studied the sum of the interior angles of the geodesic and translation 
triangles (see Subsection 5.3).

\subsection{Basic notions of the $\SLR$ geometry}
In this section we summarize the 
the real $ 2\times 2$ matrices $\begin{pmatrix}
         d&b \\
         c&a \\
         \end{pmatrix}$ with unit determinant $ad-bc=1$ which
constitute a Lie transformation group by the usual product operation, taken to act on row matrices as on point coordinates on the right as follows
\begin{equation}
\begin{gathered}
(z^0,z^1)\begin{pmatrix}
         d&b \\
         c&a \\
         \end{pmatrix}=(z^0d+z^1c, z^0 b+z^1a)=(w^0,w^1) \\
\mathrm{with} \ w=\frac{w^1}{w^0}=\frac{b+\frac{z^1}{z^0}a}{d+\frac{z^1}{z^0}c}=\frac{b+za}{d+zc} \tag{5.1}
\end{gathered}
\end{equation}
as a right action on the complex projective line $\bC^\infty$.
This group is a $3$-dimensional manifold, because of its $3$ independent real coordinates and with its usual neighbourhood topology
\cite{S}, \cite{T}.
In order to model the above structure in the projective sphere $\cP \cS^3$ and in the projective space $\cP^3$ (see \cite{M97}),
we introduce the new projective coordinates $(x^0,x^1,x^2,x^3)$ where
\begin{equation}
a:=x^0+x^3, \ b:=x^1+x^2, \ c:=-x^1+x^2, \ d:=x^0-x^3. \notag
\end{equation}
Meanwhile we turn
to the proportionality $\mathbf{SL_{\mathrm{2}}R} < \mathbf{PSL_{\mathrm{2}}R}$, natural in this context.
Then it follows that
\begin{equation}
0>bc-ad=-x^0x^0-x^1x^1+x^2x^2+x^3x^3 \tag{5.2}
\end{equation}
describes the interior of the above one-sheeted hyperboloid solid $\cH$ in the usual Euclidean coordinate simplex, with the origin
$E_0(1,0,0,0)$ and the ideal points of the axes $E_1^\infty(0,1,0,0)$, $E_2^\infty(0,0,1,0)$, 
$E_3^\infty(0,0,0,1)$.
We consider the collineation group ${\bf G}_*$ that acts on the projective sphere 
$\cS \cP^3$  and preserves a polarity, i.e. a scalar product of signature
$(- - + +)$, this group leaves the one sheeted hyperboloid solid $\cH$ invariant.
We have to choose an appropriate subgroup $\mathbf{G}$ of $\mathbf{G}_*$ as the isometry group, then the universal covering group and space
$\widetilde{\cH}$ of $\cH$ will be the hyperboloid model of $\SLR$ (\cite{M97}).

The specific isometries $\bS(\phi)$ $(\phi \in \bR )$ constitute a one parameter group given by the matrices
\begin{equation}
\begin{gathered} \bS(\phi):~(s_i^j(\phi))=
\begin{pmatrix}
\cos{\phi}&\sin{\phi}&0&0 \\
-\sin{\phi}&\cos{\phi}&0&0 \\
0&0&\cos{\phi}&-\sin{\phi} \\
0&0&\sin{\phi}&\cos{\phi}
\end{pmatrix}.
\end{gathered} \tag{5.3}
\end{equation}
The elements of $\bS(\phi)$ are the so-called {\it fibre translations}. We obtain a unique fibre line to each $X(x^0,x^1,x^2,x^3) \in \widetilde{\cH}$
as the orbit by right action of $\bS(\phi)$ on $X$. The coordinates of points lying on the fibre line through $X$ can be expressed
as the images of $X$ by $\bS(\phi)$:
\begin{equation}
\begin{gathered}
(x^0,x^1,x^2,x^3) \stackrel{\bS(\phi)}{\longrightarrow} {(x^0 \cos{\phi}-x^1 \sin{\phi}, x^0 \sin{\phi} + x^1 \cos{\phi},} \\ {x^2 \cos{\phi} + x^3 \sin{\phi},-x^2 \sin{\phi}+
x^3 \cos{\phi})}
\end{gathered} \tag{5.4}
\end{equation}
for the Euclidean coordinate set $x:=\frac{x^1}{x^0}$, $y:=\frac{x^2}{x^0}$, $z:=\frac{x^3}{x^0}$, $x^0\ne 0$.

In (5.3) and (5.4) we can see the $2\pi$ periodicity of $\phi$. Moreover, we see the (logical) extension to $\phi \in \bR$, as real parameter, to have
the universal covers $\widetilde{\cH}$ and $\SLR$, respectively, through the projective sphere $\cP\cS^3$. The elements of the isometry group of
$\mathbf{SL_{\mathrm{2}}R}$ (and so by the above extension the isometries of $\SLR$) can be described by the matrix $(a_i^j)$ (see \cite{M97} and \cite{MSz})
\begin{equation}
\begin{gathered} (a_i^j)=
\begin{pmatrix}
a_0^0&a_0^1&a_0^2&a_0^3 \\
\mp a_0^1&\pm a_0^0&\pm a_0^3&\mp a_0^2 \\
a_2^0&a_2^1&a_2^2&a_2^3 \\
\pm a_2^1&\mp a_2^0&\mp a_2^3&\pm a_2^2 \\
\end{pmatrix} \ \ \text{where} \\
-(a_0^0)^2-(a_0^1)^2+(a_0^2)^2+(a_0^3)^2=-1, \ \ -(a_2^0)^2-(a_2^1)^2+(a_2^2)^2+(a_2^3)^2=1, \\
-a_0^0a_2^0-a_0^1a_2^1+a_0^2a_2^2+a_0^3a_2^3=0=-a_0^0a_2^1+a_0^1a_2^0-a_0^2a_2^3+a_0^3a_2^2,
\end{gathered} \tag{5.5}
\end{equation}
and we allow positive proportionality, as the projective freedom.
The horizontal intersection of the hyperboloid solid $\cH$  with the plane $E_0 E_2^\infty E_3^\infty$
provides the {\it base plane} of the model $\widetilde{\cH}=\SLR$.
The fibre through $X$ intersects the hyperbolic $(\mathbf{H}^2)$ 
base plane $z^1=x=0$ at the foot point
\begin{equation}
\begin{gathered}
Z(z^0=x^0 x^0+x^1x^1; z^1=0; z^2=x^0x^2-x^1x^3;z^3=x^0x^3+x^1x^2).
\end{gathered} \tag{5.6}
\end{equation}
We introduce a so-called hyperboloid parametrization as in \cite{M97} as follows
\begin{equation}
\begin{gathered}
x^0=\cosh{r} \cos{\phi}, \\
x^1=\cosh{r} \sin{\phi}, \\
x^2=\sinh{r} \cos{(\theta-\phi)}, \\
x^3=\sinh{r} \sin{(\theta-\phi)},
\end{gathered} \tag{5.7}
\end{equation}
where $(r,\theta)$ are the polar coordinates of the $\mathbf{H}^2$ base plane, and $\phi$ is the fibre coordinate. We note that
$$-x^0x^0-x^1x^1+x^2x^2+x^3x^3=-\cosh^2{r}+\sinh^2{r}=-1<0.$$
The inhomogeneous coordinates will play an important role in the later 
$\EUC$-visualization e.g. of the prism tilings in $\SLR$, and
are given by
\begin{equation}
\begin{gathered}
x=\frac{x^1}{x^0}=\tan{\phi}, \\
y=\frac{x^2}{x^0}=\tanh{r} \frac{\cos{(\theta-\phi)}}{\cos{\phi}}, \\
z=\frac{x^3}{x^0}=\tanh{r} \frac{\sin{(\theta-\phi)}}{\cos{\phi}}.
\end{gathered} \notag
\end{equation}
\subsection{Distances and spheres}
The infinitesimal arc length square can be derived by the standard pull back method.
By the $T^{-1}$-action of (5.6) on the differentials $(\mathrm{d}x^0;\mathrm{d}x^1;\mathrm{d}x^2;\mathrm{d}x^3)$, we obtain
that in this parametrization
the infinitesimal arc length square
at any point of $\SLR$ is the following:
\begin{equation}
   \begin{gathered}
      (\mathrm{d}s)^2=(\mathrm{d}r)^2+\cosh^2{r} \sinh^2{r}(\mathrm{d}\theta)^2+\big[(\mathrm{d}\phi)+\sinh^2{r}(\mathrm{d}\theta)\big]^2.
       \end{gathered} \tag{5.8}
     \end{equation}
Hence we get the symmetric metric tensor field $g_{ij}$ on $\SLR$ by components:
     \begin{equation}
       g_{ij}:=
       \begin{pmatrix}
         1&0&0 \\
         0&\sinh^2{r}(\sinh^2{r}+\cosh^2{r})& \sinh^2{r} \\
         0&\sinh^2{r}&1 \\
         \end{pmatrix}, \tag{5.9}
     \end{equation}
and
\begin{equation}
     \mathrm{d}V=\sqrt{\det(g_{ij})}~dr ~\mathrm{d} \theta ~\mathrm{d} \phi= \frac{1}{2}\sinh(2r) \mathrm{d}r ~\mathrm{d} \theta~ \mathrm{d} \phi
     \notag
     \end{equation}
as the volume element in the hyperboloid coordinates.
The geodesic curves of $\SLR$ are generally defined as having locally minimal arc length between any two of their (close enough) points.

By (5.9) the second order differential equation system of the $\SLR$ geodesic curve is the following:
\begin{equation}
\begin{gathered}
\ddot{r}=\sinh(2r)~\! \dot{\theta}~\! \dot{\phi}+\frac12 \big( \sinh(4r)-\sinh(2r) \big)\dot{\theta} ~\! \dot{\theta},\\
\ddot{\phi}=2\dot{r}\tanh{(r)}(2\sinh^2{(r)}~\! \dot{\theta}+ \dot{\phi}),\\ \ddot{\theta}=\frac{2\dot{r}}{\sinh{(2r)}}\big((3 \cosh{(2r)}-1)
\dot{\theta}+2\dot{\phi} \big). \tag{5.10}
\end{gathered}
\end{equation}
We can assume, by the homogeneity, that the starting point of a geodesic curve is the origin $(1,0,0,0)$.
Moreover, $r(0)=0,~ \phi(0)=0,~ \theta(0)=0,~ \dot{r}(0)=\cos(\alpha),~ \dot{\phi}(0)=\sin(\alpha)=-\dot{\theta}(0)$ are the initial values
in Table 1 for the solution of (5.10),
and so the unit velocity will be achieved.
\smallbreak
\centerline{\vbox{
\halign{\strut\vrule~\hfil $#$ \hfil~\vrule
&\quad \hfil $#$ \hfil~\vrule
&\quad \hfil $#$ \hfil\quad\vrule
&\quad \hfil $#$ \hfil\quad\vrule
&\quad \hfil $#$ \hfil\quad\vrule
\cr
\noalign{\vskip2pt}
\noalign{\hrule}
\multispan2{\strut\vrule\hfill{\bf  Table 2} \hfill\vrule}%
\cr
\noalign{\hrule}
\noalign{\vskip2pt}
{\rm Types} & {}  \cr
\noalign{\vskip2pt}
\noalign{\hrule}
\noalign{\vskip2pt}
\begin{gathered} 0 \le \alpha < \frac{\pi}{4} \\ (\bH^2-{\rm like~direction}) \end{gathered}
& \begin{gathered}  r(s,\alpha)={\mathrm{arsinh}} \Big( \frac{\cos{\alpha}}{\sqrt{\cos{2\alpha}}}\sinh(s\sqrt{\cos{2\alpha}}) \Big) \\
\theta(s,\alpha)=-{\mathrm{arctan}} \Big( \frac{\sin{\alpha}}{\sqrt{\cos{2\alpha}}}\tanh(s\sqrt{\cos{2\alpha}}) \Big) \\
\phi(s,\alpha)=2\sin{\alpha} s + \theta(s,\alpha) \end{gathered} \cr
\noalign{\vskip2pt}
\noalign{\hrule}
\noalign{\vskip2pt}
\begin{gathered} \alpha=\frac{\pi}{4} \\ ({\rm light~direction}) \end{gathered} &
\begin{gathered}  r(s,\alpha)={\mathrm{arsinh}} \Big( \frac{\sqrt{2}}{2} s \Big) \\
\theta(s,\alpha)=-{\mathrm{arctan}} \Big( \frac{\sqrt{2}}{2} s \Big) \\
\phi(s,\alpha)=\sqrt{2} s +\theta(s,\alpha) \end{gathered} \cr
\noalign{\vskip1pt}
\noalign{\hrule}
\noalign{\vskip2pt}
\begin{gathered} \frac{\pi}{4}  < \alpha \le \frac{\pi}{2} \\ ({\rm fibre-like~direction}) \end{gathered} &
\begin{gathered}  r(s,\alpha)={\mathrm{arsinh}} \Big( \frac{\cos{\alpha}}{\sqrt{-\cos{2\alpha}}}\sin(s\sqrt{-\cos{2\alpha}}) \Big) \\
\theta(s,\alpha)=-{\mathrm{arctan}} \Big( \frac{\sin{\alpha}}{\sqrt{-\cos{2\alpha}}}\tan(s\sqrt{-\cos{2\alpha}}) \Big) \\
\phi(s,\alpha)=2\sin{\alpha} s + \theta(s,\alpha) \end{gathered}  \cr
\noalign{\vskip2pt}
\noalign{\hrule}
\noalign{\hrule}}}}
\smallbreak
The equation of the geodesic curve in the hyperboloid model was determined in \cite{DESS09}, with the usual geographical sphere coordiantes
$(\lambda, \alpha)$, $(-\pi < \lambda \le \pi, ~ -\frac{\pi}{2}\le \alpha \le \frac{\pi}{2})$,
and the arc length parameter $0 \le s \in \bR$. 
\begin{definition}
The {\rm distance} $d(P_1,P_2)$ between the points $P_1$ and $P_2$ is defined by the arc length of the geodesic curve
from $P_1$ to $P_2$.
\end{definition}
\begin{definition}
The {\rm geodesic sphere} of radius $\rho$ (denoted by $S_{P_1}(\rho)$) with center at point $P_1$ is defined as the set of all points
$P_2$ satisfying the condition $d(P_1,P_2)=\rho$. Moreover, we require that the geodesic sphere is a simply connected
surface without selfintersection (Fig.~16).
\end{definition}
\begin{definition}
The body of the geodesic sphere of centre $P_1$ and with radius $\rho$ is called {\rm geodesic ball}, denoted by $B_{P_1}(\rho)$,
i.e., $Q \in B_{P_1}(\rho)$ iff $0 \leq d(P_1,Q) \leq \rho$.
\end{definition}
\begin{figure}[ht]
\centering
\includegraphics[width=12cm]{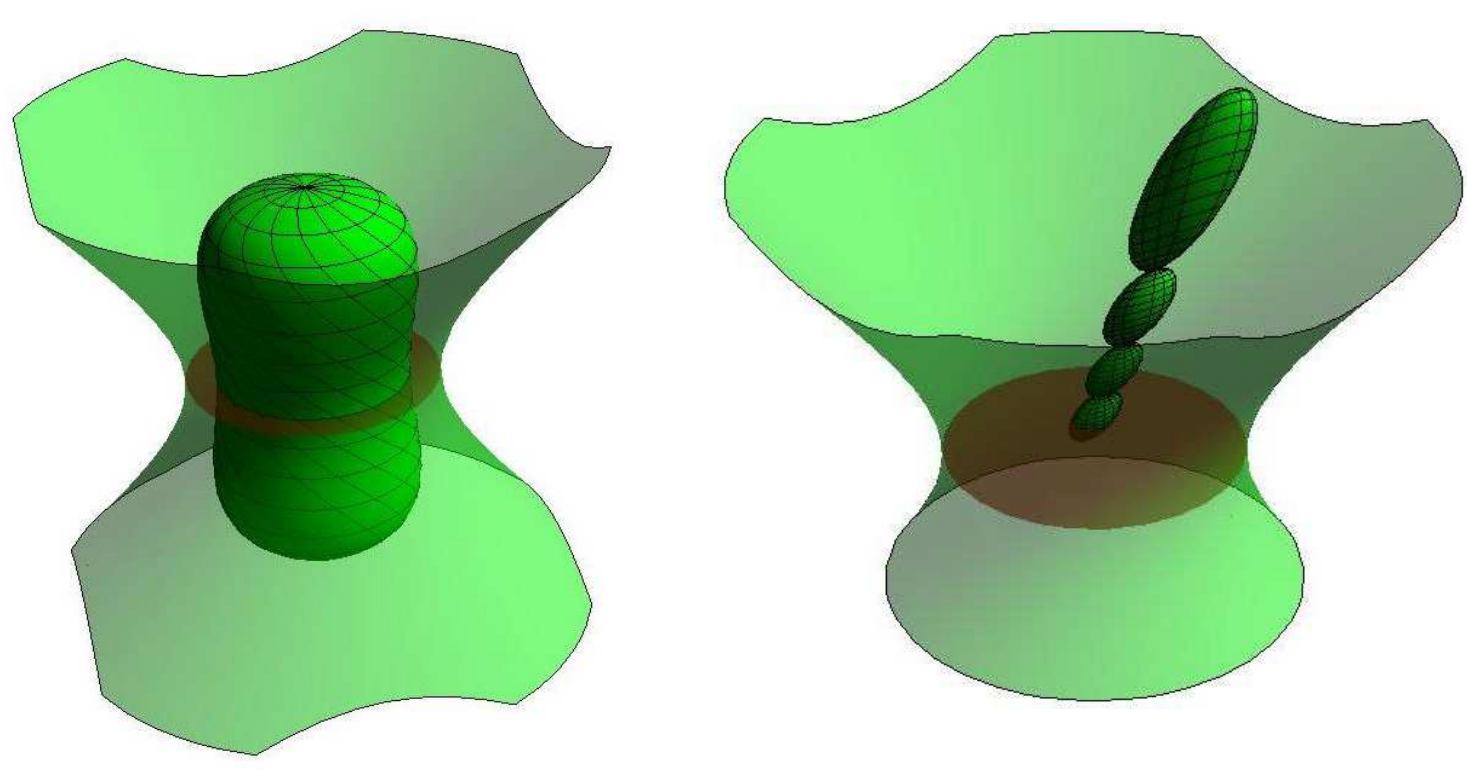}
\caption{a.~Geodesic sphere of radius 1 centred at the origin,~b. A series of tangent geodesic 
spheres of radius $\frac{1}{6}$ centred along a fibre line.}
\label{}
\end{figure}
 From (5.11) it follows that $S(\rho)$ is a simply connected surface in both $\mathbf{E}^3$ and $\SLR$
 if $\rho \in [0,\frac{\pi}{2})$.
 If $\rho\ge \frac{\pi}{2}$ then the universal cover should be discussed.
 {\it Therefore, we consider geodesic spheres and balls only with radii $\rho \in [0,\frac{\pi}{2})$ in the following}.
 These will be satisfactory in our cases.
\subsection{Geodesic triangles and their interior angle sums}
In this subsection we recall the results of \cite{CsSz16}.
\subsubsection{Fibre-like right angled triangles}
A geodesic triangle is called fibre-like if one of its edges lies on a fibre line (Fig.~17). 
We can assume without loss of generality
that the vertices $A_1$, $A_2$, $A_3$ of a fibre-like right angled triangle $T_g$ has the following coordinates:
$A_1=(1,0,0,0),~A_2=(1,0,y^2,0),~A_3=(1,x^3,0,0)$.

The geodesic segment $A_1A_2$ lies on the $y$ axis, the geodesic segment 
$A_1A_3$ lies on the $x$ axis and its angle is
$\omega_1=\frac{\pi}{2}$ in $\SLR$ (this angle is also $\frac{\pi}{2}$ in Euclidean sense since $A_1=E_0$).  

In order to determine the further interior angles of fibre-like geodesic triangle $A_1A_2A_3$ 
we define \emph{translations} $\bT_{A_i}$, $(i\in \{2,3\})$ as elements of the isometry group of $\SLR$, that 
maps the origin $E_0$ onto $A_i$. 
E.g. the isometry $\bT_{A_2}$ and its inverse (up to a positive determinant factor) and the images 
$\bT_{A_2}(A_i)$ of the vertices $A_i$ $(i \in \{1,2,3\}$ are the following:   
\begin{equation}
\begin{gathered}
\bT^{-1}_{A_2}(A_1)=A_1^2=(1,0,-y^2,0),~\bT^{-1}_{A_2}=A_2^2(A_2)=E_0=(1,0,0,0), 
\\ \bT^{-1}_{A_2}(A_3)=A_3^2=(1,x^3,-y^2,x^3y^2). 
\end{gathered}
\tag{5.11}
\end{equation}
Similarly to the above cases we obtain:
\begin{equation}
\begin{gathered}
\bT^{-1}_{A_3}(A_1)=A_1^3=(1,-x^3,0,0),~\bT^{-1}_{A_3}(A_2)=A_2^3=(1,-x^3,y^2,-x^3y^2), \\ \bT^{-1}_{A_3}(A_3)=A_3^3=E_0=(1,0,0,0), 
\end{gathered}
\tag{5.12}
\end{equation}
\begin{figure}[ht]
\centering
\includegraphics[width=11cm]{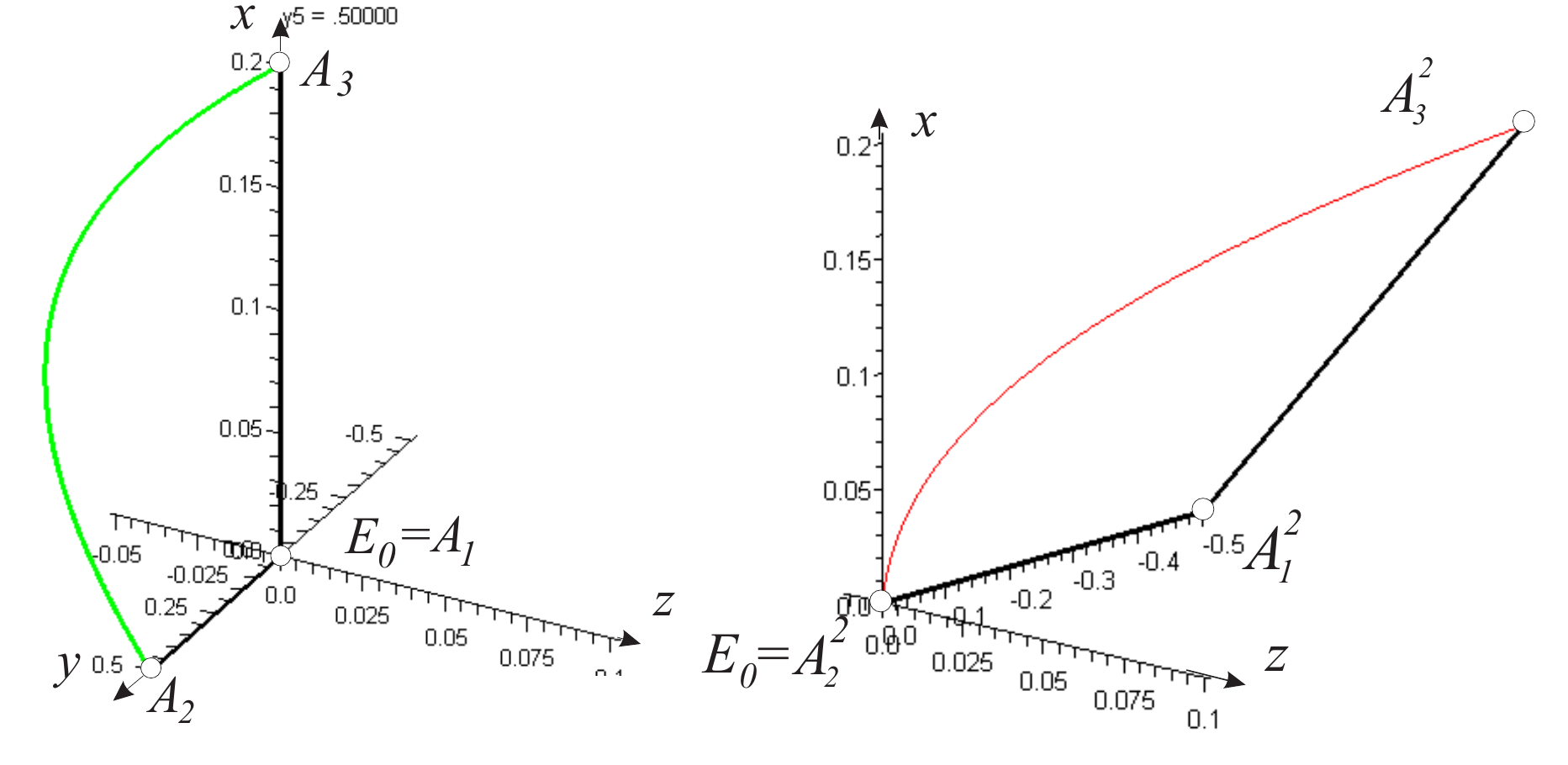}
\caption{Fibre-like geodesic triangle $A_1A_2A_3$ and its translated image $A_1^2A_2^2A_3^3$ by translation $\bT_{A_2}$}
\label{transant}
\end{figure}
We denote the oriented unit tangent vector of the oriented geodesic curves $g(E_0, A_i^j)$ with $\mathbf{t}_i^j$ where
$(i,j)=(2,3),(3,2),(1,2),(1,3)$ (Fig.~17).
The Euclidean coordinates of $\mathbf{t}_i^j$ are: 
\begin{equation}
\mathbf{t}_i^j=(\sin(\alpha_i^j), \cos(\alpha_i^j) \cos(\lambda_i^j), \cos(\alpha_i^j) \sin(\lambda_i^j)). \notag
\end{equation}
\begin{lemma}[\cite{CsSz16}]
The sum of the interior angles of a fibre-like right angled geodesic triangle is greater or equal to $\pi$. 
\end{lemma}
\subsubsection{Hyperbolic-like right angled geodesic triangles}
\begin{figure}[ht]
\centering
\includegraphics[width=5cm]{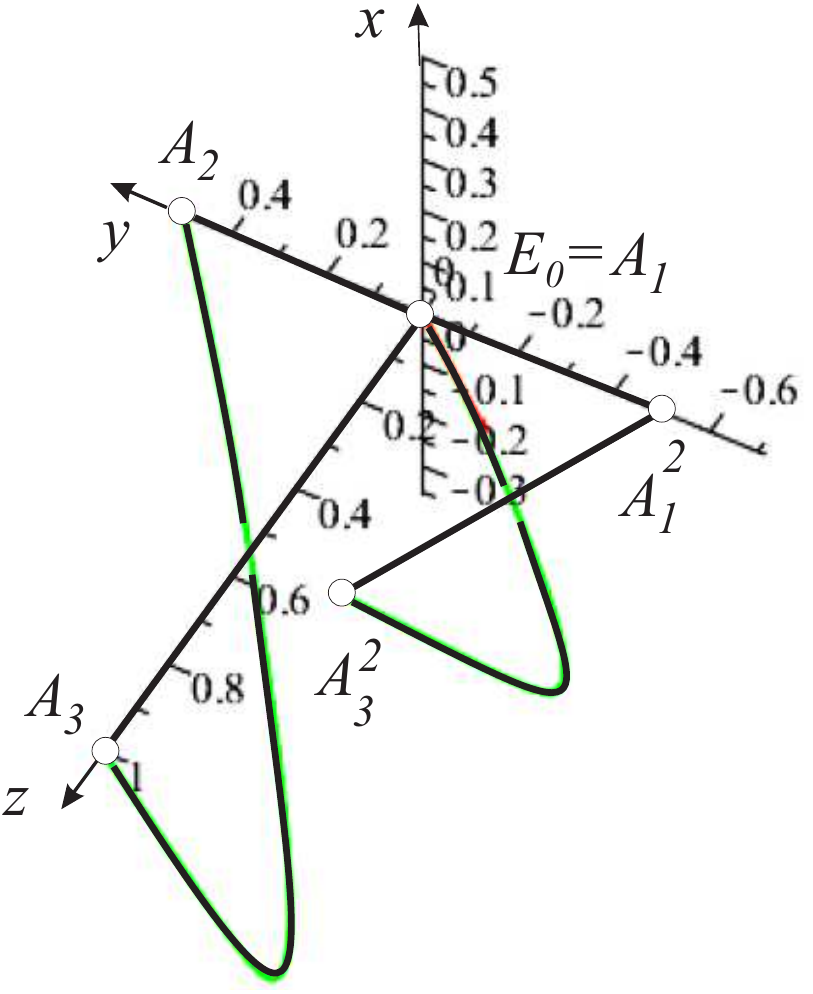}
\caption{Hyperbolic-like geodesic triangle $A_1A_2A_3$ and its translated copy $A_1^2A_2^2A_3^2$}
\label{hyph}
\end{figure}

A geodesic triangle is hyperbolic-like if its vertices lie in the base plane (i.e. $[y,z]$ coordinate plane) of the model. 
In this section we analyze the interior angle sum of the right angled hyperbolic-like
triangles. We can assume without loss of generality
that the vertices $A_1$, $A_2$, $A_3$ of a hyperbolic-like right angled triangle (see Figure \ref{hyph}) $T_g$ have the following coordinates:
$A_1=(1,0,0,0),~A_2=(1,0,y^2,0),~A_3=(1,0,0,z^3)$.
The geodesic segment $A_1A_2$ lies on the $y$ axis, the geodesic segment $A_1A_3$ lies on the $z$ axis and its angle is
$\omega_1=\frac{\pi}{2}$ in the $\SLR$ space (this angle is in Euclidean sense also $\frac{\pi}{2}$ since $A_1=E_0$).  

Similarly we get to the above cases that the images $\bT^{-1}_{A_j}(A_i)$ of the vertices $A_i$ $(i \in \{1,2,3\}, ~ j \in \{2,3\} )$  
are the following (see also Figure \ref{hyph}):   
\begin{equation}
\begin{gathered}
\bT^{-1}_{A_2}(A_1)=A_1^2=(1,0,-y^2,0),~\bT^{-1}_{A_2}=A_2^2(A_2)=E_0=(1,0,0,0), \\ \bT^{-1}_{A_2}(A_3)=A_3^2=(1,y^2z^3,-y^2,z^3), 
\end{gathered} \tag{5.13}
\end{equation}
\begin{equation}
\begin{gathered}
\bT^{-1}_{A_3}(A_1)=A_1^3=(1,-z^3,0,0),~\bT^{-1}_{A_3}(A_2)=A_2^3=(1,-y^2z^3,y^2,-z^3), \\ \bT^{-1}_{A_3}(A_3)=A_3^3=E_0=(1,0,0,0), 
\end{gathered} \tag{5.14}
\end{equation}
Finally, we obtain the following lemma:
\begin{lemma}[\cite{CsSz16}]
The interior angle sums of hyperbolic-like geodesic triangles can be less or equal to $\pi$.
\end{lemma}
\begin{conj}[\cite{CsSz16}]
The interior angle sum of any hyperbolic-like right angled geodesic triangle is less or equal to $\pi$.
\end{conj}
\subsubsection{Geodesic triangles with interior angle sum $\pi$}
In the above sections we discussed the fibre- and hyperbolic-like geodesic triangles and proved that there are right angled
geodesic triangles whose angle sum $\sum_{i=1}^{3}(\omega_i)$ is greater than or equal to $\pi$, less than or 
equal to $\pi$, but $\sum_{i=1}^{3}(\omega_i)=\pi$ is
realized only if one of the vertices of a geodesic triangle $A_1A_2A_3$ tends to the infinity (see Table 3-4). We prove the following 
\begin{lemma}[\cite{CsSz16}]
There is geodesic triangle $A_1A_2A_3$ with interior angle sum $\pi$ where its vertices are {\it proper} (i.e. $A_i \in \SLR$).   
\end{lemma}
We summarize the Lemmas of this Section as follows
\begin{Theorem}[\cite{CsSz16}]
The sum of the interior angles of a geodesic triangle of $\SLR$ space can be greater, less or equal to $\pi$. 
\end{Theorem}
\section{$\SOL$ space}
The projective model of Sol geometry was developed in \cite{M97}. 

\subsection{Basic notions of $\SOL$ geometry}
\label{sec:1}
In this Section we summarize the significant concepts and notations of real $\SOL$ geometry (see \cite{M97}, \cite{S}).

$\SOL$ is defined as a 3-dimensional real Lie group with multiplication  
\begin{equation}
     \begin{gathered}
(a,b,c)(x,y,z)=(x + a e^{-z},y + b e^z ,z + c).
     \end{gathered} \tag{6.1}
     \end{equation}
We note that conjugation by $(x,y,z)$ leaves invariant the plane $(a,b,c)$ with fixed $c$:
\begin{equation}
     \begin{gathered}
(x,y,z)^{-1}(a,b,c)(x,y,z)=(x(1-e^{-c})+a e^{-z},y(1-e^c)+b e^z ,c).
     \end{gathered} \tag{6.2}
     \end{equation}
Moreover, for $c=0$ the action of $(x,y,z)$ is only by its $z$-component, where $(x,y,z)^{-1}=(-x e^{z}, -y e^{-z} ,-z)$. Thus the $(a,b,0)$ plane is distinguished as a {\it base plane} in
$\SOL$, or by other words, $(x,y,0)$ is a normal subgroup of $\SOL$.
$\SOL$ multiplication can also be affinely (projectively) interpreted by ``right translations" 
on its points as the following matrix formula shows, according to (6.1):
     \begin{equation}
     \begin{gathered}
     (1;a,b,c) \to (1;a,b,c)
     \begin{pmatrix}
         1&x&y&z \\
         0&e^{-z}&0&0 \\
         0&0&e^z&0 \\
         0&0&0&1 \\
       \end{pmatrix}
       =(1;x + a e^{-z},y + b e^z ,z + c)
       \end{gathered} \tag{6.3}
     \end{equation}
by row-column multiplication.
This defines ``translations" $\mathbf{L}(\mathbf{R})= \{(x,y,z): x,~y,~z\in \mathbf{R} \}$ 
on the points of space $\SOL= \{(a,b,c):a,~b,~c \in \mathbf{R}\}$. 
These translations are not commutative, in general. 
Here we can consider $\mathbf{L}$ as projective collineation group with right actions in homogeneous 
coordinates as usual in classical affine-projective geometry.
We will use the Cartesian homogeneous coordinate simplex $E_0(\be_0)$,$E_1^{\infty}(\be_1)$, \ $E_2^{\infty}(\be_2)$, \ 
$E_3^{\infty}(\be_3), \ (\{\be_i\}\subset \bV^4$ \ $\text{with the unit point}$ $E(\be = \be_0 + \be_1 + \be_2 + \be_3 ))$ 
which is distinguished by an origin $E_0$ and by the ideal points of coordinate axes, respectively. 
Thus {$\SOL$} can be visualized in the affine 3-space $\bA^3$
(so in Euclidean space $\bE^3$) as well (see \cite{Sz19}).

In this affine-projective context E. Moln\'ar in \cite{M97} derived the usual infinitesimal arc length square at any point 
of $\SOL$, by pull back translation, as follows
\begin{equation}
   \begin{gathered}
      (ds)^2:=e^{2z}(dx)^2 +e^{-2z}(dy)^2 +(dz)^2.
       \end{gathered} \tag{6.4}
     \end{equation}
Hence we get the infinitesimal Riemann metric invariant under translations, by the symmetric metric tensor field $g$ on $\SOL$ by components as usual.

It will be important for us that the full isometry group $Isom(\SOL)$ has eight components, since the stabilizer of the origin 
is isomorphic to the dihedral group $\mathbf{D_4}$, generated by two involutive (involutory) transformations, preserving (6.4): 
\begin{equation}
   \begin{gathered}
      (1)  \ \ y \leftrightarrow -y; \ \ (2)  \ x \leftrightarrow y; \ \ z \leftrightarrow -z; \ \ \text{i.e. first by $3\times 3$ matrices}:\\      
     (1) \ \begin{pmatrix}
               1&0&0 \\
               0&-1&0 \\
               0&0&1 \\
     \end{pmatrix}; \ \ \ 
     (2) \ \begin{pmatrix}
               0&1&0 \\
               1&0&0 \\
               0&0&-1 \\
     \end{pmatrix}; \\
     \end{gathered} \tag{6.5}
     \end{equation}
     with its product, generating a cyclic group $\mathbf{C_4}$ of order 4
     \begin{equation}
     \begin{gathered}
     \begin{pmatrix}
                    0&1&0 \\
                    -1&0&0 \\
                    0&0&-1 \\
     \end{pmatrix};\ \ 
     \begin{pmatrix}
               -1&0&0 \\
               0&-1&0 \\
               0&0&1 \\
     \end{pmatrix}; \ \ 
     \begin{pmatrix}
               0&-1&0 \\
               1&0&0 \\
               0&0&-1 \\
     \end{pmatrix};\ \ 
     \mathbf{Id}=\begin{pmatrix}
               1&0&0 \\
               0&1&0 \\
               0&0&1 \\
     \end{pmatrix}. 
     \end{gathered} \notag
     \end{equation}
     Or we write by collineations fixing the origin $O(1,0,0,0)$:
\begin{equation}
(1) \ \begin{pmatrix}
         1&0&0&0 \\
         0&1&0&0 \\
         0&0&-1&0 \\
         0&0&0&1 \\
       \end{pmatrix}, \ \
(2) \ \begin{pmatrix}
         1&0&0&0 \\
         0&0&1&0 \\
         0&1&0&0 \\
         0&0&0&-1 \\
       \end{pmatrix} \ \ \text{of form (6.3)}. \tag{6.6}       
\end{equation}
A general isometry of $\SOL$ is defined by a product $\gamma_O \tau_X$, first $\gamma_O$ of 
form (6.6) then $\tau_X$ of (6.3). For
a general point $A(1,a,b,c)$, this will be a product $\tau_A^{-1} \gamma_O \tau_X$, mapping $A$ into $X(1,x,y,z)$. 

Conjugacy of translation $\tau$ by the above isometry $\gamma$, denoted as 
$\tau^{\gamma}=\gamma^{-1}\tau\gamma$.
\subsection{Related results in $\SOL$ geometry}
In paper \cite{BSz07} the authors considered the geodesics and Frenet formulas in $\SOL$ space, gave the differential equation system that describes geodesics but 
unfortunately the differential equation system in the main cases cannot be expressed with elementary functions. 

Therefore, the $\SOL$ geodesic sphere unfortunatelly cannot be expressed
in a closed explicite form. We could visualize the geodesic sphere by numerically solving differential equations (see Fig.~19). 
\begin{figure}[ht]
\centering
\includegraphics[width=13cm]{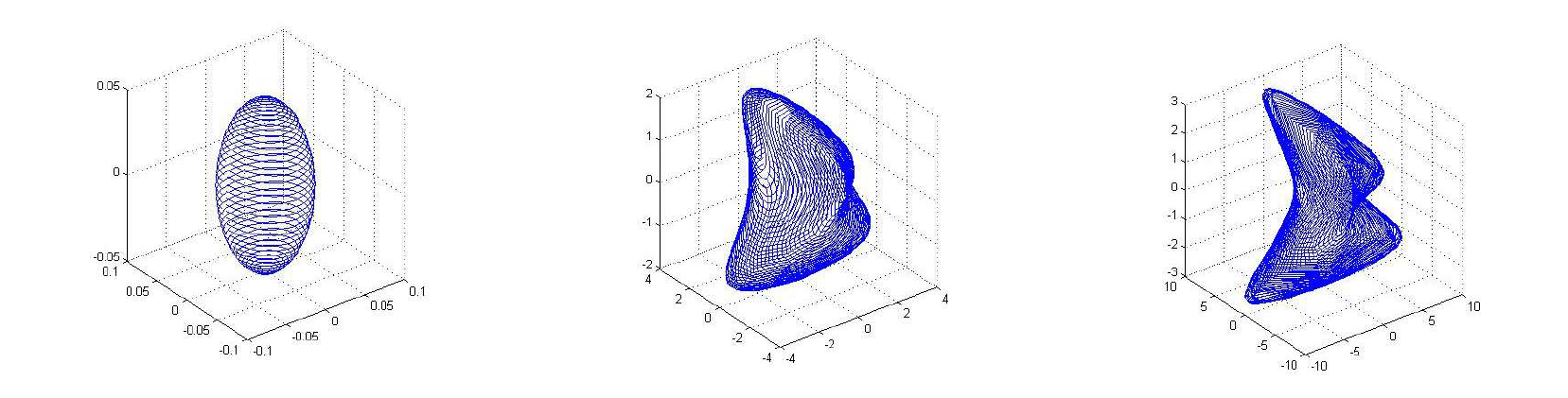}
\caption{$\SOL$ geodesic balls of radii $R=0.05$, $R=2$, $R=3$.}
\label{}
\end{figure}
In \cite{Er15} Z.~Erjavec studied certain class of Weingarten surfaces. The main result is that the only non-planar ruled Weingarten surface composed from vertical 
geodesics are surfaces $r(u,v)=(a e^{ku}, b e^{-ku},v)$. In \cite{Er18} the authors determined magnetic curves corresponding to the Killing magnetic fields in $\SOL$.
In \cite{Er20} Z. Erjavec and J. Inoguchi discussed magnetic curves with respect to the almost cosymplectic structure and investigated the curvature properties 
of these curves.

In \cite{MSz12} we classified $\SOL$ lattices in an algorithmic way with 17 types, 
in analogy of the 14 Bravais types of the Euclidean 3-lattices, but infinitely many $\SOL$ affine 
equivalence classes in each type. Then the discrete isometry groups of compact fundamental domain (crystallographic groups) can also be classified into infinitely 
many classes but finitely many types. To this we studied relations between 
$\SOL$ lattices and lattices of the pseudoeuclidean (or here rather called Minkowskian) plane. 
Moreover, we introduce the notion of $\SOL$ parallelepiped to every lattice type. 

In \cite{CMPSz} we studied a series of 2-generator Sol-manifolds depending on a
positive integer $n$ constructed them as tetrahedron manifolds, and proved that they are twofold coverings of the
3-sphere branched over specified links. 
\section{Ball packings in Thurston geometries}
\subsection{Geodesic ball packings in spaces of constant curvature}
Let $X$ denote a space of constant curvature, either the $n$-dimensional sphere $\mathbf{S}^n$,
Euclidean space $\mathbf{E}^n$, or hyperbolic space $\mathbf{H}^n$ with $n \geq 2$. An important question of discrete geometry
is to find the highest possible packing density in $X$ by congruent non-overlapping balls of a given radius \cite{Be}, \cite{G--K}.

Euclidean cases are the best explored. One major recent development has been the settling of the long-standing Kepler conjecture,
part of Hilbert's 18th problem, by Thomas Hales at the turn of the 21st century.
Hales' computer-assisted proof was largely based on a program set forth by L. Fejes T\'oth in the 1950's \cite{Ha}.

In $n$-dimensional hyperbolic geometry several new questions occur concerning the  packing and covering problems, e.g.
in $\HYN$ there are $3$ kinds of ``generalized balls (spheres)":  the {\it usual balls} (spheres),
{\it horoballs} (horospheres) and {\it hyperballs} (hyperspheres). Moreover, the definition of packing density is crucial in hyperbolic spaces
as shown by B\"or\"oczky \cite{B78}, for standard examples also see \cite{Be}, \cite{G--K}.
The most widely accepted notion of packing density considers the local densities of balls with respect to their Dirichlet--Voronoi cells
(cf. \cite{B78}). In order to consider ball packings in $\overline{\mathbf{H}}^n$, we use an extended notion of such local density.

In space $X^n$ let $d_n(r)$ be the density of $n+1$ mutually touching spheres or horospheres of radius
$r$ (in case of horosphere $r=\infty$) with respect
to the simplex spanned by their centres. L.~Fejes T\'oth and H.~S.~M.~Coxeter
conjectured that the packing density of balls of radius $r$ in $X^n$ cannot exceed $d_n(r)$.
This conjecture has been proved by C.~A.~Rogers for Euclidean space $\mathbf{E}^n$.
The 2-dimensional spherical case was settled by L.Fejes T\'oth in \cite{FTL}. 

In \cite{B78} K.~B\"or\"oczky proved the following theorem
for {\it ball and horoball} packings for any dimension $2 \le n \in \mathbf{N})$:

{\it In an $n$-dimensional space of constant curvature consider a packing of spheres of radius $r$.
In spherical space suppose that $r<\frac{\pi}{4}$.
Then the density of each sphere in its Dirichlet-Voronoi cell cannot
exceed the density of $n+1$ spheres of radius $r$ mutually
touching one another with respect to the simplex spanned by their centers.}

This density is $\approx 0.85328$ in $\mathbf{H}^3$ 
which is not realized by packings with equal balls. However, it is attained by the horoball packing
(case $r=\infty$) of
$\overline{\mathbf{H}}^3$ where the ideal centers of horoballs lie on the
absolute figure of $\overline{\mathbf{H}}^3$. This corresponds to packing an ideal regular
tetrahedron tiling is given by Coxeter-Schl\"afli symbol $\{3,3,6\}$.
Ball packings of hyperbolic $n$-space and of other Thurston geometries
are extensively discussed in the literature see e.g. \cite{Be, G--K, G--K--K, MSz18}, where the reader finds further references as well.

In this survey, we do not deal in detail with the examination of the ball 
(sphere) packings and coverings of spaces of constant curvature, 
so now we only mention that the questions regarding horosphere and 
hypersphere packings and coverings are not yet settled. 
New interesting problems have also arisen, the examination of which are related to the use of Busemann functions. 
The interested reader can read about the results of these in the papers 
\cite{ESz21, KSz, KSz14, KSz18, KSz22, Sz14, Sz17, Sz17-1, Sz17-2,
Sz07-1, Sz12, Sz12-2, Sz06-1, Sz06-2, Sz13-3, Sz13-4} and the references therein. 
\subsection{Geodesic ball packings in Thurston geometries of non-constant curvature}
Definitions of ball (sphere) packing and covering densities are already critical 
in hyperbolic geometry, 
therefore in order to introduce this concept to Thurston geometries of non-constant curvature 
we use the discrete isometry groups of the considered geometry. First, we have summarized the basic definitions and notions of this (see \cite{Sz13-1}).

Let $X$ be one of the five Thurston geometries with non-constant curvature
$$
\SXR,~\HXR,~\SLR,~\NIL,~\SOL,
$$
where the geodesic curves are generally defined as having locally minimal arc length between any two of their points
(sufficiently close to each other.
The system of equation for the parametrized geodesic curves $\gamma(\tau)$ in 
our model can be determined by the
general theory of Riemann geometry. 
{\it Then geodesic sphere and ball can usually be defined as follows.
We consider only geodesic ball packings which are transitively generated by discrete groups of isometries of $X$ and
the density of the packing is related to its Dirichlet-Voronoi cells.}

In the following let $\Gamma$ be a fixed group of isometries of $X$. Denote by $d(P_1,P_2)$ the distance of two points $P_1$, $P_2$.
\begin{definition}
We say that the point set
$$
\cD(K)=\{P\in X\,:\,d(K,P)\leq d(K^\bg,P)\text{ for all }\bg\in\ \Gamma\}
$$
is the {\rm Dirichlet--Voronoi cell ($D-V$~cell)} of $\Gamma$ around the kernel
point $K\in X$.
\end{definition}
\begin{definition}
We say that
$$
\Gamma_P=\{\bg\in\Gamma\,:\,P^\bg=P\}
$$
is the \emph{stabilizer subgroup} of $P \in X$ in $\Gamma$.
\end{definition}
\begin{definition}
Assume that the stabilizer $\Gamma_K=\bI$ is the identity, i.e., $\Gamma$ acts simply transitively on
the $\Gamma$-orbit of $K \in X$. Then let $B_K$ denote the \emph{largest ball}
with centre $K$ inside the D-V cell $\cD(K)$. Moreover, let $\rho(K)$ denote the
\emph{radius} of $B_K$. It is easy to see that
$$
\rho(K)=\min_{\bg\in\Gamma\setminus\bI}\frac12 d(K,K^\bg).
$$
\end{definition}
\begin{definition}
If the stabilizer $\Gamma_K > \bI$ then $\Gamma$ acts {\rm multiply transitively} on
the $\Gamma$-orbit of $K \in X$. Then the greatest ball radius of $\cB_K$ is
$$
\rho(K)=\min_{\bg\in\Gamma\setminus \Gamma_K}\frac12 d(K,K^\bg),
$$
where $K$ belongs to a 0-, 1-, or 2-dimensional region of $X$ (vertices, axes, reflection planes).
\end{definition}
In both cases the $\Gamma$-images of $B_K$ form a ball packing $\cB^\Gamma_K$ with centre
points $K^\bG$.
\begin{definition}
The \emph{density} of ball packing $\cB^\Gamma_K$ is
$$
\delta(K)=\frac{Vol(B_K)}{Vol(\cD(K))}.
$$
\end{definition}
It is clear that the orbit $K^\Gamma$ and the ball packing $\cB^\Gamma_K$ have the
same symmetry group. Moreover, this group contains the
crystallographic group $\Gamma$:
$$
Sym K^\Gamma=Sym\cB^\Gamma_K\geq\Gamma.
$$
\begin{definition}
We say that the orbit $K^\Gamma$ and the ball packing $\cB^\Gamma_K$ is
\emph{characteristic} if $Sym K^\Gamma=\Gamma$, otherwise the orbit is {\rm not
characteristic}.
\end{definition}
\subsubsection{Simply transitive ball packings}
 Let $\Gamma$ be a fixed group of isometries in the space $X$.
 \emph{Our goal is} to find a
 point $K\in\ X$ and the orbit $K^\Gamma$ for $\Gamma$ such that $\Gamma_K=\bI$
 and the density $\delta(K)$ of the corresponding ball packing
 $\cB^\Gamma(K)$ is maximal. In this case the ball packing $\cB^\Gamma(K)$ is
 said to be \emph{optimal.}
 
 Our aim is to determine the maximal radius $\rho(K)$ of the balls, and the maximal density $\delta(K)$.
 The space groups considered could have free parameters. So we have to find the densest ball packing for fixed
 parameters $p(\Gamma)$, and then we have to vary them to get the optimal ball packing
 \begin{equation}
 \delta(\Gamma)=\max_{K, \ p(\Gamma)}(\delta(K)). \tag{7.1}
 \end{equation}
 We look for the optimal kernel point
 in a 3-dimensional region, contained in a fundamental domain of $\Gamma$.
 
 \subsubsection{Multiply transitive ball packings}
 Similarly to the simply transitive case we must find a kernel
 point $K\in\ X$ and the orbit $K^\Gamma$ of $\Gamma$ such that
 the density $\delta(K)$ of the corresponding ball packing
 $\cB^\Gamma(K)$ is maximal, but here $\Gamma_K \ne \bI$. Such ball packing $\cB^\Gamma(K)$ is 
 also called \emph{optimal}.
 In this multiply transitive case we look for the optimal kernel point $K$
 in possible 0-, 1-, or 2-dimensional regions $\mathcal{L}$, respectively.
 Our aim is to deteremine the maximal radius $\rho(K)$ of the balls, and the maximal density $\delta(K)$.
 The space group considered can also have
 free parameters $p(\Gamma)$. Then we find the densest ball packing for fixed
 parameters, and vary them to find the optimal ball packing
 \begin{equation}
 \delta(\Gamma)=\max_{K\in \mathcal{L}, \ p(\Gamma)}(\delta(K)). \tag{7.2}
 \end{equation}
 \subsubsection{Geodesic ball packings in $\NIL$ space}
 
 {W. Heisenberg}'s famous real matrix group provides a non-commutative translation group of an affine 3-space.
 $\NIL$ geometry, which is one of the eight homogeneous
 Thurston 3-geometries, can be derived from this matrix group (see Section 4).
  
 In \cite{Sz07} we investigated the geodesic balls of $\NIL$ and computed their volume,
 introduced the notion of the $\NIL$ lattice, $\NIL$ parallelepiped and the density of the lattice-like ball packing.
 Moreover, we have determined the densest lattice-like geodesic ball packing by a family of $\NIL$ lattices.
 The density of this packing is $\approx 0.78085$, which may be surprising enough
 in comparison with the $3$-dimensional analogous Euclidean result $\frac{\pi}{\sqrt{18}} \approx 0.74048$. The kissing number of every ball
 in this packing is 14 (Fig.~20,~21). 
 {\it We conjecture that the densest geodesic ball packing belongs to the above ball arrangement in $\NIL$ space.}
 The symmetry group of this packing has also been described in \cite{MSz}.  
 \begin{figure}
 \centering
 \includegraphics[width=5cm]{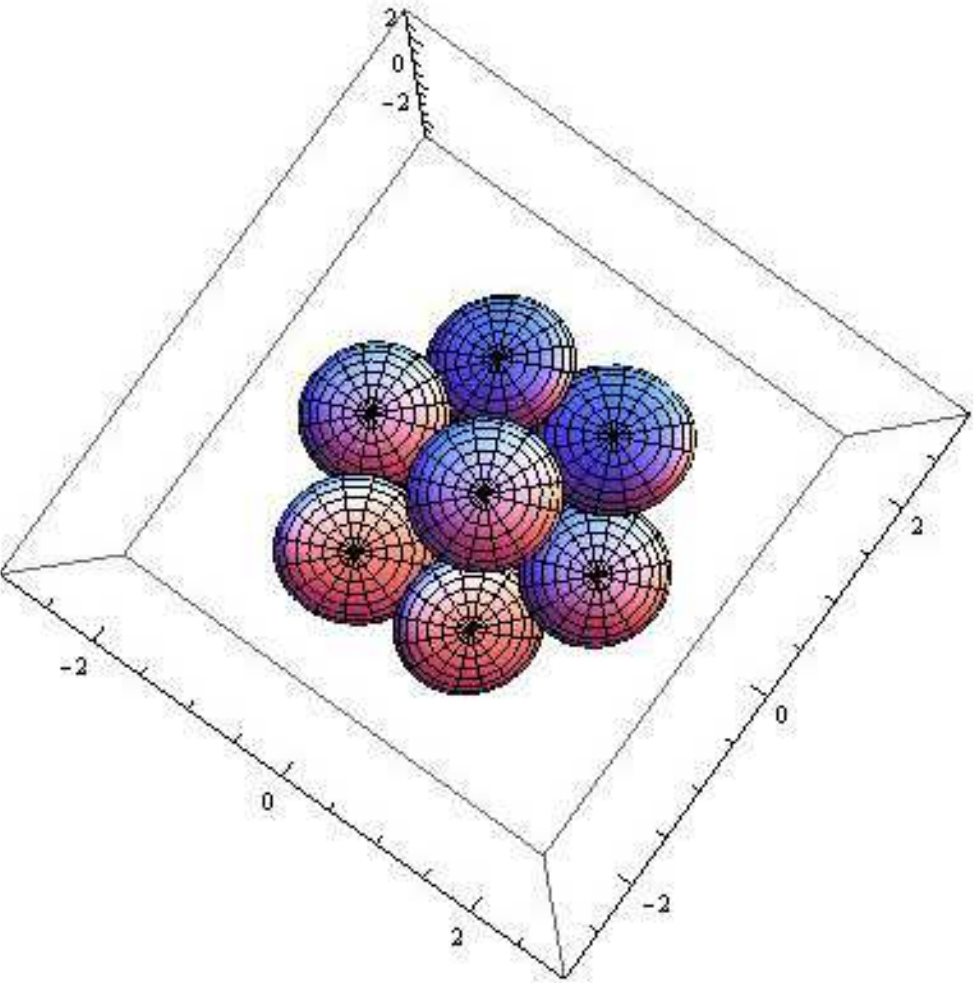},\includegraphics[width=4cm]{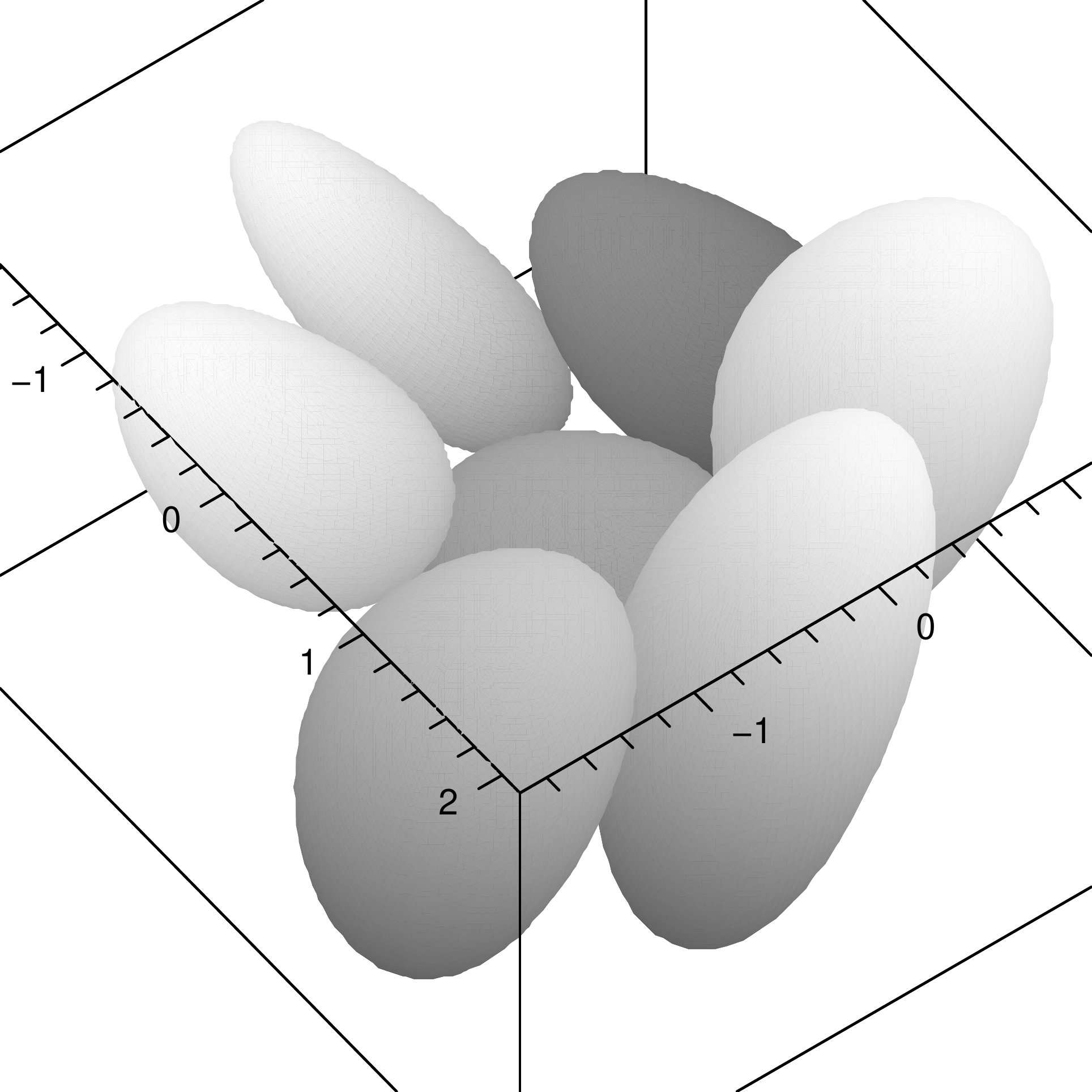}
 \caption{The densest geodesic lattice like geodesic ball packing in $\NIL$ space.}
 \label{}
 \end{figure}
 \begin{figure}
 \centering
 \includegraphics[width=12cm]{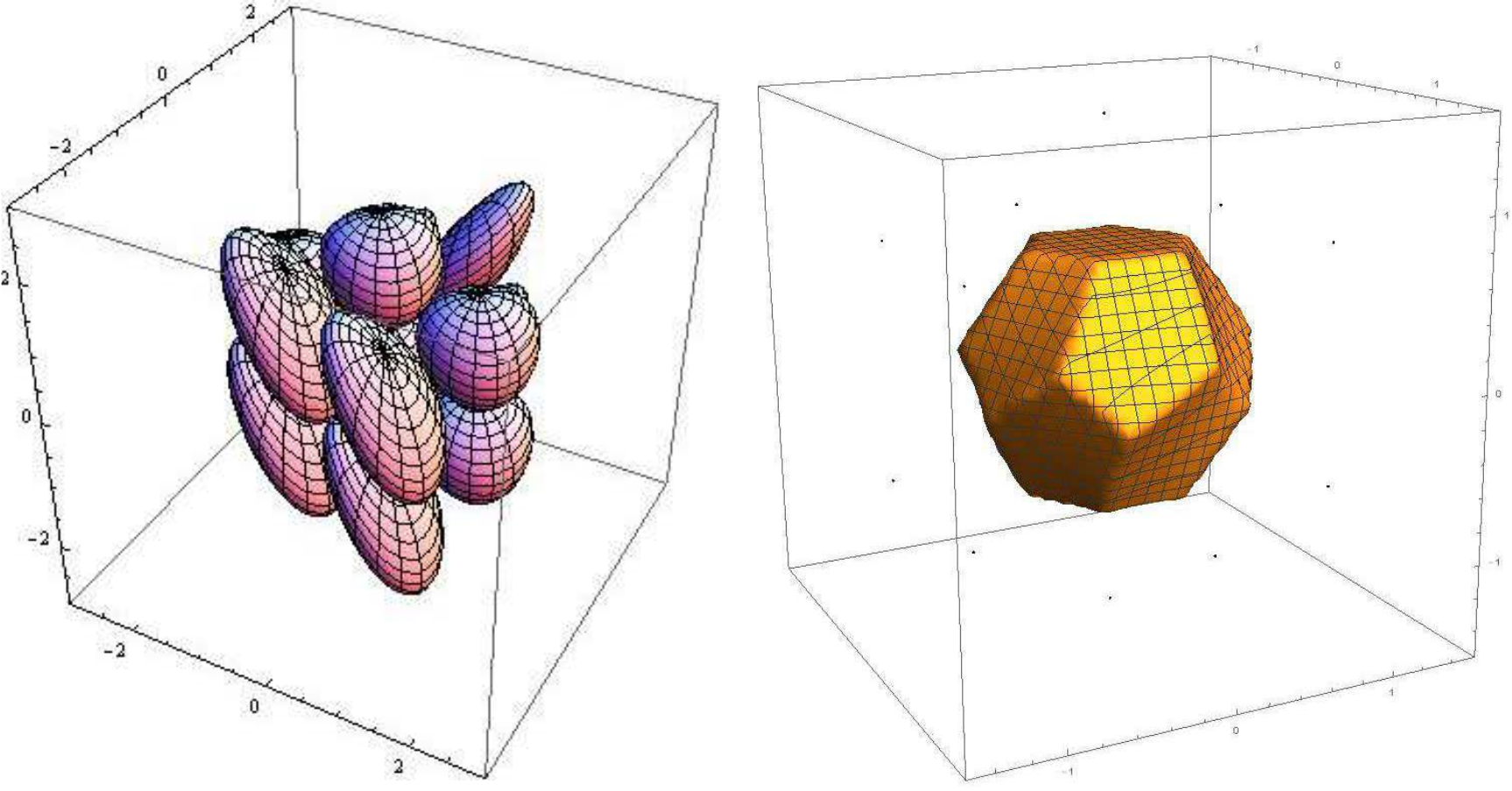}
 \caption{The densest geodesic lattice like geodesic ball packing in $\NIL$ space and the corresponding Dirichlet-Voronoi cell.}
 \label{}
 \end{figure}
In \cite{SchSz21} we investigated the geodesic ball packings related to $\NIL$ prism tilings where
we found that the largest density is $\approx 0.7272$ and the kissing number
of this ball arrangement is again $14$. 
 \subsubsection{Geodesic ball packings in $\HXR$ space}
 
 This Seifert fibre space is derived from the direct product of the hyperbolic plane $\bH^2$ and the real line $\bR$.
 In \cite{Sz12-5} we determined the geodesic balls of $\HXR$ space and computed their volume,
 defined the notion of the geodesic ball packing and its density.
 Moreover, we have developed a procedure to determine the density of the simply or multiply transitive geodesic ball packings for
 generalized Coxeter space groups of $\HXR$ and applied this algorithm to them.
 For the above space groups the Dirichlet-Voronoi cells are ``prisms" in $\HXR$ sense.
 The optimal packing density of the generalized Coxeter space groups is: $\approx 0.60726$. I am sure, that in this space there are
 denser ball packings.
\subsubsection{Geodesic ball (sphere) packings in $\SLR$ space}
In \cite{Sz14-3} we investigated the regular prisms and prism tilings in $\SLR$ 
(see Section 5) and in \cite{MSz14} we considered
the problem of geodesic ball packings related to tilings and their symmetry groups $\mathbf{pq2_1}$. 
Moreover, we computed the volumes of prisms and defined the
notion of the geodesic ball packing and its density. In \cite{MSz14} we developed a
procedure to determine the densities of the densest geodesic ball
packings for the tilings considered, more precisely, for
their generating groups $\mathbf{pq2_1}$ (for integer rotational
parameters $p,q$; $3\le p,~\frac{2p}{p-2} <q$). We looked for those
parameters $p$ and $q$ above, where the packing density as largest as possible. 
Currently our record is $0.5674$ for $(p, q) = (8,
10)$. In \cite{Sz16-1} we studied the non periodic geodesic ball packings related to the prism tilings and 
of the cases examined, the highest density that occurs is $\approx 0.6266$.
\subsection{Geodesic ball packings in $\SXR$ space}
The structure and the model of $\SXR$ geometry is described in Sections 2-3. 

In this section we briefly show the structure of the discrete isometry groups of $\SXR$ geometry, 
through which we can see that there are analogies with the Euclidean case, 
but even this geometry exhibits significant differences. 

The points in $\SXR$ geometry are described by $(P,p)$ where $P\in \bS^2$ and $p\in \bR$.
The isometry group $Isom(\SXR)$ of $\SXR$ can be derived by the direct product
of the isometry group of the spherical plane $Isom(\bS^2)$ and the isometry group of the real line $Isom(\bR)$.
The structure of an isometry group $\Gamma \subset Isom(\SXR)$  is the following: $\Gamma=\{(A_1 \times \rho_1), \dots (A_n \times \rho_n) \}$, where
$A_i \times \rho_i:=A_i \times (R_i,r_i):=(g_i,r_i)$, $(i \in \{ 1,2, \dots n \})$ and $A_i \in Isom(\bS^2)$, $R_i$ is either the identity map
$\mathbf{1_R}$ of $\bR$ or the point reflection $\overline{\mathbf{1}}_{\mathbf{R}}$. $g_i:=A_i \times R_i$ is called the linear part of the transformation
$(A_i \times \rho_i)$ and $r_i$ is its translation part.
The multiplication formula is the following:
\begin{equation}
(A_1 \times R_1,r_1) \circ (A_2 \times R_2,r_2)=((A_1A_2 \times R_1R_2,r_1R_2+r_2). \tag{8.1}
\end{equation}
{\it A group of isometries $\Gamma \subset Isom(\SXR)$ is called {\rm space group} if 
the linear parts form a finite group $\Gamma_0$ 
called the point group of
$\Gamma$. Moreover, the translation components of the identity of this point group are required to form a one dimensional lattice $L_{\Gamma}$ of $\bR$.}

In \cite{F01} J.~Z. {Farkas} classified and gave the complete list of the space groups in $\SXR$.

In \cite{Sz11-1} we studied the geodesic balls and their volumes in $\SXR$ space, 
moreover I have introduced the notion of geodesic ball packing
and its density and determined the densest simply and multiply transitive geodesic ball packings 
for generalized Coxeter space groups of $\SXR$, respectively.
The density of the densest packing is $\approx 0.82445$.

In paper \cite{Sz11-2} we studied the simply transitive locally optimal 
ball packings for the $\SXR$ space groups with Coxeter point groups such that
at least one of the generators is a non-trivial glide reflection.
We determined the densest simply transitive geodesic ball arrangements for the above space groups, 
moreover computed their optimal
densities and radii.
The density of the densest packing in this case is $\approx 0.80408$.

{\it In this survey we only recall the results from \cite{Sz13-1} where we studied the class 
of $\SXR$ space groups {\bf 4q.~I.~2} (with a natural parameter $q \ge 2$, see \cite{F01}).
Each of them belongs to the glide reflection groups, i.e.,
the generators $\bg_i \ (i=1,2,\dots m)$ of its point group $\Gamma_0$
are reflections and at least one of the possible translation components of the above 
generators differs from zero (see \cite{Sz11-2})}.

\subsubsection{A very dense multiply transitive ball packing in $\SXR$ geometry}
We considered a $\SXR$ space group (see \cite{F01, Sz11-1, Sz11-2}
with point group $\Gamma_0$ generated by three reflections $\bg_i ~ (i=1,2,3)$
\begin{equation}
\begin{gathered}
(+,~0,~[~~]~ \{(2,2,q)\}),~ q \ge 2, \ \\ \Gamma_0=(\bg_1,\bg_2,\bg_3 - \bg_1^2,\bg_2^2,\bg_3^2,(\bg_1\bg_3)^2, (\bg_2\bg_3)^2,
(\bg_1\bg_2)^q).
\end{gathered} \notag
\end{equation}
The possible translation parts $\tau_1$, $\tau_2$, $\tau_3$ of the corresponding generators of $\Gamma_0$ are derived from the so-called Frobenius congruence relations:
$$(\tau_1,\tau_2,\tau_3) \cong (0,0,0),~\big(0,0,\frac{1}{2}\big),~\big(\frac{1}{2},\frac{1}{2},\frac{1}{2}\big),~
\big(\frac{1}{2},\frac{1}{2},0 \big),~\big(0, \frac{1}{2}, 0 \big),~\big(0, \frac{1}{2}, \frac{1}{2} \big).$$
If $(\tau_1,\tau_2,\tau_3) \cong ~(0,0,\frac{1}{2})$ then we have obtained the $\SXR$ space group {\bf 4q.~I.~2} (for a fixed $q$, $2\le q \in \mathbf{N}$).

The fundamental domain of the point group of the space group considered is a spherical triangle $A_1A_2A_3$ with
angles $\frac{\pi}{q}$, $\frac{\pi}{2}$, $\frac{\pi}{2}$ in the base plane $\Pi$.
It can be assumed that the fibre coordinate of the center of the optimal ball is zero and it is a point of triangle $A_1A_2A_3$.

{\it We consider ball packings related to parameter $q=2$.}

To determine the optimal multiply transitive ball packing we studied 2-cases:
\begin{enumerate}
\item[1.] $K$ is an interior point of the spherical geodesic segment $A_2A_3$ (or $A_1A_3$).
In this situation the point $K$ and its images by $\Gamma={\bf 4q.~I.~2}$ $(q=2)$, as the centers of the optimal ball arrangement
$\mathcal{B}_{opt}(K,R)$ have to satisfy the following requirements, 
since an arbitrary ball of the optimal packing is fixed by
its neighbouring balls.
Here we got the optimal packing if $K=A_2$ (or $K=A_1$) with the following data:
\begin{equation}
\begin{gathered}
\phi_3= \frac{\pi}{2} \approx 1.57079633,\ \ \theta_3 =0, \ \ R_3 = \frac{\pi}{2}\approx 1.57079633, \\
Vol(B(R_3))\approx 13.74539472, \ \ \delta(R_3,K_3) \approx 0.69634983.
\end{gathered} \tag{8.2}
\end{equation}
\item[2.] $K=A_3$.
Fig.~22 shows the orbit of the point $K=A_3$ by the space group considered. The images of $K$ lie on a line through the origin and $A_3$.

\begin{equation}
\begin{gathered}
\phi_4= \frac{\pi}{4} \approx 0.78539816,\ \ \theta_4 =\frac{\pi}{2} \approx 1.57079633, \ \ R_4 \approx 1.81379936, \\
Vol(B(R_4))\approx 20.00238509, \ \ \delta(R_4,K_4) \approx 0.87757183.
\end{gathered} \tag{2.12}
\end{equation}
The "outwardly transformed" images of the balls surround the initial balls (see Fig.~22) thus 
the touching number of this packing is 4.
Finally, we obtain the following
\begin{Theorem}[\cite{Sz13-1}]
The ball arrangement $\mathcal{B}_{opt}(R_4,K_4)$ provides the densest multiply transitive ball 
packing of the $\SXR$ space group {\bf 4q.~I.~2} $(q=2)$.
\end{Theorem}
\end{enumerate}
\begin{figure}[ht]
\centering
\includegraphics[width=12cm]{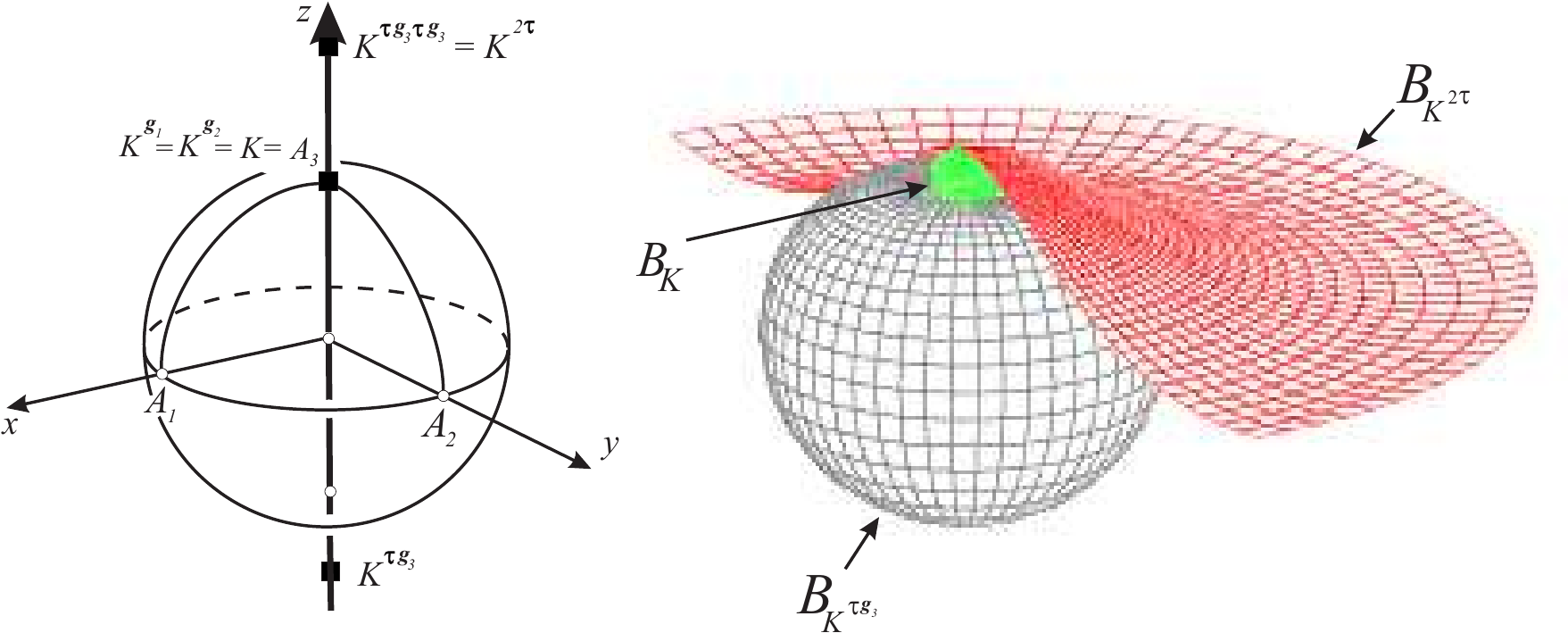}
\caption{a. The orbit of $K=A_3$ by the group $\Gamma={\bf 4q.~I.~2}$ $(q=2)$. 
b. The densest ball packing is determined by its balls $B_K$, $B_{K^{\tau g_3}}$ and a part of the sphere $B_{K^{2\tau}}$.}
\label{}
\end{figure}
\begin{rmrk}
\begin{enumerate}
\item[1.] To the author's best knowledge there are no results for the geodesic ball packings 
in $\SOL$  geometry at the time of writing.
\item[2.] In $\NIL$, $\SLR$ and $\SOL$ spaces we have studied the so-called 
{\it translation ball packings} 
\cite{MSzV, MSzV14, Sz12-1, Sz19, Sz13-2, VSz19}
but we did not consider these cases in this work.
\end{enumerate}
\end{rmrk}
\subsubsection{The conjecture for the densest ball arrangement in Thurston geometries}
We introduced the density function for the geodesic ball packings generated by a 
discrete group of isometries in a given Thurston geometry.
This density is related to the Dirichlet-Voronoi cells generated by centres of the balls.
For these ball packings we can formulate the following
\begin{conj}[\cite{Sz13-1}]
Let $\cB$ be an arbitrary congruent geodesic ball packing in a Thurston geometry $X$, where $\cB$ is generated by a discrete isometry group of $X$.
The above determined ball arrangement $\mathcal{B}_{opt}(R_4,K_4)$ with density $\delta(R_4,K_4) \approx 0.87757183$ provides
the densest congruent geodesic ball packing for the Thurston geometries.
\end{conj}
The general definition of the density of congruent geodesic ball packings for the 
Thurston geometries is not settled yet. However, by our investigation
for any ``good" definition of density following conjecture may be formulated.
\begin{conj}[\cite{Sz13-1}]
The densest congruent geodesic ball packing in the Thurston geometries is realized by
the above ball arrangement $\mathcal{B}_{opt}(R_4,K_4)$ with density $\delta(R_4,K_4) \approx 0.87757183$.
\end{conj}

\vspace{5mm}

{\it In this paper we mentioned only some classical theorems and problems related to 
Thurston spaces, but we hope that from these
the reader can appreciate that our projective method is suitable to study and solve 
similar problems that represent a huge class of open mathematical problems. 
Detailed studies are the objective of ongoing research.}
\medbreak

\end{document}